\documentclass[12pt, a4paper]{amsart}

\usepackage{amssymb}
\usepackage{amsfonts, amsthm, amsmath, amscd, flafter}
\usepackage[all]{xy}
\usepackage{graphicx}
\usepackage{color}
\usepackage[hmarginratio=1:1]{geometry}
\usepackage{multirow}
\usepackage{tabularx}
\usepackage{enumerate}

\allowdisplaybreaks
\usepackage{url}
\makeatletter
\def\url@leostyle{%
  \@ifundefined{selectfont}{\def\UrlFont{\sf}}{\def\UrlFont{\small\ttfamily}}}
\makeatother
\newtheorem{theorem}{Theorem}[section]
\newtheorem{lemma}[theorem]{Lemma}
\newtheorem{proposition}[theorem]{Proposition}

\newtheorem{thm}{Theorem}

\newenvironment{remark}[1][Remark]{\begin{trivlist}
\item[\hskip \labelsep {\bfseries #1}]}{\end{trivlist}}
\newenvironment{definition}[1][Definition]{\begin{trivlist}
\item[\hskip \labelsep {\bfseries #1}]}{\end{trivlist}}

\newcommand{\lens}[2]{L({\scriptstyle #1},{\scriptstyle #2})}

\newcommand{\seifuno}[3]{\big(#1,({\scriptstyle #2},{\scriptstyle #3})\big)}

\newcommand{\seifdue}[5]{\big(#1,({\scriptstyle #2},{\scriptstyle #3}),
                       ({\scriptstyle #4},{\scriptstyle #5})\big)}

\newcommand{\seiftre}[7]{\big(#1,({\scriptstyle #2},{\scriptstyle #3}),
                       ({\scriptstyle #4},{\scriptstyle #5}),
                       ({\scriptstyle #6},{\scriptstyle #7})\big)}

\newcommand{\bigu}[4]{\bigcup\nolimits_{{\tiny{\matr {#1} {#2} {#3} {#4}}}\phantom{\Big|}\!\!}}
\newcommand{\bigb}[4]{\big/_{{\tiny{\matr {#1} {#2} {#3} {#4}}}\phantom{\Big|}\!\!}}
\newcommand{\matr} [4] {\left(\begin{array}{@{}c@{\ }c@{}} #1 & #2 \\ #3 & #4 \\ \end{array} \right)}

\begin{document}

\title{Exceptional Slopes on Manifolds of Small Complexity} 
\author{Fionntan Wouter Murray Roukema}
\thanks{The author was supported as a member of the Italian FIRB project 
`Geometry and topology of low-dimensional manifolds' (RBFR10GHHH)} 
\maketitle 
\begin{abstract} 
It has been observed that most manifolds in the Callahan-Hildebrand-Weeks census of 
cusped hyperbolic $3$-manifolds are obtained by surgery on the minimally twisted 5-chain link. 
A full classification of the exceptional surgeries on the 5-chain link has recently been 
completed. In this article, we provide a complete classification of the sets of exceptional 
slopes and fillings for all cusped hyperbolic surgeries on the minimally twisted 5-chain link, 
thereby describing the sets of exceptional slopes and fillings for most hyperbolic manifolds 
of small complexity. The classification produces the description of exceptional fillings for 
many families of one and two cusped manifolds, and provides supporting evidence for some well-known 
conjectures. One such family that appears in the classification is an infinite family of 
1-cusped hyperbolic manifolds with four Seifert manifold fillings and a toroidal filling. 
\end{abstract}

\section{Introduction}\label{intro:sec}

The set of exceptional slopes on a boundary component of a hyperbolic manifold has generated a 
lot of interest in the literature. There are many restrictions on the set of exceptional slopes 
on a boundary component of a hyperbolic 3-manifold $M$ and its corresponding fillings. For 
example, no such $M$ has two distinct $S^3$ fillings \cite{GL} or more than ten exceptional slopes 
\cite{Lack:setsize}. However, it is still not known if there exists a hyperbolic knot exterior in 
$S^3$ with a reducible filling, or an $M$ with a pair of exceptional slopes $\beta$ and $\beta'$ 
corresponding to a lens space and toroidal space so that the \textit{distance} (minimal number 
of intersections) between $\beta$ and $\beta'$ is greater than three, or if there is a manifold not equal to the 
Figure-8 knot exterior with 10 exceptional slopes. Conjecturally, no examples exist, 
see \cite{Cabling}, \cite{Gor'}, \cite[Problem~1.77]{problems} respectively. 

The distance between two exceptional slopes on a boundary component is at most 8 \cite{Lack:setsize} 
(which is realised on both the figure-eight knot exterior and the figure-eight sister manifold) 
and it is known that only finitely many one cusped 3-manifolds have exceptional slopes at distance 
more than 5 \cite{Ago2}. It is conjectured that if an orientable 3-manifold has two exceptional slopes 
at distance greater than five then it is obtained by surgery on the Whitehead link \cite{Gor1}. 

Also of interest are manifolds with more than one exceptional reducible filling; it is known that the 
distance between the fillings is 1 \cite{GL3}, and examples are given in \cite{EW}, \cite{HM}, 
and \cite{GLi}. Eudave-Mu\~{n}oz and Wu's examples in \cite{EW} are the only known with more 
than one boundary component \cite{Gor}, and Hoffman and Matignon ask in \cite{HM} if reducible pairs 
must have at least one $\lens 21$, $\lens31$, $\lens41$ summand and whether any hyperbolic manifold 
has three reducible fillings.

%Hoffman and Matignon in their paper in Pacific J Math 209 (2003) ask a couple of questions: (i) in such examples is there always an L(2,1), L(3,1) or L(4,1) summand? and (ii) are there any examples with 3 reducible fillings? You might point out that these are true for manifolds coming from M_5. Also, Eudave-Munoz and Wu give an infinite family of manifolds with two torus boundary components having a pair of reducible fillings. But I guess there are no such examples coming from M_5.

%With the exception of small Seifert manifolds that fibre over the sphere with three exceptional fibres, a great deal is 
%known about the maximal distance between slopes corresponding to fillings in different classes of non-hyperbolic manifolds, 
%and about the manifolds with slopes achieving these maximal distances, see \cite{Gor1} for an overview. 

In this article, by classifying the sets of exceptional slopes and the corresponding fillings for all 
manifolds obtained by surgery on the minimally twisted 5-chain link (see the rightmost link in Figure 
\ref{sequence:fig}), we provide experimental evidence that supports the above conjectures of 
Gonz\'{a}lez-Acu\~{n}a, Short, and Gordon. 

\begin{thm}\label{flash}
If $M$ is a cusped hyperbolic manifold obtained by surgery on the minimally twisted 5-chain link and 
$\tau$ is a fixed boundary component of $M$ then: 
\begin{enumerate}[(i)]
\item If $M$ is the exterior of a knot in $S^3$ then $M$ does not have a reducible filling;
\item If $M$ has two exceptional slopes on $\tau$ at distance greater than 3 apart then they do not 
correspond to a lens space and a toroidal filling;
\item If $M$ has 10 exceptional slopes on $\tau$ then $M$ is the figure-8 knot exterior;
\item If $M$ is a manifold with exceptional slopes on $\tau$ at distance greater than 5 
then $M$ is obtained by surgery on the Whitehead link.
\item $M$ does not have more than one reducible filling. 
\end{enumerate}
\end{thm}

A full analysis of the exceptional fillings of surgeries on the minimally twisted 5-chain 
link is given to obtain Theorem \ref{flash}. This produces a classification of exceptional filling 
types for infinitely many 1-cusped and 2-cusped manifolds (see Tables \ref{tc}-\ref{th}). These 1-cusped 
and 2-cusped manifolds are distinct from the examples in \cite{Magic} which all have a cyclic filling 
and at least five exceptional slopes. Among other families, we highlight:
\begin{itemize}
\item An infinite family of hyperbolic knots in $S^3$ with consecutive integral toroidal, 
small Seifert manifold, toroidal surgeries;
\item An infinite family of 1-cusped manifolds with a reducible filling and a small Seifert manifold
filling at distance one apart; 
\item An infinite family of 1-cusped hyperbolic manifolds with four small Seifert manifold fillings 
and a toroidal filling; 
\item An infinite family of 2-cusped manifolds with four fillings on a fixed cusp containing an 
essential annulus.
\end{itemize}
These families are not contained in the classification given in \cite{Magic}, and I am 
unaware of them appearing elsewhere in the literature. The specific description of these 
families and their exceptional fillings can be found in Table \ref{families_table}. 

The classification of exceptional fillings in this article does not improve any of the lower bounds on 
maximal distances to small Seifert manifolds that fibre over the sphere with three exceptional fibres obtained 
in \cite{Magic}. 

\subsection{The minimally twisted 5-chain link}\label{mt5cl_sec}

\begin{figure}[h!]
\begin{center}
\includegraphics[width = 12 cm]{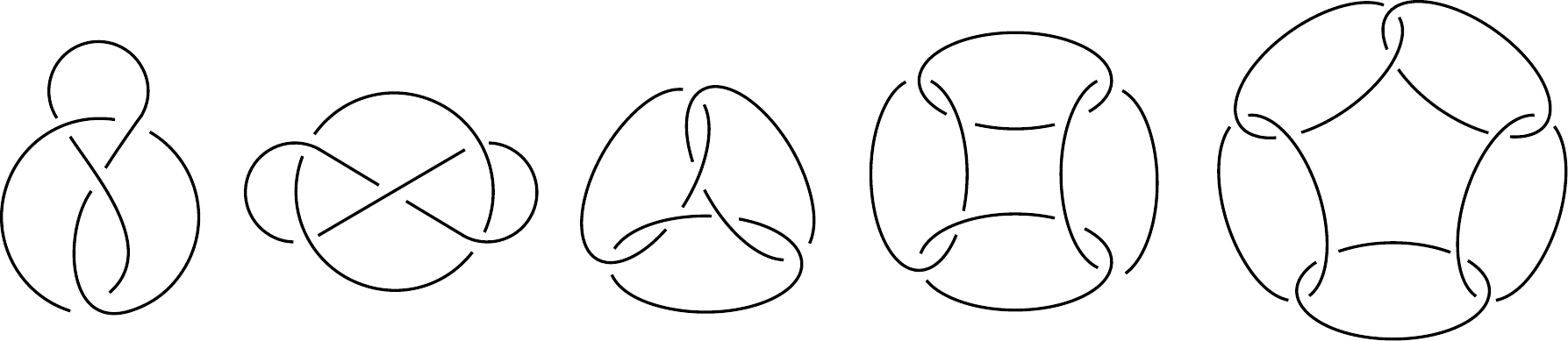}
\end{center}
\caption{The links 1CL, 2CL, 3CL, 4CL, 5CL in $S^3$ whose exteriors we denote 
by $M_1$, $M_2$, $M_3$, $M_4$, and $M_5$ respectively.}
\label{sequence:fig}
\end{figure}

A notable collection of hyperbolic chain links is described in \cite{MPR}; they are the figure-8 knot, 
the Whitehead link, the 3-chain link, the 4-chain link with a half twist, and the minimally 
twisted 5-chain link. These links are shown in Figure \ref{sequence:fig}. We follow \cite{MPR} 
and denote these chain links by 1CL, 2CL, 3CL, 4CL and 5CL, and their exteriors by $M_1$, $M_2$, 
$M_3$, $M_4$, $M_5$ respectively.  

The significance of this sequence of links comes from the following facts: each $M_i$ 
is the (or conjecturally the) smallest volume hyperbolic 3-manifold with $i$ cusps, see 
\cite{Ago} and \cite{Yoshida}; 
more than $80\%$ of the cusped hyperbolic 3-manifolds from the Callahan-Hildebrand-Weeks 
census \cite{CaHiWe} are surgeries on 5CL (personal communication with Nathan Dunfield). Furthermore, 5CL 
relates to the program of enumerating exceptional pairs at maximal distance. In particular, all 
knots realising half-integral toroidal surgeries \cite{essential:tori}, many of the knots realising 
lens space surgeries \cite{Baker}, and all cusped hyperbolic 3-manifolds with distinct reducible and 
toroidal fillings at maximal distance are obtained by surgery on 5CL \cite{Kang}.

It is easy to see that if $\partial M_n$ is equipped with the usual (meridian, longitude) 
homology basis then a $-1$ filling on any boundary component of $M_n$ results in $M_{n-1}$. 
As a result, any manifold obtained by surgery on $(n-1)$CL is obtained by surgery on $n$CL. 
A classification of the exceptional surgeries on 3CL is given in \cite{Magic} together with a 
complete description of the set of exceptional slopes, and corresponding exceptional 
fillings, on the boundary components of all hyperbolic manifolds obtained by surgery 
on 3CL. The statements of Theorem \ref{flash} are known to hold 
for any manifold obtained by surgery on 3CL (see the appendix of \cite{Magic}). 

A classification of exceptional surgeries on 5CL is found in \cite{MPR}. In this 
article we complete the description of the set of exceptional slopes and corresponding exceptional 
fillings for manifolds obtained by surgery on 5CL in Theorems \ref{5CL_thm} and \ref{4CL_thm}. 
We then use this classification to verify the statements of Theorem \ref{flash}. 

\subsection{Article structure} We start with Section \ref{basics} where we recall the classification from 
\cite{MPR}. In order to do so, we begin by recalling and introducing notation and terminology in 
Sections \ref{introbasics}-\ref{sicl}. In Section \ref{basics} we also establish some results which 
turn out to be of great use in the remainder of the paper (see Proposition \ref{fillingF} and Lemma 
\ref{sym} in Sections \ref{excep:fill:sec} and \ref{result:recap}). 

Theorems \ref{5CL_thm} and \ref{4CL_thm} are stated and proved in Section \ref{mr}. These theorems 
complete the classification of exceptional sets of slopes on cusped hyperbolic manifolds obtained 
by surgery on 5CL. The proofs are heavily reliant on Proposition \ref{fillingF} and Lemma \ref{sym}. 

Theorem \ref{flash} is proved in Section \ref{flash:proof}. The proof uses Theorems \ref{5CL_thm} 
and \ref{4CL_thm} to impose restrictions on the filling instructions that can correspond to a 
counterexample. The exceptional slopes and fillings of many families of manifolds are completely 
enumerated using Theorems \ref{5CL_thm} and \ref{4CL_thm} in Tables \ref{tc}--\ref{th} found in Section 
\ref{tables_sec}. Careful consideration of these tables is needed to complete the proof of Theorem 
\ref{flash}.

\subsection{Acknowledgements and remarks} This is the second version of this article; the initial version 
of this article omitted many details and was not clear as a result. The feedback from the anonymous referee 
on the first submission highlighted this. The current presentation has greatly benefited from the anonymous 
referee's remarks and from discussions with Marc Lackenby. The results of this article were mainly obtained 
as a graduate student at the University of Pisa under the supervision of Carlo Petronio and Bruno Martelli. 
The article has also benefited from discussions with Daniel Matignon, and from email correspondences with 
Carlo Petronio, Bruno Martelli, Cameron Gordon, and Nathan Dunfield. 

Some of the main results of this revised version have been extended and amended from the original version: 
Tables \ref{tc}-\ref{th} were not displayed in the original version, and the families of 
cusped manifolds highlighted in the introduction were not mentioned. Moreover, several typos/omitted examples 
in Tables \ref{tb1}-\ref{tb6} have been corrected.

\section{Background Terminology and Useful Results}\label{basics}

\subsection{Terminology}\label{introbasics}

We begin with some general terminology, and we introduce the notion of an 
``exceptional filling instruction" which is used throughout this paper.

Fix an orientable compact 3-manifold $X$ with $\partial X$ consisting of tori:

\begin{itemize}
\item A \emph{slope} on a boundary component $\tau$ of $X$ is the isotopy class of a non-trivial 
unoriented loop on $\tau$; 
\item A \emph{filling instruction} $\alpha$ for $X$ is a set consisting of either a slope or the 
empty set for each component of $\partial X$;
\item The \emph{filling} $X(\alpha)$ given by an instruction $\alpha$ is the manifold obtained by 
attaching one solid torus to 
$\partial X$ for each (non-empty) slope in $\alpha$, with the meridian of the solid torus attached to the slope.
\end{itemize}
We recall that if $M$ is a hyperbolic non-compact finite-volume 3-manifold then
$M=\text{int}(X)$ with $\partial X$ consisting of tori, and that $\text{int}(X(\alpha))$ 
is hyperbolic for all but finitely many $\alpha$'s consisting of one slope and $\varnothing$'s \cite{2pi}.
\begin{itemize}
\item If the interior of $X$ is hyperbolic but the interior of $X(\alpha)$ is not, we say that $\alpha$ is an \emph{exceptional filling instruction} 
for $X$ and that $X(\alpha)$ is an \emph{exceptional filling} of $X$;
\item We say that an exceptional filling instruction $\alpha'$ on a hyperbolic 3-manifold $X$ with boundary 
is \emph{properly contained} in $\alpha$, and write $\alpha' \subset \alpha$ if $\alpha'$ is contained in 
$\alpha$ and $\alpha'\neq \alpha$ (as sets of slopes).
\item We say that an exceptional filling instruction $\alpha$ on a hyperbolic $X$ is \emph{isolated} 
if $X(\alpha')$ is hyperbolic for all $g$ properly contained in $\alpha$; for such an $\alpha$ 
we call $X(\alpha)$ an \emph{isolated exceptional filling} of $X$.
\end{itemize}
A \emph{surgery} on a link $L$ corresponds to a 
filling of the exterior of $L$. That is, a surgery on $L$ is a filling of $M\backslash N(L)$ 
where $N(L)$ is an open regular neighbourhood of $L$. By a 
\emph{surgery instruction} for $L$ we mean a filling instruction 
on the exterior of $L$. 

We now recall some standard notation used in the description of the set 
of exceptional slopes on a fixed boundary component of a hyperbolic manifold, 
see for example \cite{Gor}. If $X$ is a hyperbolic 3-manifold with boundary consisting of tori and $\tau$ 
is a fixed boundary component of $\partial X$ then the set of exceptional slopes 
on $\tau$ is denoted by $E_{\tau}(X)$, and the cardinality of $E_{\tau}(X)$ by 
$e_{\tau}(X)$. The subscript $\tau$ is dropped whenever the boundary component is clear. 
To describe $E_{\tau}(M_5(\alpha))$ we introduce the following definition:

\begin{definition}
Let $\alpha$ be a filling instruction on a manifold $X$. We say that $\alpha$ 
\emph{factors through} a manifold $Y$ if there exists some filling instruction 
$\alpha' \subseteq \alpha$ such that $Y=X(\alpha')$. 
\end{definition}

We remarked above that a $-1$ filling on any boundary component of $M_n$ results in $M_{n-1}$. Therefore, 
any filling instruction on $M_n$ that contains a $-1$ slope factors through $M_{n-1}$. Note that if 
$\alpha$ is exceptional for $X$ and factors through a hyperbolic $Y$ with $Y=X(\alpha')$, then 
$\alpha \backslash \alpha'$ is exceptional for $Y$.

\subsection{Notation}\label{notation}

Our description of the exceptional fillings of $M_5(\alpha)$ will employ the notation now discussed 
for Seifert manifolds with orientable base surface. 
Given integers $p_1,\dots, p_n, q_1,\dots, q_n$, with $p_i$ and $q_i$ coprime, and $G$ an 
orientable surface with $k\geq0$ boundary components $b_1,\dots, b_k$, we let $\Sigma$ denote the 
surface obtained by removing $n$ open discs from $G$ and we denote by $b_{k+1}, \dots, b_{k+n}$ 
the $n$ newly introduced boundary circles. We fix an orientation on $\Sigma\times S^1$ and 
orient $\{\mu_i, \lambda_i\}=\{b_i\times \{\ast\}, \{\ast\}\times S^1\}$ 
so that $\mu_i$, $\lambda_i$ is a positive basis of $H_1(b_i \times S^1)$ with $b_i \times S^1$ oriented 
as $\partial (\Sigma\times S^1)$. We denote by $(G, (p_1, q_1), \dots (p_n, q_n))$, the manifold 
obtained by performing a Dehn filling on each $b_i \times S^1$ along $p_i\mu_i + q_i\lambda_i$ for $i>k$. 
In our case, $G$ will be either the disc $D$, the annulus $A$, or the sphere $S^2$.

Given Seifert manifolds $X$ and $Y$ with orientable base surfaces with boundary as described above, 
and $B\in\textrm{GL}(2,\mathbb{Z})$ with $\det(B)=-1$, we define $X\bigcup_B Y$ unambiguously to be the quotient 
manifold $X\bigcup_f Y$ where $f:T\ \rightarrow U$ for $T$ and $U$ arbitrary boundary components of $X$ 
and $Y$ respectively, and $f$ acting on homology by $B$ with respect to the bases described above.
The case $T, U\subset \partial X$, $T\neq U$ is also allowed and we write the quotient manifold 
as $X\big/_B$. 

The JSJ decomposition and Geometrization theorems tell us that every 
non-hyperbolic 3-manifold not homeomorphic to the 3-ball either contains an essential 
sphere, disc, torus, annulus, or is a closed small Seifert space. The closed small Seifert 
spaces are precisely those manifolds with Heegaard genus 0, Heegaard genus 1 or fibres 
over the sphere with exactly 3 exceptional fibres. Following \cite{Gor} we now assign names to each 
class of non-hyperbolic manifolds:

\begin{itemize}
\item The class of Heegaard genus 0 manifolds (\emph{i.e.} $\{S^3\}$) is 
denoted by $S^H$;
\item The class of all reducible 3-manifolds is denoted by $S$;
\item The class of manifolds with Heegaard genus 1 (\emph{i.e.}\ lens spaces) is denoted by $T^H$;
\item The class of manifolds containing an essential torus is denoted by $T$;
\item The class of boundary reducible manifolds is denoted by $D$;
\item The class of manifolds containing an essential annulus is denoted by $A$;
\item The class of Seifert spaces fibering over the sphere with exactly 
three exceptional fibres is denoted by $Z$.
\end{itemize}

We will say that a manifold in a class $\mathcal{C}$ is of \emph{type} $\mathcal{C}$.
We remark that the above classes are not mutually exclusive, for example 
$(D^2\times S^1) \# (D^2\times S^1)$ is of type $S$ and of type $D$.

\subsection{Surgery instructions on the chain links}\label{sicl} 

\begin{figure}[h!]
\begin{center}
\includegraphics[width = 12 cm]{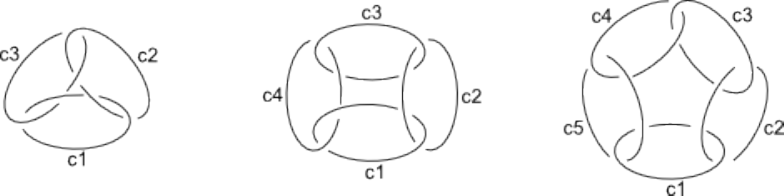}
\end{center}
\caption{A cyclic ordering of the components of 3CL, 4CL, 5CL.}
\label{nCL_ordering}
\end{figure}

We now explain the convention used to describe surgeries on the chain links. 
By ordering the components of $n$CL for $3\leq n\leq5$ cyclically as in Figure \ref{nCL_ordering}, surgery instructions on $n$CL can 
be naturally identified with $(\mathbb{Q}\cup\{\varnothing, \infty\})^n$.  By $M_n\left(x_1,\dots, x_n\right)$ we 
mean the manifold obtained by performing an $x_i$-surgery on the $i^{th}$ component of $n$CL.

To establish Theorem \ref{flash} we will examine all $M_5(\alpha)$. To avoid additional work we introduce 
the following definition which allows us to identify distinct surgery instructions that correspond to the 
same surgery.

\begin{definition}
Let $\alpha,\alpha'$ be filling instructions on a 3-manifold $X$ with toroidal boundary components. 
We will say that $\alpha$ and $\alpha'$ are \emph{equivalent} and write $\alpha\sim \alpha'$ when there exists a 
$h\in X(\alpha)\rightarrow X(\alpha')$.
\end{definition}

The appendix of \cite{Magic} contains a comprehensive analysis of the set of exceptional slopes on all $M_3(\alpha)$. 
Therefore, for the purposes of Theorem \ref{flash} we can omit the investigation of $E(M_5(\alpha))$ when $\alpha$ 
factors through $M_3$. As noted in the introduction, a positive twist about a boundary component of $M_n$ with a $-1$ 
slope results in $M_{n-1}$ for $n\geq2$. When we keep track of surgery coefficients we get (\ref{5CL_twist})-(\ref{4CL_twist}) 
(see \cite{MPR} for precise details). 
\begin{gather}
M_{5}\big(\tfrac pq, -1, \tfrac rs, \tfrac uv, \tfrac xy\big) = M_{4}\big(\tfrac{p+q}q, \tfrac {r+s}s, \tfrac {u}v, \tfrac{x}y\big)\label{5CL_twist} \\
M_{4}\big(\tfrac pq, \tfrac rs, -1, \tfrac uv\big) = M_{3}\big(\tfrac{p}q, \tfrac {r+s}s, \tfrac{u+v}{v}\big). \label{4CL_twist}
\end{gather}
Putting Identities (\ref{5CL_twist})-(\ref{4CL_twist}) we get:
\begin{equation}\label{rolfsen:thm}
M_{5}\big(\tfrac pq, -1, \tfrac rs, -1, \tfrac uv\big) = M_{3}\big(\tfrac{p+q}q, \tfrac {r+2s}s, \tfrac {u+v}v\big)
= M_{5}\big(\tfrac pq, -1, -2, \tfrac rs, \tfrac {u+v}v\big).
\end{equation}
Identities (\ref{5CL_twist})--(\ref{rolfsen:thm}) will be useful in Section \ref{mr}.

\subsection{The minimally twisted 4-chain link}\label{excep:fill:sec}

\begin{figure}[h!]
 \begin{center}
\includegraphics[width = 6 cm]{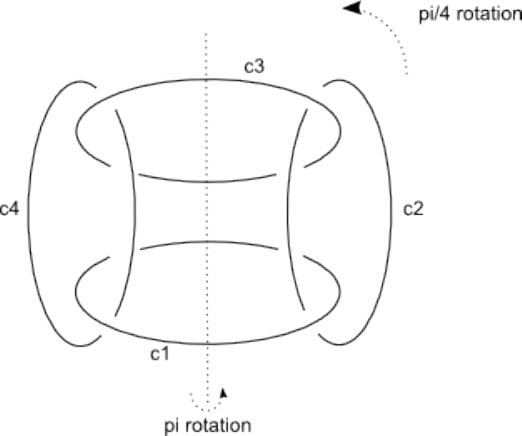}
 \end{center}
 \caption{The minimally twisted 4-chain link M4CL}
 \label{m4cl_pic}
\end{figure}

Most of the exceptional surgeries on 5CL are obtained by surgery on the minimally twisted 4-chain link 
M4CL shown in Figure \ref{m4cl_pic} (see Proposition \ref{M_3M_5:prop}). The proof of Theorem \ref{flash} will turn unto an investigation of 
the surgeries on M4CL. We will therefore need to understand the manifolds obtained by surgery on M4CL. 
Proposition \ref{fillingF} shows that all small Seifert spaces as well as many distinct reducible and 
toroidal manifolds are obtained by surgery on M4CL.
We denote the exterior of M4CL by $F$. As with $n$CL, we order the components of M4CL cyclically (see Figure \ref{m4cl_pic}) 
and equip each component of M4CL with the standard choice of meridian and longitude. 
Surgery instructions on M4CL are naturally identified with $(\mathbb{Q}\cup\{\varnothing, \infty\})^4$; by 
$F\left(\alpha_1,\alpha_2,\alpha_3,\alpha_4\right)$ we mean the manifold obtained by performing an $\alpha_i$-surgery on the $i^{th}$ component 
of M4CL.
It is easy to see from Figure \ref{m4cl_pic} that the symmetry group of M4CL contains the Dihedral group $D_4$. So, 
for any $\sigma\in D_4$ we have:
\begin{equation}\label{Fsym}
F\left(\alpha_1,\alpha_2,\alpha_3,\alpha_4\right)= F\left(\alpha_{\sigma(1)}, \alpha_{\sigma(2)}, \alpha_{\sigma(3)}, \alpha_{\sigma(4)}\right).
\end{equation} 
It is also useful for us to note that a negative twist about a boundary component of $M_5$ with a $+1$ 
slope results in $F$. When we keep track of surgery coefficients we get:
\begin{equation}\label{F_from_M_5}
M_{5}\big(\tfrac pq, 1, \tfrac rs, \tfrac uv, \tfrac xy\big) = F\big(\tfrac{p-q}q, \tfrac {r-s}s, \tfrac {u}v, \tfrac{x}y\big)
\end{equation} 
 
\begin{figure}[h!]
 \begin{center}
\includegraphics[width = 8 cm]{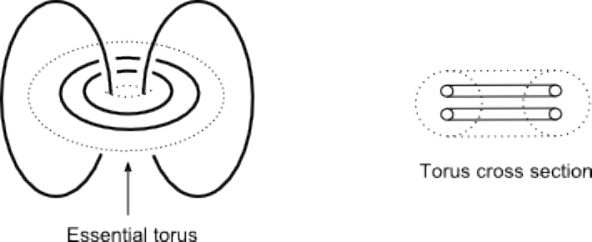}
 \end{center}
 \caption{$F$, the exterior of the minimally twisted 4-chain link.}
 \label{F_pic}
\end{figure}

\begin{figure}[h!]
 \begin{center}
\includegraphics[width = 8 cm]{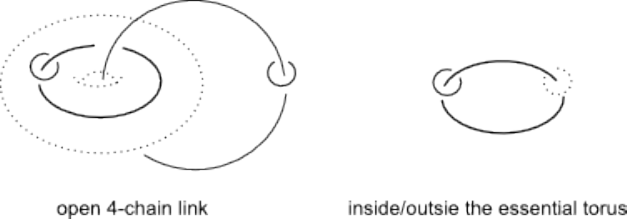}
 \end{center}
 \caption{$F$ realised as the exterior of the open 4-chain link.}
 \label{open4cl}
\end{figure}

Figure \ref{F_pic} highlights an exceptional torus $T$ in $F$, and it is clear that 
$T$ separates $F$ into two copies of $P\times S^1$ which are glued together by 
identifying a boundary component $\gamma \times S^1$ of one $P\times S^1$ to the 
other $P\times S^1$ with a horizontal loop $\gamma \times \{\ast\}$ in the former 
identified to a fibre $\{\ast\}\times S^1$ in the latter. Thus $F$ is homeomorphic 
to $P\times S^1\bigu 0110 P\times S^1$. 
We finally remark that $F$ is homeomorphic to the exterior of the open chain link with four 
components used to describe the exceptional surgeries of 5CL in \cite{MPR} (see Figure \ref{open4cl}).
  
To describe the fillings of the 4-chain link, we employ a flexible notation for 
Seifert manifolds. We will formally identify an $\varnothing$ slope with $\frac 00$ and allow all 
coprime pairs $(p_i,q_i)$ including $p_i, q_i\leq0$. 
Moreover, $S^2\times S^1$, $S^3$, $\mathbb{RP}^3$ will be written as $\lens pq$ with $p=0, 1,2$ respectively. 

\begin{proposition}\label{fillingF}
If $\alpha$ is a filling instruction on $F$ then up to (\ref{Fsym}) $\alpha$ is equivalent 
to some $(\frac pq,\frac rs, \frac uv, \frac xy)$ in Tables \ref{F1slope}-\ref{F4slopes} 
and $F(\frac pq,\frac rs, \frac uv, \frac xy)$ is described in 
Tables \ref{F1slope}-\ref{F4slopes}.
\end{proposition}

\begin{table}[h!] 
\begin{center}
\begin{tabular}{|c|c|c|c|c|c|}
\hline $\frac pq$ & $\frac rs$ & $\frac uv$ & $\frac xy$ & 
$F(\frac pq,\frac rs, \frac uv, \frac xy)$  & Type\\ \hline \hline
$0$ & $\varnothing$ & $\varnothing$ & $\varnothing$ & 
$D^2\times S^1 \# A\times S^1$ & $A$, $D$,
$S$\\ \hline 
$\frac1n$ & $\varnothing$ & $\varnothing$ & $\varnothing$ & $P\times S^1$ & 
$A$\\ \hline 
$|p|\geq2$ & $\varnothing$ & $\varnothing$ & $\varnothing$ & 
$\seifuno{A}pq\bigu0110 P\times S^1$ & $A$, $T$ \\ \hline
\end{tabular}
\end{center}
\caption{All manifolds obtained by filling $F$ along a single slope.}
\label{F1slope}
\end{table}

\begin{table}[h!] 
\begin{center}
\begin{tabular}{|c|c|c|c|c|c|}
\hline $\frac pq$ & $\frac rs$ & $\frac uv$ & $\frac xy$ & 
$F(\frac pq,\frac rs, \frac uv, \frac xy)$ & Type \\ \hline \hline

\multirow{3}{*}{0} & $\frac rs$ & $\varnothing$ & $\varnothing$ & 
$D^2\times S^1\# D^2\times S^1$ & $D$, $S$ \\ \cline{2-6}
& $\varnothing$ & $\frac uv$, $|u|\neq1$ & $\varnothing$ & $\lens uv\#A\times S^1$ & $A$, $S$ 
\\ \cline{2-6}
& $\varnothing$ & $\frac1n$ & $\varnothing$ & $A\times S^1$ & $A$\\ 
\hline

\multirow{5}{*}{$\frac1n$} & $|r|\geq2$ & $\varnothing$ & $\varnothing$ & 
$\seifuno{A}{r}{s}$ & $A$\\ \cline{2-6}
& $\frac rs=\frac 1k$ & $\varnothing$ & $\varnothing$ & $A\times S^1$ & $A$ \\ \cline{2-6}
& $\varnothing$ & $|v+nu|\geq2$ & $\varnothing$ & 
$\seifuno{A}{v+nu}{-u}$ & $A$\\ \cline{2-6}
& $\varnothing$ & $|v+nu|=1$ & $\varnothing$ & $A\times S^1$ & $A$\\ \cline{2-6}
& $\varnothing$ & $-\frac1n$ & $\varnothing$ & $S^2\times S^1\# A\times S^1$ & $A$, $S$\\ \hline 
 
\multirow{2}{*}{$|p|\geq2$} & $|r|\geq2$ & $\varnothing$ & $\varnothing$ 
& $\seifuno{A}pq\bigu0110\seifuno{A}rs$ & $A$, $T$\\ 
\cline{2-6}
& $\varnothing$ & $|u|\geq2$ & $\varnothing$ & $\seifdue{D^2}pquv\bigu0110 P\times S^1$ 
& $A$, $T$
\\ \hline
 
\end{tabular}
\end{center}
\caption{All manifolds obtained by filling $F$ along two slopes.}
\label{F2slopes}
\end{table}

\begin{table}[h!] 
\begin{center}
\begin{tabular}{|c|c|c|c|c|c|}
\hline $\frac pq$ & $\frac rs$ & $\frac uv$ & $\frac xy$ & 
$F(\frac pq,\frac rs, \frac uv, \frac xy)$ & Type \\ \hline \hline

\multirow{2}{*}{0} & \multirow{2}{*}{$\frac rs$} & $\frac uv$, $|u|\neq1$ & $\varnothing$ 
& $\lens uv\# D^2\times S^1$ & $D$, $S$ 
\\ \cline{3-6} 
& & $\frac1n$ & $\varnothing$ & $D^2\times S^1$ & $D$ \\ \hline

\multirow{6}{*}{$\frac 1n$} & $\frac1k$ & $\frac uv$ & $\varnothing$ & $D^2\times S^1$ & $D$ \\ \cline{2-6}

& \multirow{2}*{$\frac rs$, $|r|>1$} & $\frac uv=-\frac1n$ 
& $\varnothing$ & $\lens rs\# D^2\times S^1$ & $D$, $S$ \\ \cline{3-6}
&  & $|v+nu|=1$ & $\varnothing$ & $D^2\times S^1$ & $D$ \\ \cline{2-6}

& \multirow{2}{*}{$\frac rs=0$} & $\frac uv$, $|v+nu|\neq1$ & $\varnothing$ & 
$\lens{v+nu}{-u}\# D^2\times S^1$ & $D$, $S$ \\ \cline{3-6}
& & $\frac uv$, $|v+nu|=1$ & $\varnothing$ & $D^2\times S^1$ & $D$\\
\cline{2-6} 
& $|r|\geq2$ & $\frac uv$ & $\varnothing$ & 
$\seifdue{D}rs{v+nu}{-u}$ & $A$ \\ \hline

\multirow{3}{*}{$|p|\geq2$} & $|r|\geq2$ & $|u|\geq2$ & $\varnothing$ 
& $\seifdue{D}pquv\bigu0110\seifuno{A}rs$ & $A$, $T$ \\ \cline{2-6} 
& $\frac1n$ & $|u|\geq2$ & $\varnothing$ 
& $\seifdue{D}pquv$ & $A$ \\ \cline{2-6} 
& $0$ & $|u|\geq2$ & $\varnothing$ & $D\times S^1\# \lens{pv+qu}{pv'+qu'}$ & $D$, $S$ \\
& & & & where $|uv'-vu'|=1$ & \\ \hline 

\end{tabular}
\end{center}
\caption{All manifolds obtained by filling $F$ along three slopes.}
\label{F3slopes}
\end{table}

\begin{table}[h!] 
\begin{center}
\scalebox{0.86}{
\begin{tabular}{|c|c|c|c|c|c|}
\hline $\frac pq$ & $\frac rs$ & $\frac uv$ & $\frac xy$ & 
$F(\frac pq,\frac rs, \frac uv, \frac xy)$ & Type \\ \hline \hline

\multirow{6}{*}{0} & \multirow{2}{*}{$\frac rs$} & $\frac uv$ & $\frac xy$ & 
$\lens uv\# \lens{ry+sx}{xs'+yr'}$ & \multirow{2}{*}{$S$} \\
& & $|u|\neq1$ & $|ry+sx|\neq1$ & where $|rs'-sr'|=1$ & \\ \cline{2-6}

& \multirow{2}{*}{$\frac rs$} & \multirow{2}{*}{$\frac 1v$} & \multirow{2}{*}{$\frac xy$} & 
$\lens{ry+sx}{xs'+yr'}$ & \multirow{2}{*}{$T^H/S^H$} \\
& & &  & where $|rs'-sr'|=1$ & \\ \cline{2-6}

& \multirow{2}{*}{$\frac rs$} & \multirow{2}{*}{$\frac uv$} & $\frac xy$ & 
$\lens uv$ & \multirow{2}{*}{$T^H/S^H$} \\
& & & $|ry+sx|=1$ & where $|rs'-sr'|=1$ & \\ \hline

\multirow{8}{*}{$\frac 1n$} & \multirow{2}{*}{$\frac 1k$} & \multirow{2}{*}{$\frac uv$} 
& \multirow{2}{*}{$\frac xy$} & $\lens{(v+nu)(y+kx)-xu}{(v+nu)j-ui}$ & \multirow{2}{*}{$T^H/S^H$}\\ 
& & & & where $|xj-(y+kx)i|=1$ &\\ \cline{2-6}

& \multirow{2}{*}{$\frac rs$} & $v+nu=\epsilon$ 
& \multirow{2}{*}{$\frac xy$} & $\lens{ry+(x-\epsilon ru)x}{rj+(s-\epsilon ru)i}$ & \multirow{2}{*}{$T^H/S^H$} \\ 
& & $(\epsilon=\pm1)$ &  & where $|xj-yi|=1$ & \\ \cline{2-6}

& $\frac rs$ & \multirow{2}{*}{$-\frac 1n$} 
& $\frac xy$ & \multirow{2}{*}{$\lens rs\#\lens xy$} & \multirow{2}{*}{$S$} \\ 

& $|r|>1$ &  & $|x|>1$ & & \\ \cline{2-6}

& $\frac rs$ & $\frac uv$ 
& $\frac xy$ & \multirow{2}{*}{$\seiftre{S^2}{v+nu}{-u}rsxy$} & \multirow{2}{*}{$Z$}\\ 
& $\frac rs\neq\frac1k$ & $|v+nu|\geq2$ & $|x|\geq2$ & & \\ \hline

\multirow{2}{*}{$|p|\geq2$} & \multirow{2}{*}{$|r|\geq2$} & \multirow{2}{*}{$|u|\geq2$} 
& \multirow{2}{*}{$|x|\geq2$} 
& $\seifdue{D}pquv\bigu0110\seifdue{D}rsxy$ & \multirow{2}{*}{$T$} \\ \hline 
\end{tabular}
}
\end{center}
\caption{All manifolds obtained by filling $F$ along four slopes.}
\label{F4slopes}
\end{table}

\begin{proof} As noted above, the symmetry group of M4CL has a $D_4$ 
action on the boundary components of M4CL. Therefore, by (\ref{Fsym}), we may assume that 
when a filling instruction on $F$ has a single slope it is of the form 
$(\frac pq,\varnothing, \varnothing, \varnothing)$, when a filling instruction 
on $F$ has two slopes it is of the form $(\frac pq,\frac rs, \varnothing, \varnothing)$ or 
$(\frac pq, \varnothing, \frac uv, \varnothing)$, and when a filling instruction 
on $F$ has three slopes it is of the form $(\frac pq,\frac rs, \frac uv, \varnothing)$. 

The slopes of a filling instruction $\alpha=(\frac pq,\frac rs, \frac uv, \frac xy)$ 
can be split according to whether the numerator is 0, $\pm1$, or greater than one in 
absolute value. Using the symmetry group of M4CL we see that if $0\in \alpha$ 
we may assume $\tfrac pq= 0$. If $0\not\subset\alpha$ and $\tfrac1n$ is a slope in $\alpha$ we 
may assume $\tfrac pq=\frac1n$. It is now easy to see that the description of 
possible slopes in Tables \ref{F1slope}-\ref{F4slopes} is exhaustive.

We have already seen that $F$ is the union of two copies of 
$P\times S^1$ glued together by the orientation reversing map sending 
a meridian to a longitude and a longitude to a meridian. Therefore, the 
Identities (\ref{F_filling1})-(\ref{F_filling6}) hold:
\begin{eqnarray}
F & = & P\times S^1\bigu0110 P\times S^1, \label{F_filling1} \\ 
F(\tfrac pq,\varnothing,\varnothing, \varnothing) & = & \seifuno{A}pq\bigu0110P\times S^1  \label{F_filling2}\\
F(\tfrac pq,\tfrac rs,\varnothing, \varnothing) & = & \seifuno{A}pq\bigu0110\seifuno{A}rs  \label{F_filling4}\\
F(\tfrac pq,\varnothing,\tfrac uv, \varnothing) & = & \seifdue{D}pquv\bigu0110P\times S^1  \label{F_filling3}\\
F(\tfrac pq,\tfrac rs,\tfrac uv, \varnothing) & = & \seifdue{D}pquv\bigu0110\seifuno{A}rs  \label{F_filling5}\\
F(\tfrac pq, \tfrac rs, \tfrac uv, \tfrac xy) & = & \seifdue{D}pquv\bigu0110\seifdue{D}rsxy. \label{F_filling6}
\end{eqnarray}

Keeping the conventions set out in Section \ref{notation}, the description of manifolds given in Tables 
\ref{F1slope}-\ref{F4slopes} is obtained by repeated use of well-known Identities (\ref{seifeq1})-(\ref{grapheq6}) between 
graph manifolds (see \cite{FomMat} for details) to Identities (\ref{F_filling1})-(\ref{F_filling6}) until each 
$F(\frac pq,\frac rs, \frac uv, \frac xy)$ is written as a $\lens pq$, or a graph manifold with exceptional 
fibres of the form $(a,b)$ with $|a|\geq2$ and $|b|\geq1$.

%with $\partial E\neq \varnothing$ we denote by $E'$ the surface obtained by gluing a copy of 
%$D^2$ to one of the boundary components of $E$. Details of the choice of bases on $\partial (E\times S^1)$ 
%can be found in Section \ref{excep:fill:sec}.

%Let $X$ and $Y$ be $3$-manifolds with fixed homology bases on $\partial X$ and $\partial Y$, 
%and $E$ be a surface with boundary, $G$ be a closed surface, and $p_i, q_i\in\mathbb{Z}$. 
%Then the following identities hold:

\textbf{Seifert manifolds:}
{\small \begin{equation}\label{seifeq1}
\big (G, (p_1, q_1), (p_2, q_2), \dots, (p_k, q_k)\big) =  \big (G, (p_1, q_1 -np_1), (p_2, q_2 +np_2), 
\dots, (p_k, q_k)\big),
\end{equation}
\begin{equation}\label{seifeq2}
\big (G, (1, q_1), (p_2, q_2), \dots, (p_k, q_k)\big) =  \big (G, (p_2, q_2+q_1p_2), \dots, (p_k, q_k)\big),
\end{equation}
\begin{equation}\label{seifeq3}
\big (G, (p_1, q_1), (p_2, q_2), \dots \big) =  \big (G, (p_1, q_1-np_1), (p_2, q_2), \dots\big) \text{ if } \partial G\neq\varnothing.
\end{equation}}

\textbf{Small Seifert manifolds:}
{\small \begin{equation}\label{smseifeq1}
(S^2, (p, q))=\lens qp,
\end{equation}
\begin{equation}\label{smseifeq2}
(S^2, (p, q), (r, s))=\lens{ps+rq}{ps'+r'q} \text{ where } rs'-sr'=1,
\end{equation}
\begin{equation}\label{smseifeq3}
\seiftre{S^2}01pqrs = \lens pq \# \lens rs.
\end{equation}}

\textbf{Graph Manifolds:} 
{\small \begin{equation}\label{grapheq1}
X\bigcup\nolimits_B\phantom{\Big|}\!\! Y=Y\bigcup\nolimits_{B^{-1}}\phantom{\Big|}\!\!X,
\end{equation}
\begin{equation}\label{grapheq2}
(D, (p, q))\bigu abcd (G, \dots)= (G', (ap-bq, cp-dq), \dots) \text{ where } G'\backslash\text{disc} = G,
\end{equation}
\begin{equation}\label{grapheq3}
\big(G, (p, q), \dots\big)\bigu abcd X = \big(G, (p, q+kp), \dots\big)\bigu {a+kb}b{c+kd}d X,
\end{equation}
\begin{equation}\label{grapheq4}
\seifdue{D^2}01pq \bigu abcd \big( E, \dots\big)= \lens pq\#\big(E', (b, d) \dots\big),
\end{equation}
\begin{equation}\label{grapheq6}
\big (E, (p_1, q_1), \dots, (p_k, q_k)\big) \bigu acbd X = \big (E, (p_1, -q_1), \dots, (p_k, -q_k)\big)\bigu {-a}c{-b}d X.
\end{equation}}

\end{proof}

\subsection{Exceptional surgery instructions on 5CL}\label{result:recap}

We now return to 5CL and present the concise description of exceptional surgery instructions given in \cite{MPR}. 
Given a surgery instruction $\alpha$ on 5CL, the symmetry group of $M_5(\alpha)$ induces a natural action on the 
boundary components of $M_5(\alpha)$. This action induces an action on the filling instructions on $M_5(\alpha)$. 
Among the most significant actions arising we mention those coming from the symmetry group of $M_5$, see 
(\ref{first:eqn})-(\ref{3rd:eqn}) below, a symmetry of $M_4$ which may be deduced from the Fenn-Rourke blow-down 
move on 5CL, see (\ref{blow:eqn}), and from the amphichirality of the Figure-8 knot $M_5(-1,-2,-2,-2,\varnothing)$, 
see (\ref{figure8:eqn}) (see \cite{MPR} for full details).
\begin{align}
(\alpha_1, \alpha_2,\alpha_3,\alpha_4,\alpha_5) & \longmapsto (\alpha_5, \alpha_1, \alpha_2, \alpha_3, \alpha_4) \label{first:eqn} \\
(\alpha_1, \alpha_2, \alpha_3, \alpha_4, \alpha_5) & \longmapsto (\alpha_5, \alpha_4, \alpha_3, \alpha_2, \alpha_1)
\label{second:eqn} \\
(\alpha_1, \alpha_2,\alpha_3,\alpha_4,\alpha_5) & \longmapsto (\alpha_2^{-1}, \alpha_1^{-1}, 1-\alpha_3, (1-\alpha_4)^{-1}, 1-\alpha_5)\label{3rd:eqn}\\
(-1, \alpha_1, \alpha_2, \alpha_3, \alpha_4) & \longmapsto (\alpha_1, -1, \alpha_2-1, \alpha_3, \alpha_4+1)\label{blow:eqn}\\
\left(-1, -2, -2, -2, \alpha \right) & \longmapsto \left(-1, -2, -2, -2, -\alpha - 6 \right)\label{figure8:eqn}
\end{align}

For a filling instruction $\alpha$ on $M_5$ we will often simplify notation by omitting empty slopes but leaving 
the subscripts on non-empty slopes. For example, $((-1)_2, (-1)_4)$ corresponds to the filling instruction 
$(\varnothing, \mu_2-\lambda_2, \varnothing, \mu_4-\lambda_4, \varnothing)$ with $(\mu_i, \lambda_i)$ the 
(meridian, longitude) basis of the homology of the $i^{\text{th}}$ cusp. Note that for any 
$\frac pq\in\mathbb{Q}\cup \{\infty\}$ and $i\neq j$ one has $[((\frac pq)_i)]=[((\frac pq)_j)]$ by (\ref{first:eqn}), 
so the fillings $M_5(\frac pq)$ are defined without ambiguity. Our convention will be that 
a filling instruction $(\frac pq, \frac rs, \frac uv, \frac xy)$ on $M_5$ with four non-empty 
slopes and no subscripts represents $({\frac pq}, {\frac rs}, {\frac uv}, {\frac xy}, \varnothing)$.
We now state the main result of \cite{MPR}.

\begin{thm} \label{main:teo}
Every exceptional filling instruction on $M_5$ is equivalent up to a composition of the symmetries 
\emph{(\ref{first:eqn})-(\ref{figure8:eqn})} to a filling instruction containing one of:
\begin{align*}
1, (-1, -2, -2, -1), & (-1, -2, -3, -2, -4), (-1, -2, -2, -3, -5),\\
(-1, -3, -2, -2, -3),& \left(-2, -\tfrac 12, 3, 3, -\tfrac 12\right), (-2, -2, -2, -2, -2).
\end{align*}
Moreover, the following equalities hold:
\begin{align*}
M_5\left(\tfrac ab, \tfrac cd, \tfrac ef, \tfrac gh, 1 \right) & = F(\tfrac{a-b}b, \tfrac cd, \tfrac ef, \tfrac{g-h}h)\\
M_5\left(-1,-2, -2, -1, \tfrac ab \right) & =  \seifuno{A}{b}{-a-b} \bigb 0110\\
M_5\left(-1, -2, -2, -3, -5\right) & = \seifdue{D}212{-1} \bigu {-1}21{-1} \seifdue{D}2131 \\
M_5\left(-1, -2, -3, -2, -4\right) & = \seifdue{D}212{-1} \bigu {-1}31{-2} \seifdue{D}2131 \\
M_5\left(-1, -3, -2, -2, -3\right) & = \seifdue{D}212{-1} \bigu {-1}41{-3} \seifdue{D}2131 \\
M_5\left(-2, -2, -2, -2, -2\right) & =\seifdue{D}212{-1} \bigu {-1}51{-4} \seifdue{D}2131 \\
M_5\left(-2, -\tfrac 12, 3, 3, -\tfrac 12\right) & = \seifuno A 2{-1} \bigb 1211.
\end{align*}
\end{thm}

\begin{remark} The surgery instructions 
$(-1,-2, -2, -1,\varnothing)$, $(-1, -2, -2, -3, -5)$, and $(-1, -2, -3, -2, -4)$ 
factor through $M_3$, and every non-toroidal exceptional filling of $M_5$ is obtained by filling $F$. 
Furthermore, each isolated exceptional filling of $M_5$ is a graph manifold and therefore a filling 
instruction $\alpha$ on $M_5$ is exceptional if and only if $\alpha$ contains an isolated exceptional 
surgery instruction.
\end{remark}

We now recall the classification of isolated exceptional surgeries on 3CL found in \cite{Magic}. It is easy to 
see that the symmetry group of 3CL is $S_3$ and so for a filling instruction $\alpha$ on $M_3$ we can write $M_3(\alpha)$ 
unambiguously. 

\begin{thm}\label{3CL_thm}\emph{[Martelli, Petronio]}
Up to the $S_3$ action on the components of \emph{3CL}, a surgery instruction on \emph{3CL} is an isolated exceptional surgery 
instruction if and only if it is one of
\begin{gather*}
 \infty, 0, 1, 2, 3, (-1, -1), (4, \tfrac12), (\tfrac32,\tfrac52), (5,5,\tfrac12), (4,4,\tfrac23), (4,\tfrac32, \tfrac32), 
(4,\tfrac13, -1), (\tfrac83, \tfrac32, \tfrac32),\\ 
(\tfrac52,\tfrac52, \tfrac43), (\tfrac52,\tfrac53,\tfrac53), (\tfrac73,\tfrac73,\tfrac32), (-1,-2,-2), (-1,-2,-3), 
(-1,-2,-4),(-1,-2,-5), \\ 
(-1,-3,-3), (-2,-2,-2).
\end{gather*}
\end{thm}

We remark that all non-hyperbolic $M_3(\alpha)$ are described in Theorems 1.1--1.3 in \cite{Magic}. To keep the presentation in this 
article self-contained we describe the exceptional fillings in terms of the manifold $F$ (see Proposition \ref{M_3M_5:prop}). 
To do this the following lemma will be very helpful.

\begin{lemma}\label{sym} The action of \emph{Aut}$(M_5)$ on surgery instructions on \emph{5CL} is generated by 
\emph{(\ref{linksymeq1})--(\ref{symeq11})}.
Moreover, for $\ref{symeq1}\leq n \leq \ref{symeq11}$ each \emph{($n$)} corresponds to the action of a distinct element of 
\emph{Aut}$(M_5)/G$ where $G$ is the subgroup generated by the elements \emph{(\ref{linksymeq1})--(\ref{linksymeq2})} 
corresponding to the generators of the link symmetry group of \emph{5CL}.
\begin{equation}\label{linksymeq1}
\left( \alpha_1, \alpha_2, \alpha_3, \alpha_4, \alpha_5 \right) \, \, \, \longmapsto \, \, \, 
\left( \alpha_5, \alpha_1, \alpha_2, \alpha_3, \alpha_4 \right) 
\end{equation}
\begin{equation}\label{linksymeq2}
\left( \alpha_1, \alpha_2, \alpha_3, \alpha_4, \alpha_5 \right) \, \, \, \longmapsto \, \, \, \left( \alpha_5, \alpha_4, \alpha_3, \alpha_2, \alpha_1 \right)\end{equation}
\begin{equation}\label{symeq1}
\big( \tfrac{a}{b}, \tfrac{c}{d}, \tfrac{e}{f}, \tfrac{g}{h}, \tfrac{i}{j} \big) \, \, \, \longmapsto \, \, \, \big( \tfrac{f}{e}, \tfrac{j-i}{j}, \tfrac{a}{a-b}, \tfrac{d-c}{d}, \tfrac{h}{g} \big) 
\end{equation}
\begin{equation}\label{symeq2}
\big( \tfrac{a}{b}, \tfrac{c}{d}, \tfrac{e}{f}, \tfrac{g}{h}, \tfrac{i}{j} \big) \, \, \, \longmapsto \, \, \, \big( \tfrac{b}{b-a}, \tfrac{i-j}{i}, \tfrac{e-f}{e}, \tfrac{d}{d-c}, \tfrac{g}{h} \big) 
\end{equation}
\begin{equation}\label{symeq3}
\big( \tfrac{a}{b}, \tfrac{c}{d}, \tfrac{e}{f}, \tfrac{g}{h}, \tfrac{i}{j} \big) \, \, \, \longmapsto \, \, \, \big( \tfrac{i}{i-j}, \tfrac{b-a}{b}, \tfrac{f}{e}, \tfrac{d}{c}, \tfrac{h-g}{h} \big) 
\end{equation}
\begin{equation}\label{symeq4}
\big( \tfrac{a}{b}, \tfrac{c}{d}, \tfrac{e}{f}, \tfrac{g}{h}, \tfrac{i}{j} \big) \, \, \, \longmapsto \, \, \, \big( \tfrac{j}{j-i}, \tfrac{e}{f}, \tfrac{b}{b-a}, \tfrac{c-d}{c}, \tfrac{g-h}{g} \big)  
\end{equation}
\begin{equation}\label{symeq5}
\big( \tfrac{a}{b}, \tfrac{c}{d}, \tfrac{e}{f}, \tfrac{g}{h}, \tfrac{i}{j} \big) \, \, \, \longmapsto \, \, \, \big( \tfrac{a}{a-b}, \tfrac{e}{e-f}, \tfrac{i}{i-j}, \tfrac{c}{c-d}, \tfrac{g}{g-h} \big) 
\end{equation}
\begin{equation}\label{symeq6}
\big( \tfrac{a}{b}, \tfrac{c}{d}, \tfrac{e}{f}, \tfrac{g}{h}, \tfrac{i}{j} \big) \, \, \, \longmapsto \, \, \, \big( \tfrac{h}{g}, \tfrac{j}{i}, \tfrac{f-e}{f}, \tfrac{c}{c-d}, \tfrac{b-a}{b} \big)
\end{equation}
\begin{equation}\label{symeq7}
\big( \tfrac{a}{b}, \tfrac{c}{d}, \tfrac{e}{f}, \tfrac{g}{h}, \tfrac{i}{j} \big) \, \, \, \longmapsto \, \, \, \big( \tfrac{h}{h-g}, \tfrac{a}{b}, \tfrac{f}{f-e}, \tfrac{c-d}{c}, \tfrac{i-j}{i} \big) 
\end{equation}
\begin{equation}\label{symeq8}
\big( \tfrac{a}{b}, \tfrac{c}{d}, \tfrac{e}{f}, \tfrac{g}{h}, \tfrac{i}{j} \big) \, \, \, \longmapsto \, \, \, \big( \tfrac{g}{g-h}, \tfrac{f-e}{f}, \tfrac{b}{a}, \tfrac{d}{c}, \tfrac{j-i}{j} \big)
\end{equation}
\begin{equation}\label{symeq9}
\big( \tfrac{a}{b}, \tfrac{c}{d}, \tfrac{e}{f}, \tfrac{g}{h}, \tfrac{i}{j} \big) \, \, \, \longmapsto \, \, \, \big( \tfrac{g-h}{g}, \tfrac{f}{f-e}, \tfrac{i}{j}, \tfrac{d}{d-c}, \tfrac{a-b}{a} \big) 
\end{equation}
\begin{equation}\label{symeq10}
\big( \tfrac{a}{b}, \tfrac{c}{d}, \tfrac{e}{f}, \tfrac{g}{h}, \tfrac{i}{j} \big) \, \, \, \longmapsto \, \, \, \big( \tfrac{h-g}{h}, \tfrac{b}{a}, \tfrac{j}{i}, \tfrac{d-c}{d}, \tfrac{e}{e-f} \big)
\end{equation}
\begin{equation}\label{symeq11}
\big( \tfrac{a}{b}, \tfrac{c}{d}, \tfrac{e}{f}, \tfrac{g}{h}, \tfrac{i}{j} \big) \, \, \, \longmapsto \, \, \, \big( \tfrac{a-b}{a}, \tfrac{e-f}{e}, \tfrac{h}{h-g}, \tfrac{c}{d}, \tfrac{j}{j-i} \big).
\end{equation}
\end{lemma}

\begin{remark}
When a symbol appearing in the argument of one of these maps is $\varnothing$, there is a corresponding symbol $\varnothing$ in the image of the map.
\end{remark}

\begin{proof}
SnapPy computes the symmetry group of $M_5$ to be $S_5\times \mathbb{Z}/2\mathbb{Z}$ 
and determines its action on slopes (see \cite{Snappy}). 
The maps (\ref{linksymeq1})-(\ref{symeq11}) come directly from SnapPy. 

Each element of the symmetry group of $M_5$ acts on the set of 
boundary components of $M_5$, and on the set of filling instructions. We will first demonstrate 
that the action on the set of boundary components of $M_5$ 
generated by the maps (\ref{linksymeq1})-(\ref{symeq11}) is that 
of the full $S_5$, and then we will use this fact to conclude that 
the action on surgery instructions on 5CL induced by the symmetry group of $M_5$ 
is generated by the maps (\ref{linksymeq1})-(\ref{symeq11}). 

No two of (\ref{symeq1})-(\ref{symeq11}), considered as permutations of the 
boundary components, are equal up to the $D_5$ action from the link symmetry 
group of 5CL generated by (\ref{linksymeq1}) and (\ref{linksymeq2}). 
Thus, each of (\ref{symeq1})-(\ref{symeq11}) corresponds 
to a representative element of a distinct left coset of $D_5$ in $S_5$. 
Since there are $\frac{5!}{10}=12$ such cosets, and our list of maps 
(\ref{symeq1})-(\ref{symeq11}) consists of 11 items, the symmetries of 
$M_5$ corresponding to (\ref{linksymeq1})-(\ref{symeq11}) generate the 
full $S_5$ action. 

We recall that SnapPy computes the order of the symmetry group of $M_5$ to be 240; 
the generator of the remaining $\mathbb{Z}/2\mathbb{Z}$ corresponds to a strong 
involution of 5CL with a trivial action on the set of filling instructions. Thus the equivalence 
relation on filling instructions induced by the action of the symmetry group of $M_5$ 
is generated by the maps (\ref{linksymeq1})-(\ref{symeq11}). 
\end{proof}

Since the action of the maps (\ref{first:eqn}) and (\ref{second:eqn}) is very easily understandable 
while that of (\ref{3rd:eqn})-(\ref{figure8:eqn}) is more involved, we introduce the symbol 
$[\alpha]$ for the equivalence class of a filling instruction $f$ under (\ref{first:eqn}) 
and (\ref{second:eqn}) only, and the symbol $[\![\alpha]\!]$ for the equivalence class of $\alpha$ 
under (\ref{first:eqn})-(\ref{figure8:eqn}). Note that, $[\alpha]\subseteq[\![\alpha]\!]$ and if 
$\alpha_1, \alpha_2\in [\![\alpha]\!]$ then $M_5(\alpha_1)=M_5(\alpha_2)$.

\section{Main Results}\label{mr}

We now precisely describe the classification of exceptional slopes of every hyperbolic surgery on 5CL by 
describing $E_{\tau}(M_5(\alpha))$ for every $\alpha$ not factoring through $M_3$. We split the result according 
to whether or not $\alpha$ factors through $M_4$.

\begin{thm}\label{5CL_thm} Let $\alpha$ be a hyperbolic filling instruction on $M_5$ containing 
at least one $\varnothing$ and not factoring through $M_4$, and let $\tau$ be a boundary component 
of $\partial M_5(\alpha)$. Either $e_{\tau}(M_5(\alpha))=3$ and, with respect to the basis induced 
from $M_5$, we have $E_{\tau}(M_5(\alpha))=\{0,1,\infty\}$ and 
\begin{gather}
M_5(\tfrac ab, \tfrac cd, \tfrac ef, \tfrac gh)(\infty)= F(-\tfrac ab, \tfrac fe, \tfrac dc, -\tfrac gh) \label{5CLinf} \\
M_5(\tfrac ab, \tfrac cd, \tfrac ef, \tfrac gh)(1)= F(\tfrac{a-b}b, \tfrac cd, \tfrac ef, \tfrac{g-h}h) \label{5CL1} \\
M_5(\tfrac ab, \tfrac cd, \tfrac ef, \tfrac gh)(0)= F(\tfrac{b}{b-a},\tfrac{c-d}c,-\tfrac hg, \tfrac{e-f}f)
\label{5CL0}
\end{gather}
or $4\leqslant e_{\tau}(M_5(\alpha)) \leqslant 5$ and 
\begin{itemize}
\item $M_5(\alpha)$ is homeomorphic to $M_5(f)$ for some $f$ in Tables \emph{\ref{tb1}}-\emph{\ref{tb5}} and, 
with respect to the basis induced from $M_5$, we have
\begin{equation*}
E_{\tau}(M_5(f))=
\begin{cases}
\{\beta_1,\beta_2, 0,1,\infty\} \mbox{ if $f$ is found in Table \emph{\ref{tb1}}};\\ 
\{\beta, 0,1,\infty\} \mbox{ otherwise},
\end{cases}
\end{equation*}
where $\beta, \beta_1,\beta_2$ are found in Tables \emph{\ref{tb1}}-\emph{\ref{tb5}}.
\item Identities \emph{(\ref{5CLinf})-(\ref{5CL0})} hold and $M_5(f)(\beta_i)$, 
$M_5(f)(\beta)$ are explicitly described in Tables \emph{\ref{tb1}}-\emph{\ref{tb5}}.
\end{itemize} 
\end{thm}

Adopting the convention that a filling instruction $(\frac ab, \frac cd, \frac ef)$ on $M_4$ 
with three slopes and no subscripts represents $({\frac ab}, {\frac cd}, {\frac ef}, \varnothing)$, 
the corresponding result for $M_4$ is:

\begin{thm}\label{4CL_thm} 
Let $\alpha$ be a hyperbolic filling instruction on $M_4$ containing at least one $\varnothing$ 
and not factoring through $M_3$, and let $\tau$ a boundary component of $\partial M_4(\alpha)$. 
Either $e_{\tau}(M_4(\alpha))=4$ and, with respect to the basis induced from $M_4$, we have 
$E_{\tau}(M_4(\alpha))=\{0,1,2,\infty\}$ and 
\begin{gather}
M_4(\tfrac ab, \tfrac cd, \tfrac ef)(\infty)= F(\tfrac{b-a}b, \tfrac{d}{c-d}, -1, -\tfrac ef) \label{4CLinf} \\
M_4(\tfrac ab, \tfrac cd, \tfrac ef)(2)= F(\tfrac{a-b}b, \tfrac cd, \tfrac{e-f}f,-2) \label{4CL2} \\
M_4(\tfrac ab, \tfrac cd, \tfrac ef)(1)= F(\tfrac{a-2b}b,-1,\tfrac{c-d}d,\tfrac{e-f}f), \label{4CL1} \\
M_4(\tfrac ab, \tfrac cd, \tfrac ef)(0)= F(\tfrac{2d-c}{c-d},\tfrac{b}{a-2b}, 2,\tfrac fe), \label{4CL0}
\end{gather}
or $5\leqslant e_\tau(M_4(\alpha))\leqslant 6$ and 
\begin{itemize}
\item $M_4(\alpha)$ is homeomorphic to $M_4(f)$ for some $f$ in Table \emph{\ref{tb6}}, 
and with respect to the basis induced from $M_4$, we have
\begin{equation*}
E_{\tau}(M_4(f))=
\begin{cases}
\{\beta_1,\beta_2, 0,1,2,\infty\} \mbox{ if }\alpha=(-2,-2,-2),(2,-\frac12,2);\\ 
\{\beta, 0,1,2,\infty\} \mbox{ otherwise},
\end{cases}
\end{equation*}
where $\beta, \beta_1,\beta_2$ are found in Table \emph{\ref{tb6}}.
\item Identities \emph{(\ref{4CLinf})-(\ref{4CL0})} hold, and the $M_5(f)(\beta_i)$, 
$M_4(f)(\beta)$ are explicitly described in Table \emph{\ref{tb6}}.
\end{itemize} 
\end{thm}

%While proving Theorems \ref{5CL_thm} and \ref{4CL_thm} we will often see 
%identities between fillings of the $M_i$. When a homeomorphism $M_5(\alpha)\simeq M_5(\beta)$ follows 
%as a consequence of maps (\ref{first:eqn})-(\ref{figure8:eqn}) we will write $M_5(\alpha)\sim M_5(\beta)$. 
%In addition to realising homeomorphisms between fillings of $M_5$, we will often use the Rolfsen twist 
%(see \cite{Rolfsen}) to obtain homeomorphisms between fillings of the $M_i$; when 
%$M_5(\alpha)\simeq M_k(\beta)$ is induced from a (composition of) Rolfsen twist(s), we will write 
%$M_5(\alpha)\simeq M_k(\beta)$.

\subsection{Proofs of main results}\label{mr_proofs}

In Theorem \ref{main:teo} we have a complete description of the exceptional instructions and fillings 
of $M_5$. As $M_3$ is obtained by surgery on 5CL, Theorem \ref{main:teo} contains an opaque classification 
of exceptional surgery instructions contained in Theorem \ref{3CL_thm}. Having an explicit classification 
of the exceptional fillings of $M_3$ will be important for the proof of Theorem \ref{flash}. Proposition 
\ref{fillingF} can be used to give a complete description of all exceptional fillings of $M_3$. The description 
of exceptional fillings in Proposition \ref{M_3M_5:prop} is the same as that given in \cite{Magic} up to 
(\ref{seifeq1})-(\ref{grapheq6}). The description of the exceptional fillings in Tables 
\ref{tb1}-\ref{families_table} comes from Proposition \ref{M_3M_5:prop} not \cite{Magic}.

To prove Proposition \ref{M_3M_5:prop} we will often use (\ref{linksymeq1})--(\ref{linksymeq2}) in conjunction 
with previous identities to prove the main results in Section \ref{mr}. As (\ref{linksymeq1})--(\ref{linksymeq2}) 
are easy to understand  we do not indicate when (\ref{linksymeq1})--(\ref{linksymeq2}) have been used. For 
example instead of writing
\begin{gather*}
M_5(\tfrac ab,\tfrac cd, \tfrac ef, \tfrac gh, \infty)\mathop{=}\limits_{(\ref{symeq2})}
M_5(\tfrac b{b-a}, 1, \tfrac{e-f}f, \tfrac d{d-c},\tfrac gh) \mathop{=}\limits_{(\ref{linksymeq2})}
M_5(\tfrac gh,\tfrac d{d-c}, \tfrac{e-f}f, 1, \tfrac b{b-a})\\ \mathop{=}\limits_{(\ref{linksymeq1})}
M_5(\tfrac b{b-a}, \tfrac gh,\tfrac d{d-c}, \tfrac{e-f}f, 1)
\end{gather*}
we will simply write
\[
M_5(\tfrac ab,\tfrac cd, \tfrac ef, \tfrac gh, \infty)\mathop{=}\limits_{(\ref{symeq2})}
M_5(\tfrac b{b-a}, \tfrac gh,\tfrac d{d-c}, \tfrac{e-f}f, 1).
\]

\begin{proposition}\label{M_3M_5:prop}
The following identities hold:
\begin{equation}\label{-1-1}
M_3\left(-1,-1\right) = P\times S^1 \bigb 0110
\end{equation}
\begin{equation}\label{-1-1fillings}
M_3\left(-1,-1,\tfrac ab\right) = \seifuno A{b}{b-a} \bigb 0110
\end{equation}
\begin{equation}\label{-2-2}
M_3(-1,-3,-3)=\seifdue{D}212{-1}
        \bigcup\nolimits_{{\tiny{\matr 121{-1}}}\phantom{\Big|}\!\!}
        \seifdue{D}2131\phantom{\Big|}
\end{equation}
\begin{equation}\label{-2-2-2}
M_3(-2,-2,-2)=\seifdue{D}212{-1}
        \bigcup\nolimits_{{\tiny{\matr {-1}31{-2}}}\phantom{\Big|}\!\!}
        \seifdue{D}2131
\end{equation}
All other exceptional $M_3(\alpha)$ can be expressed as some filling of 
$F$ or of $M_3(-1, -1)$ and the correspondence is found in Table \ref{Magic:ids}. 
\end{proposition}

\begin{table}[htbp]
\begin{center}
\begin{tabular}{|c|c|}
\hline
\phantom{\Big|}$M_3\big( \frac{p}{q}, \frac{r}{s}, \infty\big) = F\big(-\frac{1}{2}, \frac{q}{q-p}, \frac{1}{3}, \frac{r}{s}\big)$ 
\phantom{\Big|}& $M_3(\frac pq, \frac rs, 0\big) = F\big(\frac{s}{s-r}, 2, \frac{q}{3q-p}, -3\big)$ \\ \hline

$M_3\left( \frac{p}{q}, \frac{r}{s}, 1\big) = F\big(\frac{p-3q}{q}, -1, -2, \frac{r-2s}{s}\right)$ \phantom{\Big|} & 
$M_3\big( \frac{p}{q}, \frac{r}{s}, 2\big) = F\big(\frac{s}{r}, \frac{2q-p}{p-q}, -\frac{1}{2}, 3\big)$ \phantom{\Big|}\\ \hline

$M_3\big( \frac{p}{q}, \frac{r}{s}, 3\big) = F\big(-2, \frac{p-q}{q}, \frac{r-s}{s},-2\big)$ \phantom{\Big|}
& $M_3\big( \frac{p}{q}, \frac{3}{2}, \frac{5}{2}\big) = F\big(-3, \frac{p-2q}{p-q}, 2, -2\big)$  \phantom{\Big|}\\ \hline

$M_3\big(4, \frac{1}{2}, \frac{p}{q}\big) = F\big(2, \frac{3}{2}, \frac{q}{p-q}, -2\big)$ & 
$M_3\big(-1, \frac{1}{3}, 4\big) = M_3\big(\frac{3}{2}, \frac 32, \frac 83\big)$ \\ 
& $= F\big(-3, \frac 23, 2, -2 \big)$ \phantom{\Big|}\\ \hline

$M_3\big(\frac 52, \frac 53, \frac 53\big) = F\big(2, \frac 32, \frac 23, -2\big)$ & $M_3(-1,-2,-2) = F(\tfrac13,2,\tfrac14,-3)$ \phantom{\Big|} \\ \hline 

$M_3(-1,-2,-3) = F(-\tfrac43,-\tfrac13,2,-\tfrac12)$ & $M_3(-1,-2,-4) = F(-\tfrac32,-\tfrac12,3,-\tfrac12)$ \phantom{\Big|} \\ \hline

$M_3(-1,-2,-5) = F(-2,-2,-2,-3)$ & \phantom{\Big|} \\ \hline \hline

$M_3\big(\frac{3}{2}, \frac{7}{3}, \frac{7}{3}\big) = M_3\big(4, 4, \frac{2}{3}\big)$ & 
$M_3\big(5, 5, \frac{1}{2}\big) = M_3\big( -1, -1, \frac{1}{2}\big)$ \phantom{\Big|}\\
$= M_3\big(-1, -1, \frac{3}{2}\big)$ & \\ \hline

$M_3\big(\frac{5}{2}, \frac{5}{2}, \frac{4}{3}\big) = M_3\big(-1, -1, \frac{5}{2}\big)$
& $M_3\big(4, \frac 32, \frac 32\big) = M_3\big(-1, -1, 4\big)$ \phantom{\Big|}\\ \hline
\end{tabular}
\end{center}
\caption{\,\,\,\, Homeomorphisms between fillings of $M_3$ and of $F$ or $M_3(-1, -1)$.}
\label{Magic:ids}
\end{table}

\begin{proof} (\ref{-1-1})--(\ref{-2-2-2}) come directly from \cite{Magic}.
We prove the equalities in Table \ref{Magic:ids} using the identities in Lemma \ref{sym} 
as well as Identities (\ref{rolfsen:thm}), (\ref{F_from_M_5}), (\ref{blow:eqn}) and 
(\ref{figure8:eqn}). We indicate which Identity is being employed below the equality 
sign throughout the proof. 

We start with the case where $\alpha$ is an exceptional filling instruction on $M_3$ found in Theorem \ref{3CL_thm} 
and $M_3(\alpha)$ is homeomorphic to a filling of $F$:

\begin{equation} \label{M_3eq1} 
\begin{array}{c}
M_3\big( \frac{p}{q}, \frac{r}{s}, \infty\big) 
\mathop{=}\limits_{(\ref{rolfsen:thm})} M_5\big(-1, -2, \frac{p-q}{q}, \frac{r}{s}, \infty\big)  
\mathop{=}\limits_{(\ref{symeq2})}M_5\big(\frac{1}{2}, 1, \frac{p-2q}{p-q}, \frac{1}{3}, \frac{r}{s}\big)\\ 
\mathop{=}\limits_{(\ref{F_from_M_5})} F\big(-\frac{1}{2}, \frac{q}{q-p}, \frac{1}{3}, \frac{r}{s}\big) 
\end{array} 
\end{equation}

\begin{equation} \label{M_3eq2} \begin{array}{c} M_3\big(\frac pq, \frac rs, 0\big) \mathop{=}\limits_{(\ref{rolfsen:thm})} M_5(\frac{p-2q}{q}, -1, -2, \frac{r-s}{s}, 0) 
\mathop{=}\limits_{(\ref{symeq4})} M_5 (\frac{r-2s}{r-s}, 2, \frac{q}{3q-p}, -2, 1) \\ 
\mathop{=}\limits_{(\ref{F_from_M_5})} F\big(\frac{s}{s-r}, 2, \frac{q}{3q-p}, -3\big) 
\end{array} \end{equation}

\begin{equation} \label{M_3eq3} \begin{array}{c} M_3\big( \frac{p}{q}, \frac{r}{s}, 1\big) \mathop{=}\limits_{(\ref{rolfsen:thm})} 
M_5\big(\frac{p-2q}{q}, -1, -2, \frac{r-s}{s}, 1\big) 
\mathop{=}\limits_{(\ref{F_from_M_5})}F\big(\frac{p-3q}{q}, -1, -2, \frac{r-2s}{s}\big) 
\end{array} \end{equation}

\begin{equation} \label{M_3eq4} 
\begin{array}{c} M_3\big( \frac{p}{q}, \frac{r}{s}, 2\big) \mathop{=}\limits_{(\ref{rolfsen:thm})} M_5\big(0, -1, -2, \frac{p-q}{q}, \frac{r}{s}\big) 
\mathop{=}\limits_{(\ref{symeq6})} M_5\big(\frac{s}{r}, \frac{q}{p-q}, 1, \frac{1}{2}, 3\big) \\ 
\mathop{=}\limits_{(\ref{F_from_M_5})} F\big(\frac{s}{r}, \frac{2q-p}{p-q}, -\frac{1}{2}, 3\big) \end{array}\end{equation}

\begin{equation} \label{M_3eq5}
\begin{array}{l} M_3\big( \frac{p}{q}, \frac{r}{s}, 3\big)\mathop{=}\limits_{(\ref{rolfsen:thm})} M_5\big(-1, \frac{p-q}{q}, \frac{r-s}{s}, -1, 1\big) 
\mathop{=}\limits_{(\ref{F_from_M_5})} F\big(-2, -2, \frac{p-q}{q}, \frac{r-s}{s}\big) 
\end{array} \end{equation}

\begin{equation} \label{M_3eq6}
\begin{array}{l} \qquad M_3\big( \frac{p}{q}, \frac{3}{2}, \frac{5}{2}\big)
\mathop{=}\limits_{(\ref{rolfsen:thm})} M_5\big(\frac{1}{2}, -2, -1, \frac{1}{2}, \frac{p}{q}\big) 
\mathop{=}\limits_{(\ref{symeq9})} M_5\big(-1, \frac{1}{2}, \frac{p}{q}, \frac{1}{3}, -1\big) \\ 
\mathop{=}\limits_{(\ref{rolfsen:thm})} M_5\big(0, \frac{1}{2}, \frac{p-q}{q}, -1, \frac{1}{3}\big) \mathop{=}\limits_{(\ref{symeq2})} M_5\big(1, -2, \frac{p-2q}{p-q}, 2, -1\big) 
\mathop{=}\limits_{(\ref{F_from_M_5})} F\big(-3, \frac{p-2q}{p-q}, 2, -2\big) 
\end{array} \end{equation}

\begin{equation} \label{M_3eq7}
\begin{array}{l} \,\,\,\,\, M_3\big(4, \frac{1}{2}, \frac{p}{q}\big) 
\mathop{=}\limits_{(\ref{rolfsen:thm})}M_5\big(2, \frac{1}{2}, \frac{p-q}{q}, -2, -1\big) 
\mathop{=}\limits_{(\ref{symeq4})} M_5\big(\frac{1}{2}, \frac{p-q}{q}, -1, -1, \frac{3}{2}\big) \\ 
\mathop{=}\limits_{(\ref{rolfsen:thm})}
M_5\big(-\frac{1}{2}, -1, \frac{p-q}{q}, 0, \frac{3}{2}\big) 
\mathop{=}\limits_{(\ref{symeq3})}M_5\big(3, \frac{3}{2}, \frac{q}{p-q}, -1, 1\big) 
\mathop{=}\limits_{(\ref{F_from_M_5})} F\big(2, \frac{3}{2}, \frac{q}{p-q}, -2\big)
\end{array} \end{equation}

\begin{equation} \label{M_3eq8} \begin{array}{l}
M_3\big(-1, \frac{1}{3}, 4\big) \mathop{=}\limits_{(\ref{rolfsen:thm})}M_5\big(-2, -1, -\frac{5}{3}, -1, 3\big) 
\mathop{=}\limits_{(\ref{symeq8})}M_5\big(\frac 12, \frac 83, -\frac 12, -1, -2\big) \\ 
\mathop{=}\limits_{(\ref{rolfsen:thm})}M_3\big(\frac{3}{2}, \frac 32, \frac 83\big) 
\mathop{=}\limits_{(\ref{rolfsen:thm})} 
M_5\big(-1, \frac 23, -1, \frac 12, \frac 12\big) \mathop{=}
\limits_{(\ref{symeq3})} M_5\big(2, -1, \frac 12, \frac 32, -1\big) \\
\qquad \qquad \qquad \qquad \mathop{=}\limits_{(\ref{rolfsen:thm})}M_3\big(4, \frac 32, \frac 52\big) \mathop{=}\limits_{(\ref{M_3eq6})} 
F\big(-3, \frac 23, 2, -2\big)
\end{array}
\end{equation}

\begin{equation} \label{M_3eq9}
\begin{array}{l}
M_3\big(\frac 52, \frac 53, \frac 53\big) \mathop{=}\limits_{(\ref{rolfsen:thm})} 
M_5\big(\frac 32, -1, -\frac{1}{3}, -1, \frac 23\big) 
\mathop{=}\limits_{(\ref{symeq2})} M_5\big( -2, -\frac 12, 4, \frac 12, -1\big) \\ 
\qquad \qquad \qquad \qquad \mathop{=}\limits_{(\ref{rolfsen:thm})} M_3\big(\frac 12, 4, \frac 52\big) 
\mathop{=}\limits_{(\ref{M_3eq7})}
F\big(2, \frac 32, \frac 23, -2\big)
\end{array}
\end{equation}

\begin{equation} \label{M_3eq10}
\begin{array}{l}
M_3\big(-1, -2, -2\big) \mathop{=}\limits_{(\ref{rolfsen:thm})} M_5\big(-1, -3, -1, -3, -1\big) 
\mathop{=}\limits_{(\ref{blow:eqn})} M_5\big( -3, -1, -2, -3, 0\big) \\ 
\qquad \qquad \qquad \qquad \mathop{=}\limits_{(\ref{rolfsen:thm})} M_3\big(-1, -2, 0\big) \mathop{=}\limits_{(\ref{M_3eq2})}
F\big(\tfrac13, 2, \frac 14, -3\big)
\end{array}
\end{equation}

\begin{equation} \label{M_3eq11}
\begin{array}{c}
M_3\big(-1, -2, -3\big) \mathop{=}\limits_{(\ref{rolfsen:thm})} M_5\big(-1, -2, -2, -2, -5\big)\\ 
\mathop{=}\limits_{(\ref{figure8:eqn})} M_5\big(-1, -2, -2, -2, -1\big) 
\mathop{=}\limits_{(\ref{blow:eqn})} M_5\big( -2, -1, -3, -2, 0\big)\\  
\mathop{=}\limits_{(\ref{symeq1})} M_5\big(-\tfrac13, 1, \tfrac23,2,-\tfrac12\big)
\mathop{=}\limits_{(\ref{F_from_M_5})}F\big(-\tfrac43, -\tfrac13, 2, -\tfrac12\big)
\end{array}
\end{equation}

\begin{equation} \label{M_3eq12}
\begin{array}{l}
M_3\big(-1, -2, -4\big) \mathop{=}\limits_{(\ref{rolfsen:thm})} M_5\big(-1, -2, -2, -2, -6\big) 
\mathop{=}\limits_{(\ref{figure8:eqn})} M_5\big( -1, -2, -2, -2, 0\big) \\
\qquad \qquad \qquad \qquad \mathop{=}\limits_{(\ref{symeq1})} M_5\big(-\tfrac12, 1, \tfrac12, 3, -\tfrac12\big) 
\mathop{=}\limits_{(\ref{F_from_M_5})}
F\big(-\tfrac32, -\tfrac12, 3, -\tfrac12\big)
\end{array}
\end{equation}

\begin{equation} \label{M_3eq13}
\begin{array}{l}
M_3\big(-1, -2, -5\big) \mathop{=}\limits_{(\ref{rolfsen:thm})} M_5\big(-1, -2, -2, -2, -7\big) 
\mathop{=}\limits_{(\ref{figure8:eqn})} M_5\big( -1, -2, -2, -2, 1\big) \\
\qquad \qquad \qquad \qquad \mathop{=}\limits_{(\ref{F_from_M_5})}
F\big(-2, -2, -2, -3\big)
\end{array}
\end{equation}

We now turn to the case where $f$ is a filling instruction of $M_3$ and $M_3(f)$ is homeomorphic to a filling of $M_3(-1, -1)$:

\begin{equation} \label{M_3eq14}\begin{array}{l}
M_3\big(\frac{3}{2}, \frac{7}{3}, \frac{7}{3}\big) \mathop{=}\limits_{(\ref{rolfsen:thm})} M_5\big(-2, -1, \frac{1}{3}, \frac{7}{3}, \frac{1}{2}\big) 
\mathop{=}\limits_{(\ref{symeq3})} M_5\big(-1, -\frac{4}{3}, -1, 3, 3\big) 
\\ \mathop{=}\limits_{(\ref{rolfsen:thm})} M_3\big(4, 4, \frac{2}{3}\big) \mathop{=}\limits_{(\ref{rolfsen:thm})} M_5\big(2, \frac{2}{3}, 3, -2, -1\big) 
\mathop{=}\limits_{(\ref{symeq6})}M_5\big(-\frac{1}{2}, -1, -2, -2, -1\big) 
\\ \qquad \qquad \qquad \qquad \mathop{=}\limits_{(\ref{rolfsen:thm})} M_3\big(-1, -1, \frac{3}{2}\big)
\end{array} \end{equation}

\begin{equation} \label{M_3eq15}\begin{array}{l}
M_3\big(5, 5, \frac{1}{2}\big) \mathop{=}\limits_{(\ref{rolfsen:thm})} M_5\big(3, \frac{1}{2}, 4, -2, -1\big) 
\mathop{=}\limits_{(\ref{symeq6})} M_5\big(-\frac 12, -1, -3, -1, -2\big) \\ 
\qquad \qquad \qquad \qquad \mathop{=}\limits_{(\ref{rolfsen:thm})} M_3\big(-1, -1, \frac{1}{2}\big)
\end{array}
\end{equation}

\begin{equation} \label{M_3eq16}\begin{array}{l}
M_3\big(\frac{5}{2}, \frac{5}{2}, \frac{4}{3}\big) \mathop{=}\limits_{(\ref{rolfsen:thm})} M_5\big(\frac{1}{2}, -1, \frac{3}{2}, \frac{1}{3}, -1\big) 
\mathop{=}\limits_{(\ref{symeq9})} M_5\big(-2, -2, -1, \frac{1}{2}, -1\big) \\ 
\qquad \qquad \qquad \qquad \mathop{=}\limits_{(\ref{rolfsen:thm})} M_3\big(-1, -1, \frac{5}{2}\big)
\end{array}
\end{equation}

\begin{equation} \label{M_3eq17}\begin{array}{l}
M_3\big(4, \frac 32, \frac 32\big) \mathop{=}\limits_{(\ref{rolfsen:thm})} M_5\big(4, \frac 12, -2, -1, -\frac 12\big) 
\mathop{=}\limits_{(\ref{symeq6})} M_5\big(-1, -2, 3, -1, -3\big) \\ 
\qquad \qquad \qquad \qquad \mathop{=}\limits_{(\ref{rolfsen:thm})} M_3\big(4, -1, -1\big)
\end{array}
\end{equation}

This completes the proof.
\end{proof}

\begin{remark}
Identities (\ref{seifeq1})--(\ref{grapheq6}) can be used with Proposition \ref{fillingF} to show that 
Proposition \ref{M_3M_5:prop} is consistent with classification of exceptional fillings of the mirror 
of 3CL given in \cite{Magic}.
\end{remark}

Consideration of the ``intersection index" defined in Section 3.3\ of \cite{thesis} together with 
the classification of exceptional fillings of $M_3$ given in \cite{Magic} demonstrates that 
$(-1,-3,-2,-2,-3)$, $(-2,-\frac12,3,3,-\frac12)$, $(-2,-2,-2,-2,-2)$ do not factor through $M_3$. 
This leads to the following useful fact that will be used liberally throughout the 
proof of Theorems \ref{5CL_thm} and \ref{4CL_thm}.
By Theorem \ref{main:teo} we see that an exceptional filling instruction $\alpha$ on $M_5$ with no $0,1,\infty$ slopes 
factors through $M_3$ if and only if there is a $\gamma\subset \alpha$ with $\gamma\in [\![(-1, -2)]\!]$. 
By Lemma \ref{sym} and (\ref{blow:eqn}) we have
\begin{equation}\label{factor}
\begin{array}{c}
[\![((-1)_1, (-2)_2)]\!]=\bigl\{ 
[((\tfrac 12)_1, (\tfrac 23)_3)], [((\tfrac 12)_1, 3_2)], 
[((\tfrac 12)_1, (\tfrac 32)_2)], [((\tfrac 23)_1, 2_2)],\\ 
\,[ ( (\tfrac12)_1, (\tfrac13)_3 )],
 [(2_1, (-\tfrac 12)_3)], [ ( (\tfrac 13)_1, 2_2)],  [((-1)_1, 3_3)],[((-1)_1, (-2)_2)],\\ 
\, [((-1)_1, (\tfrac 32)_3)], [(2_1, (-2)_3)],  [((-1)_1, (-\tfrac 12)_2)], [(-1_1,-1_3)], [(-1_1,2_2)], \\
\, [(-1_1,(\tfrac12)_2)], [(-1_1,(\tfrac12)_3)], [(2_1,2_2)], [(2_1,2_3)], 
[(2_1,(\tfrac12)_3], [((\tfrac12)_1,(\tfrac12)_2)]\bigr\}
\end{array}
\end{equation}

There is an obvious strategy to prove Theorem \ref{5CL_thm}. 
Namely, we examine all hyperbolic $M_5(\alpha_1,\alpha_2,\alpha_3,\alpha_4,\varnothing)$ and look 
for all $\alpha_i$ so that $\alpha=(\alpha_1,\alpha_2,\alpha_3,\alpha_4,\alpha_5)$ contains an isolated 
exceptional filling instruction. This simplifies matters greatly. If $\alpha$ factors through 
$M_3$ then $\alpha_5$ must be a slope in $[\![(-1)]\!]$ or $[\![(-2)]\!]$, and if $\alpha$ does not factor 
through $M_3$ then $\alpha$ contains a slope in $[\![(1)]\!]$ or $\alpha$ is equivalent to one of 
$(-2,-\tfrac12,3,3,-\tfrac12)$, $(-1,-3,-2,-2,-3)$, $(-2,-2,-2,-2,-2)$. 
We now prove Theorems \ref{5CL_thm} and \ref{4CL_thm}.

\subsubsection{Proof of Theorem \ref{5CL_thm}.} We start by establishing (\ref{5CLinf})-(\ref{5CL0}). 
Identity (\ref{5CLinf}) is established in \cite{MPR}, and (\ref{5CL1}) is exactly the same as 
(\ref{F_from_M_5}). For (\ref{5CL0}) we have
\begin{gather*}
M_5(\tfrac ab,\tfrac cd, \tfrac ef, \tfrac gh, 0)\mathop{=}\limits_{(\ref{symeq1})}
M_5(\tfrac fe,1, \tfrac a{a-b}, \tfrac{d-c}d, \tfrac hg)\mathop{=}\limits_{(\ref{F_from_M_5})}
F(\tfrac{f-e}e,\tfrac b{a-b},\tfrac{d-c}d,\tfrac hg) 
\mathop{=}\limits_{(\ref{Fsym})} F(\tfrac b{a-b},\tfrac{d-c}d,\tfrac hg, \tfrac{f-e}e)\\ 
\mathop{=}\limits_{(\ref{grapheq1})\,\, \&\,\, (\ref{grapheq6})} 
F(\tfrac b{b-a},\tfrac{c-d}d,-\tfrac hg, \tfrac{e-f}e)
\mathop{=}\limits_{(\ref{seifeq3})} 
F(\tfrac b{b-a},\tfrac{c-d}c,-\tfrac hg, \tfrac{e-f}f). 
\end{gather*}

We let $\alpha$ be a hyperbolic surgery 
instruction on 5CL containing at least one $\varnothing$ not factoring through $M_4$. We 
let $\tau$ be a boundary component of $M_5(\alpha)$. Lemma \ref{sym} 
allows us to assume that $\tau$ comes from the 5th component of 5CL. So, we assume that 
$\alpha=(\alpha_1,\alpha_2, \alpha_3,\alpha_4,\varnothing)$ with $\alpha_i\in\mathbb{Q}\cup\{\varnothing,\infty\}$. 

Our assumption that $\alpha$ is hyperbolic and does not factor through $M_4$ imposes restrictions on the $\alpha_i$. 
Theorem \ref{main:teo} and the remark that follows tells us that $\alpha$ being hyperbolic means no 
$\gamma\subseteq\alpha$ 
contains an isolated exceptional filling instructions on $M_5$. Identity (\ref{5CL_twist}) means that if $\alpha$ does 
not factor through $M_4$ then no $\alpha_i=-1$. Lemma \ref{sym}  tells us that 
\begin{equation}\label{image-1} 
[\![(-1)]\!]=[(-1)]\sqcup[(\tfrac12)]\sqcup[(2)].
\end{equation}
Thus, no $\alpha_i\in\{-1,\frac12,2\}$ if $\alpha$ does not factor through $M_4$.

We now examine $E_{\tau}(M_5(\alpha))$. Theorem \ref{main:teo} implies that all slopes in $[\![(1)]\!]$ are exceptional. 
By Lemma \ref{sym} 
\begin{equation}\label{image1}
[\![(1)]\!]=[(1)]\sqcup[(\infty)]\sqcup[(0)]
\end{equation} 
Therefore, (\ref{image1}) implies that no $\alpha_i\in\{0,1,\infty\}$, that 
$\{0,1,\infty\}\subseteq E_{\tau}(M_5(\alpha))$ and $e(M_5(\alpha)) \geq 3$. 

We will now describe all such $\alpha$s not factoring 
through $M_4$ with $e(M_5(\alpha))>3$. We define $(\alpha, \beta)$ to be $(\alpha_1, \alpha_2, \alpha_3, \alpha_4,\beta)$. 
If $\beta$ is an exceptional slope of $M_5(\alpha)$ then $(\alpha, \beta)$ contains an isolated exceptional filling instruction. 
Lemma \ref{sym} tells us that 
\begin{equation}\label{image-2}
[\![(-2)]\!]=[(-2)]\sqcup[(-\tfrac12)]\sqcup[(\tfrac13)]\sqcup[(\tfrac23)]\sqcup[(\tfrac32)]\sqcup[(3)].
\end{equation} 
By (\ref{image-1}) and Theorem \ref{main:teo} with (\ref{image-2}), any such isolated exceptional 
filling instruction contains at most one slope in $[\![(-1)]\!]$, and contains two slopes in $[\![(-2)]\!]$. 
Thus at least one of the slopes in $\alpha$ belongs to $[\![(-2)]\!]$. It is a 
routine consequence of Lemma \ref{sym} that we may assume without loss of generality that 
$\alpha_1=-2$ and that $\tau$ remains as the 5th component of $\partial M_5$. 

Finally, the remarks following Theorem \ref{main:teo} allow us to conclude that if $(\alpha, \beta)$ 
contains an isolated exceptional filling instruction then either $\beta$ is a slope in $[\![(1)]\!]$, 
or $\beta$ is a slope in $[\![(-1)]\!]$ and $(\alpha, \beta)$ factors through $M_3$, or 
$M_5(\alpha, \beta)$ is one of 
\begin{gather*}
m_1:=M_5(-2,-\tfrac12,3,3,-\tfrac12),\qquad m_2:=M_5(-1,-3,-2,-2,-3),
\\ m_3:=M_5(-2,-2,-2,-2,-2).
\end{gather*}
 
So, we define the following sets of $(\alpha, \beta)$s;
\begin{itemize}
\item We define $l$ to be the set of exceptional $(\alpha, \beta)$s such that $M_5(\alpha, \beta)$ is in $\{m_1, m_2, m_3\}$;
\item We define $l_{-1}$ to be the set of exceptional $(\alpha, \beta)$s factoring through $M_3$ with $\beta=-1$;
\item We define $l_{\frac 12}$ to be the set of exceptional $(\alpha, \beta)$s factoring through $M_3$ with $\beta=\frac 12$;
\item We define $l_{2}$ to be the set of exceptional $(\alpha, \beta)$s factoring through $M_3$ with $\beta=2$;
\end{itemize}

Let $p:l\cup l_{-1}\cup l_{\frac 12}\cup l_{2}\rightarrow A$ be defined by $(\alpha, \beta) \mapsto \alpha$. 
For $\alpha$ in $p(l\cup l_{-1}\cup l_{\frac 12}\cup l_{2})$ define $B_{\alpha}$ to be the set of all $\beta$'s 
such that $(\alpha, \beta)$ is contained in $l\cup l_{-1}\cup l_{\frac 12}\cup l_2$. It is clear that $E\big(M_5(\alpha)\big)=
\{0, 1, \infty\}\cup B_{\alpha}$ and that $\{\big(M_5(\alpha), E(M_5(\alpha))\big)\}_{\alpha}$ 
is a complete list of all $\big(M_5(\alpha), E(M_5(\alpha))\big)$ pairs with $\alpha=
\big(-2, \frac pq, \frac rs, \frac uv, \varnothing\big)$, $M_5(\alpha)$ hyperbolic and $e\big(M_5(\alpha)\big)>3$.

We now explicitly construct the sets $l, l_{-1}, l_{\frac 12}, l_2$. Throughout these constructions 
\begin{equation}\label{alf}
\alpha_i\not\in\{\infty,0,1,-1,\tfrac12,2\}
\end{equation} 
is used to reduce the number of possible cases (see (\ref{image-1}) and (\ref{image1})) .

\textbf{Construction of the set $l$:} We can see using Lemma \ref{sym} that every 
$(\alpha, \beta)$ with $M_5(\alpha, \beta)=m_i$ with $\alpha_1=-2$ and $\alpha$ 
hyperbolic not factoring through $M_4$ is contained in the following set
\begin{equation*}
\begin{array}{l} 
l=\bigl\{\big( -2, \tfrac 14, \tfrac 32, \tfrac 43, \tfrac 12 \big) , \big(-2, -\tfrac 12, 3, 3, -\tfrac 12\big), \big(-2, \tfrac 13, 3, \tfrac 13, -2\big), 
\big(-2, -2, \tfrac 13, 3, \tfrac 13\big),\\ 
\big(-2, -\tfrac 12, -2, \tfrac 32, \tfrac 32\big), \big(-2, \tfrac 32, \tfrac 32, -2, -\tfrac 12\big), \big(-2, -2, -2, -2, -2\big), 
\big(-2, \tfrac 13, \tfrac 32, \tfrac 32, \tfrac 13\big)\bigr\}.
\end{array}
\end{equation*}

Before constructing the sets $l_{-1}, l_{\frac 12}, l_2$ we recall that  
$(\alpha, \beta)$ is not in $l$ then $(\alpha, \beta)$ factors through $M_3$ if and only 
if $(\alpha, \beta)$ contains one of the elements of $[\![((-1)_1, (-2)_2)]\!]$ (see 
(\ref{factor}) and the preceding remarks). We will see that the requirements of $\alpha$ 
being hyperbolic and not factoring through $M_4$ (see (\ref{alf})) allow us to completely construct 
$l_{-1}, l_{\frac 12}, l_2$. 

\textbf{Construction of the set $l_{-1}$:} In this case, for $\alpha=(-2, \frac pq, \frac rs, \frac uv)$, we have 
$M_5(\alpha)(-1)=M_3(\frac{p+q}q, \frac rs, \frac {u+2v}v)$ by (\ref{rolfsen:thm}). So 
$(\alpha, \beta)$ is an exceptional filling of $M_5$ if and only if $(\frac{p+q}q, \frac rs, \frac {u+2v}v)$ 
contains an isolated exceptional filling instruction on $M_3$. As noted, (\ref{alf}) tells us that if 
$\alpha$ is hyperbolic and does not factoring through $M_4$ then 
$\frac{p+q}q\not\in\{\infty, 2, 1, 0, \frac 32, 3\}$, $\frac rs \not \in \{\infty, 1, 0, -1, \frac 12, 2\}$, 
$\frac{u+2v}v\not\in \{\infty, 3, 2, 1, \frac 52, 4\}$. With these conditions and Theorem \ref{3CL_thm} it 
is easy to see that $\beta$ is an exceptional slope on a hyperbolic $M_5(\alpha)$ if and only if 
one of the following holds:
\begin{itemize}
\item $\frac rs =3$ \\
\item $\frac uv+2=0$ \\
\item $(\frac pq+1, \frac uv+2)$ belongs to $\{(-1, -1), (\frac 52, \frac 32), (4, \frac 12)\}$ \\
\item $(\frac rs, \frac pq +1)$ belongs to $\{(\frac 32, \frac 52), (4, \frac 12)\}$ \\
\item $(\frac rs, \frac uv+2)$ belongs to $\{(\frac 52, \frac 32), (4, \frac 12)\}$ \\
\item $(\frac pq +1, \frac rs, \frac uv +2)$ belongs to 
\begin{gather*}
\Big\{(5, 5, \tfrac 12), (\tfrac 12, 5, 5),(4, 4, \tfrac 23), (4, \tfrac 32,  \tfrac 32), (4,  \tfrac13, -1), (\tfrac 13, 4, -1), 
(-1, 4, \tfrac 13), (-1, \tfrac13, 4), (\tfrac 83, \tfrac 32, \tfrac 32), \\
(\tfrac 52, \tfrac 52, \tfrac 43), (\tfrac 52, \tfrac 53, \tfrac 53), (\tfrac 53, \tfrac 52, \tfrac 53), (\tfrac 73, \tfrac 73, \tfrac 32), 
(\tfrac 73, \tfrac 32, \tfrac 73), (-1, -2, -2), (-2, -2, -1), (-1, -2, -3),\\ 
(-3, -2, -1), (-2, -3, -1), (-1, -3, -2), (-1, -2, -4), (-4, -2, -1),\\ 
(-1, -4, -2), (-2, -4, -1), (-1, -2, -5), (-5, -2, -1), (-1, -5, -2), (-2, -5, -1), \\
(-1, -3, -3), (-3, -3, -1), (-2, -2, -2))\Big\}. 
\end{gather*}
\end{itemize}

Thus, the set of all $(-2, \frac pq, \frac rs, \frac uv, -1)$ constructed in the above analysis is the 
set $l_{-1}$.

\textbf{Construction of the set $l_{2}$:} 
By (\ref{factor}) for $(\alpha, \beta)=(-2, \frac pq, \frac rs, \frac uv, 2)$ to factor through 
$M_3$ we need one of $\frac pq, \frac rs$ to be in $\{-2, -\frac 12\}$ or $\frac uv \in \{\frac 13, \frac 23\}$. 
We construct all $(\alpha, \beta)$ for each of these 6 cases individually. The number of cases is controlled 
by (\ref{alf}).

If $\frac pq = -2$ then 
\[M_5(\alpha)(\beta)\mathop{=}\limits_{(\ref{symeq11})} M_5\big(\tfrac 32, \tfrac{r-s}r, \tfrac{v}{v-u}, -2, -1\big)
\mathop{=}\limits_{(\ref{rolfsen:thm})} M_3\big(\tfrac 72, \tfrac{r-s}r, \tfrac{2v-u}{v-u}\big).\]
By (\ref{alf}) we have 
$\frac {r-s}{r}\not\in\{\infty, 0, 1, 2, \frac 12, -1\}$ and 
$\frac{2v-u}{v-u}\not\in\{2, \infty, 1, \frac 32, 0, 3\}$. 
From Theorem \ref{3CL_thm} and (\ref{alf}) we see that $\beta=2$ is an exceptional slope on $M_5(\alpha)$ 
if and only if one of the following holds;
$\frac{r-s}{r}=3$, $(\frac{r-s}{r}, \frac{2v-u}{v-u})=(\frac 32, \frac 52)$, 
$(\frac{r-s}{r}, \frac{2v-u}{v-u})=(4, \frac 12)$.

If $\frac pq =-\frac 12$ then 
\[ M_5(\alpha)(\beta)\mathop{=}\limits_{(\ref{symeq8})} M_5\big(\tfrac u{u-v}, \tfrac {s-r}s, -\tfrac 12, -2, -1\big)
\mathop{=}\limits_{(\ref{rolfsen:thm})} M_3\big(\tfrac{3u-2v}{u-v}, \tfrac{s-r}s, \tfrac 12\big).\]
From Theorem \ref{3CL_thm} and (\ref{alf}) we see that $\beta=2$ is an exceptional slope on $M_5(\alpha)$ 
if and only if one of the following holds;
$\frac{s-r}s=3$, $\frac{3u-2v}{u-v}=0$, $\frac{s-r}s = 4$, $(\frac{s-r}s, 
\frac{3u-2v}{u-v})=(\frac 52, \frac 32)$, $(\frac{s-r}s, 
\frac{3u-2v}{u-v})=(5, 5)$.

If $\frac rs = -2$ then 
\[M_5(\alpha)(\beta)\mathop{=}\limits_{(\ref{symeq4})} M_5\big(-1, -2, \tfrac 13, \tfrac {p-q}p, \tfrac{u-v}u\big)
\mathop{=}\limits_{(\ref{rolfsen:thm})} M_3\big(\tfrac 43, \tfrac{p-q}p, \tfrac{3u-v}u\big).\]
From Theorem \ref{3CL_thm} and (\ref{alf}) we see that $\beta=2$ is an exceptional slope on $M_5(\alpha)$ 
if and only if one of the following holds; 
$\frac{p-q}p=3$, $\frac{3u-v}u=0$, $(\frac{p-q}p, \frac{3u-v}u)=(4, \frac 12)$, 
$(\frac{p-q}p, \frac{3u-v}u)=(\frac 52, \frac 32)$.

If $\frac rs = -\frac 12$ then 
\[M_5(\alpha)(\beta) \mathop{=}\limits_{(\ref{symeq1})} M_5\big(-2, -1, \tfrac 23, \tfrac{q-p}q, \tfrac vu\big) 
\mathop{=}\limits_{(\ref{rolfsen:thm})} M_3\big(\tfrac 83, \tfrac{q-p}q, \tfrac{v+u}u\big).\]
From Theorem \ref{3CL_thm} and (\ref{alf}) we see that $\beta=2$ is an exceptional slope on $M_5(\alpha)$ 
if and only if one of the following holds;
$\frac{q-p}q=3$, $(\frac{q-p}q, \frac{v+u}u) = (4, \frac 12)$, $(\frac{q-p}q, \frac{v+u}u) = (\frac 32, \frac 52)$.

If $\frac uv = \frac 13$ then 
\[M_5(\alpha)(\beta) \mathop{=}\limits_{(\ref{symeq4})} M_5\big(-1, \tfrac rs, \tfrac 13, \tfrac{p-q}p, -2\big) 
\mathop{=}\limits_{(\ref{rolfsen:thm})} M_3\big(\tfrac{r+2s}{s}, \tfrac 13, \tfrac{2p-q}p\big).\]
From Theorem \ref{3CL_thm} and (\ref{alf}) we see that $\beta=2$ is an exceptional slope on $M_5(\alpha)$ 
if and only if one of the following holds;
$\frac{r+2s}s=0$, $(\frac{r+2s}{s}, \frac{2p-q}p)=(\frac 12, 4)$, 
$(\frac{r+2s}{s}, \frac{2p-q}p)= (\frac 32, \frac 52)$, 
$(\frac{r+2s}{s}, \frac{2p-q}p)=(-1, -1)$, $(\frac{r+2s}{s}, \frac{2p-q}p)=(-1, 4)$.

If $\frac uv = \frac 23$ then 
\[M_5(\alpha)(\beta) \mathop{=}\limits_{(\ref{symeq8})} M_5\big(-2, \tfrac{s-r}s, -\tfrac 12, \tfrac qp, -1\big)
\mathop{=}\limits_{(\ref{rolfsen:thm})} M_3\big(\tfrac{2s-r}s, -\tfrac 12, \tfrac{2p+q}p\big).\]
From Theorem \ref{3CL_thm} and (\ref{alf}) we see that $\beta=2$ is an exceptional slope on $M_5(\alpha)$ 
if and only if one of the following holds;
$\frac{2p+q}p=0$, $(\frac{2s-r}s, \frac{2p+q}p)=(-1, -1)$, $(\frac{2s-r}s, \frac{2p+q}p)=(4, \frac 12)$, 
$(\frac{2s-r}s, \frac{2q+p}p)=(\frac 52, \frac 32)$.

Thus, the set of all $(-2, \frac pq, \frac rs, \frac uv, 2)$ constructed in the above analysis is the set $l_2$.

\textbf{Construction of the set $l_{\frac 12}$:} 
Reasoning as in the case $\beta=2$ we see $(\alpha, \beta)$ factors through 
$M_3$ when one of $\frac pq, \frac rs$ is in $\{\frac 13, \frac 23\}$ or $\frac uv \in \{\frac 32, 3\}$. 
We examine each of these 6 cases individually and enumerate all $(-2, \frac pq, \frac rs, \frac uv, \frac 12)$ 
that satisfy (\ref{alf}).

If $\frac pq =\frac 13$ then 
\[M_5(\alpha)(\beta)\mathop{=}\limits_{(\ref{symeq7})} M_5\big(\tfrac{v}{v-u}, -2, \tfrac{s}{s-r}, -2, -1\big) 
\mathop{=}\limits_{(\ref{rolfsen:thm})} M_3\big(\tfrac{3v-2u}{v-u}-2, \tfrac{2s-r}{s-r}\big).\]

Theorem \ref{3CL_thm} and (\ref{alf}) tell us that $\beta=\frac 12$ is an exceptional slope on 
$M_5\big(-2, \frac pq, \frac rs, \frac uv\big)$ 
if and only if one of the following holds; 
$\frac{3v-2u}{v-u}=0$, $(\frac{3v-2u}{v-u}, \frac{2s-r}{s-r})=(-1, -1)$, 
$(\frac 12, 4)$, $(\frac 32, \frac 52)$, $(-1,-2)$, $(-2,-1)$, $(-1,-3)$, $(-3,-1)$, $(-1,-4)$, 
$(-4,-1)$, $(-1,-5)$, $(-5,-1)$, $(-2,-2)$, $(-3,-3)$.

If $\frac pq =\frac 23$ then 
\[M_5(\alpha)(\beta)\mathop{=}\limits_{(\ref{symeq5})} M_5\big(\tfrac 23, \tfrac r{r-s}, -1, -2, \tfrac u{u-v}\big)
\mathop{=}\limits_{(\ref{rolfsen:thm})} M_3\big(\tfrac 23, \tfrac{3r-2s}{r-s}, \tfrac{2u-v}{u-v}\big).\]
Theorem \ref{3CL_thm} and (\ref{alf}) tell us that $\beta=\frac 12$ is an exceptional slope on 
$M_5\big(-2, \frac pq, \frac rs, \frac uv\big)$ 
if and only if one of the following holds;
$\frac{3r-2s}{r-s}=0$, $(\frac{3r-2s}{r-s}, \frac{2u-v}{u-v})=(-1, -1)$, $(\frac 32, \frac 52)$, 
$(\tfrac12,4)$.

If $\frac rs=\frac 13$ then 
\[M_5(\alpha)(\beta)
\mathop{=}\limits_{(\ref{symeq2})}M_5\big(\tfrac 13, -1, -2, \tfrac q{q-p}, \tfrac uv\big)
\mathop{=}\limits_{(\ref{rolfsen:thm})} M_3\big(\tfrac 73, \tfrac{2q-p}{q-p}, \tfrac uv\big).\]
Theorem \ref{3CL_thm} and (\ref{alf}) tell us that $\beta=\frac 12$ 
is an exceptional slope on $M_5\big(-2, \frac pq, \frac rs, \frac uv\big)$ 
if and only if one of the following holds; $\frac uv=3$, $(\frac uv, \frac{2q-p}{q-p})$=
$(\frac 32, \frac 52)$, $(\frac 32, \frac 73)$, $(4,\tfrac12)$.

If $\frac rs =\frac 23$ then 
\[M_5(\alpha)(\beta)\mathop{=}\limits_{(\ref{symeq5})} M_5\big(\tfrac 23, -2, -1, \tfrac p{p-q}, \tfrac u{u-v}\big)
\mathop{=}\limits_{(\ref{rolfsen:thm})} M_3\big(\tfrac 53, \tfrac{3p-2q}{p-q}, \tfrac u{u-v}\big).\]
Theorem \ref{3CL_thm} and (\ref{alf}) tell us that $\beta=\frac 12$ 
is an exceptional slope on $M_5\big(-2, \frac pq, \frac rs, \frac uv\big)$ 
if and only if one of the following holds; 
$\frac{3p-2q}{p-q}=0$, $\frac{u}{u-v}=3$, $(\frac{3p-2q}{p-q}, \frac{u}{u-v})=(\frac 32, \frac 52)$, 
$(\frac 12, 4)$, $(\frac 53, \frac 52)$.

If $\frac uv = \frac 32$ then 
\[M_5(\alpha)(\beta)\mathop{=}\limits_{(\ref{symeq7})}
M_5\big(-2, -2, \tfrac{s}{s-r}, \tfrac{p-q}p, -1\big) 
\mathop{=}\limits_{(\ref{rolfsen:thm})} M_3\big(-1, \tfrac{s}{s-r}, \tfrac{3p-q}p \big).\]
Theorem \ref{3CL_thm} and (\ref{alf}) tell us that $\beta=\frac 12$ 
is an exceptional slope on $M_5\big(-2, \frac pq, \frac rs, \frac uv\big)$ 
if and only if one of the following holds; $\frac s{s-r}=3$, $\frac{3p-q}p=0, -1$, 
$(\frac{s}{s-r}, \frac{3p-q}p\big)=(\frac 52, \frac 32)$, $(4, \frac 12)$, 
$(-2, -2)$, $(-2, -3)$, $(-3, -2)$, $(-2, -4)$, $(-4,-2)$, $(-2, -5)$, $(-5, -2)$, $(4, \frac 13)$, $(-3,-3)$.

If $\frac uv = 3$ then 
\[M_5(\alpha)(\beta)\mathop{=}\limits_{(\ref{symeq3})}
M_5\big(-1, 3, \tfrac sr, \tfrac qp, -2\big) \mathop{=}\limits_{(\ref{rolfsen:thm})}
 M_3\big(5, \tfrac sr, \tfrac {q+p}p\big).\]
Theorem \ref{3CL_thm} and (\ref{alf}) tell us that $\beta=\frac 12$ 
is an exceptional slope on $M_5\big(-2, \frac pq, \frac rs, \frac uv\big)$ 
if and only if one of the following holds; 
$\frac sr =3$, $\big(\frac sr, \frac{q+p}p\big)=\big(\frac 32, \frac 52\big)$, $(4, \frac 12)$,  
$(5, \frac 12)$.

Thus, the set of all $(-2, \frac pq, \frac rs, \frac uv, \frac 12)$ constructed in the 
above analysis is the set $l_{\frac 12}$.

This completes the construction of $l\cup l_{-1}\cup l_{\frac 12}\cup l_2$ and the sets 
$\{(M_5(\alpha), E(M_5(\alpha))\}$ are now easily computed. The last step is to reduce 
the size of $\{(M_5(\alpha), E(M_5(\alpha))\}$ using (\ref{first:eqn})--(\ref{figure8:eqn}). 
Namely, only one $(M_5(\alpha), E(M_5(\alpha))$ is shown for each $[\![\alpha]\!]$. This is 
done using Lemma \ref{sym} with the help of an ad hoc Python script \cite{code}. 
%The lists of filling instructions before and after using the \cite{code} are available at \cite{lists}. 

The reduced list of fillings is shown in Tables \ref{tb1}--\ref{tb4}. This completes the 
proof of Theorem \ref{5CL_thm}. $\square$

\vspace{5mm}

The elementary techniques used to prove Theorem \ref{5CL_thm} can obviously be applied to 
describe all $E(M_5(\alpha))$ when $\alpha$ factors through $M_4$ but $\alpha$ does not 
factor through $M_3$.

\vspace{5mm}

\subsubsection{Proof of Theorem \ref{4CL_thm}. } We first establish 
(\ref{4CLinf})-(\ref{4CL0}).\\ For (\ref{4CLinf}) we have 
\[
M_4(\tfrac ab,\tfrac cd,\tfrac ef, \infty) \mathop{=}\limits_{(\ref{5CL_twist})} 
M_5(\tfrac{a-b}b,-1,\tfrac{c-d}d,\tfrac ef,\infty) \mathop{=}\limits_{(\ref{5CLinf})}
F(\tfrac{b-a}b,\tfrac d{c-d},-1,-\tfrac ef). 
\]
For (\ref{4CL2}) we have 
\[
M_4(\tfrac ab,\tfrac cd,\tfrac ef, 2) \mathop{=}\limits_{(\ref{5CL_twist})}
M_5(-1, \tfrac{a-b}b,\tfrac cd,\tfrac ef,1) \mathop{=}\limits_{(\ref{F_from_M_5})}
F(-2, \tfrac{a-b}b,\tfrac cd,\tfrac{e-f}f) \mathop{=}\limits_{(\ref{Fsym})} 
F(\tfrac{a-b}b,\tfrac cd,\tfrac{e-f}f, -2). 
\]
For (\ref{4CL1}) we have 
\[
M_4(\tfrac ab,\tfrac cd,\tfrac ef, 1) \mathop{=}\limits_{(\ref{5CL_twist})}
M_5(\tfrac{a-b}b,-1, \tfrac{c-d}d,\tfrac ef,1) \mathop{=}\limits_{(\ref{F_from_M_5})}
F(\tfrac{a-2b}b,-1,\tfrac{c-d}d,\tfrac{e-f}f). 
\]
For (\ref{4CL0}) we have 
\[
M_4(\tfrac ab,\tfrac cd,\tfrac ef, 0) \mathop{=}\limits_{(\ref{5CL_twist})}
M_5(\tfrac{a-b}b,-1, \tfrac{c-d}d,\tfrac ef,0) \mathop{=}\limits_{(\ref{symeq1})}
M_5(\tfrac{d}{c-d},1,\tfrac{a-b}{a-2b},2, \tfrac fe) \mathop{=}\limits_{(\ref{F_from_M_5})}
F(\tfrac{2d-c}{c-d},\tfrac{b}{a-2b},2, \tfrac fe).
\]

We now examine the set of exceptional slopes on fillings of $M_4$ not factoring through $M_3$. 
We let $\alpha$ be a hyperbolic surgery instruction on 4CL containing at least one $\varnothing$ not factoring through $M_3$. 
We know from $(\ref{5CL_twist})$ that $M_5(-1)=M_4$. So we let 
$\alpha'=(-1, \alpha'_2, \alpha'_3, \alpha'_4, \varnothing)$ be such that $M_4(\alpha)=M_5(\alpha')$.
The argument can now proceed exactly as in the proof of Theorem \ref{5CL_thm} to enumerate 
the $E(M_5(\alpha'))$. 

Theorem \ref{5CL_thm} implies that every $\beta'\in\{-1, 0, 1, \infty\}$ is an exceptional slope on 
$M_5(\alpha')$. So, $e(M_5(\alpha'))\geq4$ and, with respect to the choice of 
basis induced from $M_5$, $\{-1,0,1,\infty\}\subseteq E(M_5(\alpha'))$. 

As in the proof of Theorem \ref{5CL_thm}, (\ref{factor}) imposes restrictions on the $\alpha'_i$. 
The condition that $M_5(\alpha')$ is hyperbolic with $\alpha$ not factoring through $M_3$ means that 
no instruction properly contained in $\alpha$ contains an instruction in $[\![(1)]\!]$ (see Theorem 
\ref{5CL_thm}) or $[\![(-1_1,-2_2)]\!]$ (see (\ref{factor})) or $[\![(-1_1,-1_2)]\!]$ (see Table 1.1.4 
from \cite{thesis}). This imposes the restrictions 
\begin{equation}\label{slope_restriction4CL}
\alpha'_2\not\in\{-2,-1, -\tfrac12, 0, \tfrac12, 1, 2, \infty\}
\qquad\text{and}\qquad 
\alpha'_3, \alpha'_4\not\in\{-1, 0, \tfrac12, 1, \tfrac32, 2, 3,\infty\}.
\end{equation}

If $\beta'$ is an exceptional slope on $M_5(\alpha')$, then $(\alpha',\beta')$ contains an 
isolated filling instruction. By Theorem \ref{5CL_thm}, $\beta'\in\{-1, 0, 1, \infty\}$, or 
$(\alpha', \beta')$ factors through $M_3$, or $(\alpha',\beta')$ is equivalent to $(-1,-3,-2,-2,-3)$. 
This implies that $\beta'$ is in one of $[\![(1)]\!]$, $[\![(-1)]\!]$, $[\![(-2)]\!]$ or, by Lemma 
\ref{sym}, that $(\alpha', \beta')$ is one of 
$(-1, -3, -2, -2, -3)$, $(-1, \frac 13, \frac 43, \frac 43, \frac 13)$, 
$(-1, -\frac 13, 4, \frac 23, 3)$, $(-1, 3, \frac 23, 4, -\frac 13)$. 

If $(\alpha', \beta')$ factors through $M_3$, then (\ref{factor}) in conjunction with (\ref{slope_restriction4CL}) 
tells us that $\beta'$ is in $\{-2,-\tfrac12,\tfrac12,2\}$. Every slope in $\{-2,-\tfrac12,\tfrac12,2\}$ 
is examined individually as in the proof of Theorem \ref{5CL_thm} to obtain a complete list of all 
$\big(M_5(\alpha'), E(M_5(\alpha'))\big)$ pairs that have 
$\alpha'=\big(-1, \frac pq, \frac rs, \frac uv, \varnothing\big)$, $M_5(\alpha')$ hyperbolic with $\alpha'$ 
not factoring through $M_3$ and $e\big(M_5(\alpha')\big)>4$. The result of the enumeration is that $\beta'=-2$ is 
an exceptional slope on $M_5(\alpha')$ if and only if $\alpha'$ is one of 
\begin{gather*}
(-1,-\tfrac32,4,\tfrac uv), (-1,\tfrac pq, 4,-\tfrac12), (-1,-3,\tfrac rs,-2), (-1,-\tfrac32,5,4),\\  
(-1, 3, 5,-\tfrac12), (-1,-3,4,-\tfrac 23), (-1,-\tfrac53, 4, -2), (-1,-4,-2,-3),\\ 
(-1,-\tfrac13,\tfrac52,\tfrac23), (-1,-3,-2,-3), (-1,-4,-2,-2), (-1,-3,-2,-4),\\ 
(-1,-3,-3,-3), (-1,-4,-3,-2), (-1,-5,-2,-2), (-1,-3,-2,-5),\\ 
(-1,-3,-4,-3), (-1,-4,-4,-2), (-1,-6,-2,-2), (-1,-3,-2,-6), \\ 
(-1,-3,-5,-3), (-1,-4,-5,-2), (-1,-7,-2,-2), (-1,-3,-3,-4), \\
(-1,-5,-3,-2)  
\end{gather*}
$\beta'=-\tfrac12$ is an exceptional slope on $M_5(\alpha')$ if and only if $\alpha'$ is one of 
\begin{gather*}
(-1,3, \tfrac rs, 4), (-1,\tfrac32, \tfrac43, \tfrac uv), (-1,\tfrac pq, \tfrac43, \tfrac52), (-1,\tfrac32, \tfrac54, -2), 
(-1,-3, \tfrac54, \tfrac52), (-1,3, \tfrac43, \tfrac83),\\ (-1,\tfrac53, \tfrac43, 4), (-1,\tfrac13, \tfrac53, \tfrac43), 
(-1,3, \tfrac23,5), (-1,4, \tfrac23, 4), (-1,3, \tfrac23,6), (-1,3, \tfrac34, 5),\\ 
(-1,4, \tfrac34,4), (-1,5, \tfrac23,4), (-1,3, \tfrac23,7), (-1,3, \tfrac45,5), 
(-1,6, \tfrac23,4), (-1,4, \tfrac45,4),\\ (-1,3, \tfrac23, 8), (-1,3, \tfrac56,5), 
(-1,4, \tfrac56, 4), (-1,7, \tfrac23, 4), (-1,3, \tfrac34, 6), (-1,5, \tfrac34, 4),\\ 
(-1,4, \tfrac23, 5), 
\end{gather*} 
$\beta'=\tfrac12$ is an exceptional slope on $M_5(\alpha')$ if and only if $\alpha'$ is one of 
\begin{gather*}
(-1,\tfrac13, \tfrac rs, \tfrac43), (-1,\tfrac23, \tfrac23, \tfrac uv), (-1,\tfrac pq, \tfrac23, \tfrac53), 
(-1,\tfrac23, \tfrac54, \tfrac23), 
(-1,-\tfrac13, \tfrac54, \tfrac53), (-1,\tfrac13, \tfrac23, \tfrac85),\\ (-1,\tfrac35, \tfrac23, \tfrac43), (-1,3, \tfrac13, 4), 
(-1,\tfrac13, \tfrac43, \tfrac54), (-1,\tfrac14, \tfrac43, \tfrac43), (-1,\tfrac13, \tfrac43, \tfrac65), (-1,\tfrac13,\tfrac54,\tfrac54),\\ 
(-1,\tfrac14,\tfrac54,\tfrac43), (-1,\tfrac15,\tfrac43,\tfrac43), (-1,\tfrac13,\tfrac54,\tfrac65), (-1,\tfrac15, \tfrac54, \tfrac43), 
(-1,\tfrac14, \tfrac43, \tfrac54), (-1,\tfrac13, \tfrac43, \tfrac76),\\ (-1,\tfrac13, \tfrac65, \tfrac54), (-1,\tfrac14, \tfrac65, \tfrac43), 
(-1,\tfrac16,\tfrac43, \tfrac43), (-1,\tfrac13, \tfrac43, \tfrac87), 
(-1,\tfrac13, \tfrac76, \tfrac54), (-1,\tfrac14, \tfrac76, \tfrac43),\\ (-1,\tfrac17, \tfrac43, \tfrac43), 
\end{gather*}
and $\beta'=2$ is an exceptional slope on $M_5(\alpha')$ if and only if $\alpha'$ is one of 
\begin{gather*}
(-1,-\tfrac13,\tfrac rs,\tfrac23), (-1,-\tfrac23,-2,\tfrac uv), (-1,\tfrac pq,-2,\tfrac13), (-1,-\tfrac23,-3,\tfrac43), 
(-1,\tfrac13,-3,\tfrac13),\\ (-1,-\tfrac13,-2,\tfrac25), (-1,-\tfrac35,-2,\tfrac23), (-1,-3,-\tfrac12,-2), 
(-1,-\tfrac13,4,\tfrac34), (-1,-\tfrac14,4,\tfrac23),\\ (-1,-\tfrac13,4,\tfrac45), (-1,-\tfrac13,5,\tfrac34), 
(-1,-\tfrac14,5,\tfrac23), (-1,-\tfrac15,4,\tfrac23), (-1,-\tfrac13,4,\tfrac56),\\ (-1,-\tfrac13,6,\tfrac34), 
(-1,-\tfrac14,6,\tfrac23), (-1,-\tfrac16,4,\tfrac23),  (-1,-\tfrac13,4,\tfrac67), (-1,-\tfrac13,7,\tfrac34),\\ 
(-1,-\tfrac14,7,\tfrac23), (-1,-\tfrac17,4,\tfrac23), (-1,-\tfrac13,5,\tfrac45), (-1,-\tfrac15,5,\tfrac23), 
(-1,-\tfrac14,4,\tfrac34).
\end{gather*}

As with the proof of Theorem \ref{5CL_thm}, Identities (\ref{first:eqn})--(\ref{figure8:eqn}) are used to 
identify equivalent filling instructions. 

The final step is to use (\ref{5CL_twist}) to obtain the filling 
instructions on $M_4$ and exceptional slopes shown in Table \ref{tb5}. Namely, the enumerated 
$M_5(\alpha')=M_5\big(-1, \frac pq, \frac rs, \frac uv, \varnothing\big)$ and $E(M_5(\alpha'))=\{\beta_i\}$ 
are identified with $M_4(\frac pq+1, \frac rs, \frac uv, \varnothing)$ and 
$E\big(M_4(\frac pq+1, \frac rs, \frac uv, \varnothing)\big)=\{\beta_i+1\}$ 
respectively. These $M_4(\alpha)$, $E(M_4(\alpha))$ are shown in Table \ref{tb6}. $\square$

%A=[[(1,-1),(3,-2),(5,1),(4,1),(0,0)], [(1,-1),(3,1),(5,1),(1,-2),(0,0)], [(1,-1),(3,-1),(4,1),(2,-3),(0,0)], [(1,-1),(5,-3),(4,1),(2,-1),(0,0)], [(1,-1),(4,-1),(2,-1),(3,-1),(0,0)], [(1,-1),(1,-3),(5,2),(2,3),(0,0)], [(1,-1),(3,-1),(2,-1),(3,-1),(0,0)], [(1,-1),(4,-1),(2,-1),(2,-1),(0,0)], [(1,-1),(3,-1),(2,-1),(4,-1),(0,0)], [(1,-1),(3,-1),(3,-1),(3,-1),(0,0)], [(1,-1),(4,-1),(3,-1),(2,-1),(0,0)], [(1,-1),(5,-1),(2,-1),(2,-1),(0,0)], [(1,-1),(3,-1),(2,-1),(5,-1),(0,0)], [(1,-1),(3,-1),(4,-1),(3,-1),(0,0)], [(1,-1),(4,-1),(4,-1),(2,-1),(0,0)], [(1,-1),(6,-1),(2,-1),(2,-1),(0,0)], [(1,-1),(3,-1),(2,-1),(6,-1),(0,0)], [(1,-1),(3,-1),(5,-1),(3,-1),(0,0)], [(1,-1),(4,-1),(5,-1),(2,-1),(0,0)], [(1,-1),(7,-1),(2,-1),(2,-1),(0,0)], [(1,-1),(3,-1),(3,-1),(4,-1),(0,0)], [(1,-1),(5,-1),(3,-1),(2,-1),(0,0)]]

\begin{remark}
The reduced list in Theorem \ref{4CL_thm} is surprisingly small (see Table \ref{tb6}). This occurs as 
the fillings listed above with $\pm2$, $\pm\frac12$ are equivalent. This can be seen by setting $i=2$ 
and $j=1$ below;   
\begin{gather*}
M_5(-1, \tfrac{v-u}{v}, \tfrac{r}{r-s}, \tfrac{q-p}{q}, \tfrac ji)\mathop{=}\limits_{(\ref{symeq10})}
M_5(-1, \tfrac pq, \tfrac rs, \tfrac uv, \tfrac ij) \mathop{=}\limits_{(\ref{blow:eqn})}
M_5(\tfrac{p+q}q, \tfrac{r-s}s, -1, \tfrac{u-v}v, \tfrac{i+j}j) \\ 
\mathop{=}\limits_{(\ref{symeq1})} M_5(-1, \tfrac{v}{u-v}, \tfrac{2s-r}s, \tfrac{p+q}p,-\tfrac{i}{j})
\mathop{=}\limits_{(\ref{symeq10})} M_5(-1,-\tfrac qp, \tfrac{2s-r}{s-r}, \tfrac{u-2v}{u-v},-\tfrac ji)
\end{gather*}
The only filling instructions that appear in more than one of the above lists are $(-1, -3,-\frac12,-2)$ and $(-1,-3,-2,-2)$.
\end{remark}

\begin{table}[h!]
\begin{center}
\begin{tabular}{|c|c|c|}
\hline & Additional & Exceptional\\ 
$f$ & exceptional & filling $M_5(f)(\beta_i)$\\ 
& slopes $\beta_i$ & \\ \hline \hline

\phantom{\Big|}$(-2, -\frac{1}{2}, 3, 3)$ & $\beta_1=-1$ & $\seiftre{S^2}2{-1}5241$\\
& $\beta_2=-\frac12$ &$\seifuno{A}2{-1}\bigb1211$  \\ \hline

\phantom{\Big|}$(-2, \frac 32, \frac 32, -2)$ & $\beta_1=-1$ & $\seifdue D2{-1}21\bigu0110 \seifdue D213{-1}$\\
& $\beta_2=-\frac12$ & $\seifuno{A}2{-1}\bigb1211$ \\ \hline

\phantom{\Big|}$(-2, -3, -\frac 12, -2)$& $\beta_1=-1$ & $\seiftre{S^2}213{-1}{11}{-2}$ \\ 
& $\beta_2=2$ & $\seifdue D2135\bigu0110 \seifdue D322{-1}$\\ \hline
\phantom{\Big|}$(-2, -\frac{1}{3}, 3, \frac{2}{3})$ & $\beta_1=-1$ & $\seiftre{S^2}2{-1}7{2}5{3}$\\ 
& $\beta_2=2$ & $\seifuno{A}2{3}\bigb0110$\\ \hline
\phantom{\Big|}$(-2, -\frac{1}{2}, 3, \frac{2}{3})$ & $\beta_1=-1$ & $\seiftre{S^2}2{-1}5{2}5{3}$\\ 
& $\beta_2=2$ & $\seiftre{S^2}2{1}3{-1}{11}{-2}$ \\ \hline

\phantom{\Big|}$(-2, -2, -2, -2)$ & $\beta_1=-1$ & $\seiftre{S^2}2{1}3{-1}7{-1}$\\ 
& $\beta_2=-2$ & $\seifdue D212{-1}\bigu{-1}51{-4} \seifdue D2{1}31$\\ \hline

\end{tabular}
\end{center}
\caption{All $M_5(f)$ with $f$ not factoring through $M_4$ and 
$e_{\tau}\big(M_5(f)\big) = 5$, $E_{\tau}\big(M_5(f)\big) = \{\beta_1, \beta_2, 0,1,\infty\}$. \label{tb1}}
\end{table}

\begin{table}[h!]
\begin{center}
\begin{tabular}{|c|c|c|}
\hline $f$ & Additional exceptional & Exceptional filling $M_5(f)$\\ 
& slopes $\beta$ & \\ \hline \hline

\phantom{\Big|}$(-2, \frac{p}{q}, 3, \frac{u}{v})$ & $-1$ & $F(-2,\tfrac pq,\tfrac{u+v}v,-2)$\\ \hline
$(-2, \frac{p}{q}, \frac{r}{s}, -2)$ & $-1$ & $F(\tfrac s{s-r},2,\tfrac q{2q-p},-3)$ \\ \hline

\phantom{\Big|}$(-2, \frac{3}{2}, \frac{3}{2}, \frac{u}{v})$ & $-1$ & $F(-3,\tfrac u{u+v}, 2,-2)$ \\ \hline

\phantom{\Big|}$(-2, \frac{p}{q}, \frac{5}{2}, -\frac{1}{2})$ & $-1$ & 
$F(-3,\tfrac{p-q}p,2,-2)$ \\ \hline 
\phantom{\Big|}$(-2, -2, \frac{r}{s}, -3)$ & $-1$ &$\seifuno{A}s{s-r}\bigb0110$ \\ \hline  
\phantom{\Big|}$(-2, -\frac{1}{2}, 4, \frac{u}{v})$ & $-1$ &$F(2,\tfrac32,\tfrac v{u+v},-2)$ \\ \hline 
\phantom{\Big|}$(-2, \frac{p}{q}, 4, -\frac{3}{2})$ & $-1$ &$F(2,\tfrac32,\tfrac qp,-2)$\\ \hline

\end{tabular}
\end{center}
\caption{All hyperbolic $M_5(f)$ with $f$ not factoring through $M_4$ and 
$e_{\tau}\big(M_5(f)\big) =4$, $E_{\tau}\big(M_5(f)\big) = \{\beta, 0,1,\infty\}$, part 1/5. \label{tb1param}}
\end{table}

\begin{table}[h!]
\begin{center}
\begin{tabular}{|c|c|c|}
\hline & Additional & Exceptional\\ 
$f$ & exceptional & filling $M_5(f)(\beta)$\\ 
& slopes $\beta$ & \\ \hline \hline

\phantom{\Big|}$(-2, 4, 5, -\frac{3}{2})$& $-1$ &$\seifuno{A}21\bigb0110$\\ 
$(-2, -\frac{1}{2}, 5, 3)$& & \\ \hline
$(-2, 3, 4, -\frac{4}{3})$&$-1$&$\seifuno{A}2{-1}\bigb0110$\\  \hline
$(-2, 3, \frac{3}{2}, -\frac{1}{2})$&$-1$ & $\seifuno{A}1{-3}\bigb0110$\\  \hline
$(-2, -2, 4, -\frac{5}{3})$&$-1$&$\seifdue D3{-1}21\bigu0110 \seifdue D232{-1}$\\
$(-2, -\frac{2}{3}, 4, -3)$& &\phantom{\Big|}\\ \hline

$(-2, \frac{2}{3}, \frac{5}{2}, -\frac{1}{3})$&$-1$&
$\seifdue D2123\bigu0110 \seifdue D322{-1}$ \phantom{\Big|}\\ \hline

$(-2, \frac{4}{3}, \frac{3}{2}, \frac{1}{3})$&$-1$&$\seifuno{A}2{-1}\bigb0110$ \phantom{\Big|}\\ \hline

\end{tabular}
\end{center}

\caption{All hyperbolic $M_5(f)$ with $f$ not factoring through $M_4$ 
and $e_{\tau}\big(M_5(f)\big) =4$, $E_{\tau}\big(M_5(f)\big) = \{\beta, 0,1,\infty\}$, part 2/5. 
\label{tb2}}
\end{table}

\begin{table}[h!] 
\begin{center}
\begin{tabular}{|c|c|c|}
\hline & Additional & Exceptional\\ 
$f$ & exceptional & filling \\ 
& slope $\beta$ & $M_5(f)(\beta)$ \\ \hline \hline

$(-2, -3, -2, -3)$& $-1$ & $\seiftre{S^2}2{1}3{-1}7{-1}$ \phantom{\Big|} \\
$(-2, -2, -2, -4)$& & \\ \hline

$(-2, -4, -2, -3)$ & & \\
$(-2, -2, -2, -5)$ & $-1$ & $\seiftre{S^2}2{1}4{-3}51$ \phantom{\Big|} \\
$(-2, -3, -3, -3)$ & & \\
$(-2, -2, -3, -4)$ & \phantom{\Big|} & \\ \hline

$(-2, -2, -4, -4)$ & & \\
$(-2, -5, -2, -3)$& $-1$ & $\seiftre{S^2}3{-2}3141$ \phantom{\Big|}  \\
$(-2, -3, -4, -3)$& & \\ 
$(-2, -2, -2, -6)$ & \phantom{\Big|} & \\ \hline

$(-2, -2, -2, -7)$ & & \\
$(-2, -6, -2, -3)$ &$-1$&$\seifdue D2{-1}2{-1}\bigu0110 \seifdue D2{-1}3{-1}$ \phantom{\Big|} \\
$(-2, -3, -5, -3)$ & & \\
$(-2, -2, -5, -4)$ & \phantom{\Big|} & \\ \hline

$(-2, -4, -3, -3)$ & $-1$ & $\seifdue D212{-1}\bigu{-1}21{-1} \seifdue D2131$ \phantom{\Big|}  \\
$(-2, -2, -3, -5)$& & \\ \hline 

$(-2, -3, -2, -4)$ & $-1$ & $\seifdue D212{-1}\bigu{-1}31{-2} \seifdue D2131$ \phantom{\Big|}\\ \hline
$(-2, \frac 32, \frac 52, -\frac 23)$ & $-1$ &  $\seifuno{A}2{-3}\bigb0110$ \phantom{\Big|}\\ \hline

%$(-2, -\tfrac12, 4, 3)$& $-1$ & $\seiftre{S^2}213{-1}9{-2}$ \phantom{\Big|}\\ \hline

%$(-2, 3, -2, -\frac 32) $ & $-1$ & $\seiftre{S^2}213{-1}52$ \phantom{\Big|} \\  \hline
\end{tabular}
\end{center}

\caption{All hyperbolic $M_5(f)$ with $f$ not factoring through $M_4$ and $e\big(M_5(f)\big) =4$, 
$E_{\tau}\big(M_5(f)\big) = \{\beta, 0,1,\infty\}$, part 3/5.
\label{tb3}}
\end{table}

\begin{table}[h!]
\begin{center}
\begin{tabular}{|c|c|c|}
\hline & Additional & Exceptional\\ 
$f$ & exceptional & filling $M_5(f)(\beta)$\\ 
& slopes $\beta$ & \\ \hline \hline

\phantom{\Big|} $(-2, -2, \tfrac14, 3)$ & $\tfrac12$ & $\seiftre{S^2}322{-1}9{-2}$\\ \hline
\phantom{\Big|} $(-2, \frac25, \frac34, \tfrac32)$ & $\tfrac12$ & $\seiftre{S^2}322{-1}3{2}$ \\ \hline
\phantom{\Big|} $(-2, \frac15, \frac43, \frac32)$ & $\tfrac12$ & $\seiftre{S^2}214{-3}51$ \\ \hline
\phantom{\Big|}$(-2, \frac14, \frac23, \frac53)$ & $\tfrac12$ & $\seifdue{D}2123\bigb0110\seifdue{D}322{-1}$ \\ \hline 
\phantom{\Big|}$(-2, \tfrac15, \frac32, \tfrac32)$ & $\tfrac12$ &$\seiftre{S^2}213{-1}7{-1}$ \\ \hline  
\phantom{\Big|}$(-2, 3, \tfrac13, 4)$ & $\tfrac12$ & $\seifdue{D}2134\bigb0110\seifdue{D}322{-1}$ \\ \hline 
\phantom{\Big|}$(-2, \frac13, \tfrac32, \frac43)$ & $\tfrac12$ & $\seifuno{A}13\bigb0110$\\ \hline
\phantom{\Big|}$(-2, -2, \tfrac15, 3)$& $\tfrac12$ & $\seifuno{A}21\bigb0110$\\ \hline
\phantom{\Big|}$(-2, \tfrac23, \tfrac35, \tfrac32)$& $\tfrac12$ & $\seifdue{D}213{-1}\bigb0110\seifdue{D}322{-1}$
\\ \hline
\phantom{\Big|} $(-2, \tfrac16, \tfrac32, \frac32)$ & $\tfrac12$ & $\seiftre{S^2}214{-3}51$\\  \hline
\phantom{\Big|} $(-2, \tfrac17, \frac{3}{2}, \frac{3}{2})$ & $\tfrac12$ & $\seiftre{S^2}3{-2}3141$\\  \hline
\phantom{\Big|} $(-2, \tfrac18, \frac{3}{2}, \frac{3}{2})$ & $\tfrac12$ & 
$\seifdue{D}2{-1}2{-1}\bigb0110\seifdue{D}2{-1}3{-1}$ \\  \hline
\phantom{\Big|} $(-2, \tfrac15, \frac65, \frac32)$ & $\tfrac12$ & 
$\seifdue{D}2{-1}2{-1}\bigb0110\seifdue{D}2{-1}3{-1}$ \\  \hline
\phantom{\Big|} $(-2, \tfrac38, \frac34, \frac{3}{2})$ & $\tfrac12$ & 
$\seifdue{D}3{-1}2{1}\bigb0110\seifdue{D}2{3}2{-1}$ \\  \hline
\phantom{\Big|} $(-2, \tfrac13, \tfrac23, \frac{5}{3})$& $\tfrac12$ & $\seiftre{S^2}3{-1}2152$\\ \hline

\phantom{\Big|}$(-2, \frac{2}{3}, \tfrac34, \tfrac23)$& $\tfrac12$ & $\seifuno{A}31\bigb0110$ \\ \hline

\phantom{\Big|} $(-2, \frac35, \frac23, \frac43)$& $\tfrac12$ &
$\seifdue D2132\bigu0110 \seifdue D322{-1}$ \\ \hline

$(-2, \frac 14, \frac 32, \frac 43)$ & $\frac12$ & $\seifdue D212{-1}\bigu{-1}41{-3} \seifdue D2131$ \phantom{\Big|}\\ \hline 

\end{tabular}
\end{center}

\caption{All hyperbolic $M_5(f)$ with $f$ not factoring through $M_4$ with $e\big(M_5(f)\big) =4$, 
$E_{\tau}\big(M_5(f)\big) = \{\beta, 0,1,\infty\}$, part 4/5. \label{tb4}}
\end{table}

\begin{table}[h!]
\begin{center}
\begin{tabular}{|c|c|c|}
\hline & Additional & Exceptional\\ 
$f$ & exceptional & filling $M_5(f)(\beta)$\\ 
& slopes $\beta$ & \\ \hline \hline

\phantom{\Big|} $(-2, -\tfrac13, 3, \tfrac23)$ & $2$ & $\seifuno{A}23\bigb0110$ \\ \hline
\phantom{\Big|} $(-2, -\frac23, -2, \tfrac23)$ & $2$ & $\seifdue{D}3{-1}2{1}\bigb0110\seifdue{D}2{-1}2{-1}$ \\ \hline
\phantom{\Big|} $(-2, -2, -\frac13, 3)$ & $2$ & $\seifdue{D}2125\bigb0110\seifdue{D}322{-1}$ \\ \hline
\phantom{\Big|}$(-2, -\frac13, -2, \frac25)$ & $2$ & $\seifdue{D}2131\bigb0110\seifdue{D}322{-1}$ \\ \hline
\phantom{\Big|}$(-2, -3, -\frac12, -2)$ & $2$ & $\seifdue{D}2135\bigb0110\seifdue{D}322{-1}$ \\ \hline 
\phantom{\Big|}$(-2, -\tfrac12, -\tfrac32, \tfrac13)$ & $2$ & $\seifdue{D}3{-1}21\bigb0110\seifdue{D}312{-1}$ 
\\ \hline 
\phantom{\Big|}$(-2, \frac13, -3, \frac13)$ & $2$ & $\seifuno{A}32\bigb0110$ \\ \hline

\phantom{\Big|}$(-2, -\tfrac12, -3, \tfrac13)$& $2$ & $\seifdue{D}2121\bigb0110\seifdue{D}322{-1}$ \\ \hline
\phantom{\Big|}$(-2, -\tfrac12, -3, \tfrac35)$& $2$ & $\seifdue{D}212{-1}\bigb0110\seifdue{D}322{-1}$
\\ \hline
\end{tabular}
\end{center}

\caption{All hyperbolic $M_5(f)$ with $f$ not factoring through $M_4$ with $e\big(M_5(f)\big) =4$, 
$E_{\tau}\big(M_5(f)\big) = \{\beta, 0,1,\infty\}$, part 5/5. \label{tb5}}
\end{table}

\begin{table}[h!] 
\begin{center}
\begin{tabular}{|c|c|c|}
\hline & Additional & Exceptional\\ 
$f$ & exceptional & filling \\ 
& slopes $\beta_i$ & $M_4(f)(\beta_i)$ \\ \hline \hline

$(-2, -2, -2)$ & $\beta_1=-2$ & $\seifdue D212{-1}\bigu{-1}41{-3} \seifdue D2{1}31$ \phantom{\Big|}\\ 
& $\beta_2=-1$ & $\seifuno{A}13\bigb0110$ \phantom{\Big|} \\ \hline

$(-2, -\frac{1}{2}, -2)$ & $\beta_1=-1$ & $\seifuno{A}23\bigb0110$ \phantom{\Big|} \\ 
& $\beta_2=3$ & $\seifdue D2123\bigu0110 \seifdue D322{-1}$ \\ \hline

$(-2, \frac{r}{s}, -2)$ & $-1$ & $\seifuno{A}s{s-r}\bigb0110$ \phantom{\Big|}\\ \hline 
$(\frac{p}{q}, 4, -\frac{1}{2})$ & $-1$ & $F(2,\tfrac32,\tfrac qp,-2)$ \phantom{\Big|}\\  \hline

$(4, 5, -\frac{1}{2})$ & & $\seifuno{A}21\bigb0110$ \phantom{\Big|}\\
$(-2, 4, -\frac{2}{3})$& & $\seifdue D3{-1}21\bigu0110 \seifdue D232{-1}$  \phantom{\Big|}\\
$(-2, -5, -3)$& $\beta=-1$ &$\seifdue D2{-1}2{-1}\bigu0110 \seifdue D2{-1}3{-1}$ \phantom{\Big|}\\
$(-2, -2, -6)$& &$\seifdue D2{-1}2{-1}\bigu0110 \seifdue D2{-1}3{-1}$ \phantom{\Big|}\\
$(\tfrac23,\tfrac52,\tfrac23)$ & \phantom{\Big|} & $\seifdue D2123\bigu0110 \seifdue D322{-1}$\\
$(-2, -2, -3)$ & \phantom{\Big|} & $\seiftre{S^2}2{1}3{-1}7{-1}$ \\ 
$(-2, -3, -3)$,& \phantom{\Big|} & $\seiftre{S^2}214{-3}51$\\
$(-2,-2,-4)$ & \phantom{\Big|} & $\seiftre{S^2}214{-3}51$ \\
$(-3, -4, -2)$&\phantom{\Big|} & $\seiftre{S^2}3{-2}3141$ \\ 
$(-2,-2,-5)$ &\phantom{\Big|} & $\seiftre{S^2}3{-2}3141$\\ 
$(-2,-3,-4)$ & \phantom{\Big|} & $\seifdue D212{-1}\bigu{-1}21{-1} \seifdue D2{1}31$\\
$(-3, -2, -3)$ &\phantom{\Big|} & $\seifdue D212{-1}\bigu{-1}31{-2} \seifdue D2{1}31$\\ \hline
\end{tabular}
\end{center}

\caption{All hyperbolic $M_4(f)$ with $f$ not factoring through $M_3$ and $e_{\tau}(M_4(f)\big)\geq 5$, 
together with $E_{\tau}\big(M_4(f)\big) = \{\beta_1,\beta_2, 0,1,2,\infty\}$ 
if $e_{\tau}(M_4(f)\big)=6$ and $E_{\tau}\big(M_4(f)\big)=\{\beta, 0,1,2,\infty\}$ 
if $e_{\tau}(M_4(f)\big) = 5$. 
\label{tb6}}
\end{table}

\clearpage

\section{Families of Cusped Manifolds and Proof of Theorem 1}\label{flash:proof}
 
We finish by showing that Theorem \ref{flash} will fall out as a consequence of 
Theorems \ref{5CL_thm} and \ref{4CL_thm}. We remarked in Section \ref{mt5cl_sec} that 
statements (\emph{i})--(\emph{iv}) from Theorem \ref{flash} are shown to hold for all hyperbolic 
$M_3(\alpha)$ in the Appendix of \cite{Magic}. Thus, we must show that Theorem \ref{flash} 
(\emph{i})--(\emph{iv}) hold for all hyperbolic $M_5(\alpha)$ when $\alpha$ does not factor 
through $M_3$. 

To prove Theorem \ref{flash} we will need to know the class of every $M_5(\alpha)(\beta)$ 
for $\beta\in E(M_5(\alpha))$ when $\alpha$ is found in Tables \ref{tb1}--\ref{tb5} and 
$M_4(\alpha)(\beta)$ for $\beta\in E(M_4(\alpha))$ when $\alpha$ is found in Tables \ref{tb6}. 
Theorems \ref{5CL_thm} and \ref{4CL_thm} together with Proposition \ref{fillingF} make this 
straightforward. To simplify matters we will say that a set of exceptional slopes 
$\{\alpha_1, \dots, \alpha_k\}$ is of type $\{\mathcal{C}_1,\dots,\mathcal{C}_k\}$ 
when $\alpha_i$ is of type $\mathcal{C}_i$ for each $i$.
The results are shown in Tables \ref{tc}-\ref{th}.

As highlighted in the introduction, Tables \ref{tc}-\ref{th} are of interest in their own right. Among other 
families, we highlight the family $M_5(-2,k,3,\frac{k+1}{3k+2})$ of hyperbolic knots in $S^3$ 
two integral toroidal surgeries and a type $Z$ surgery (see Table \ref{te}), the family $M_5(-2,\frac1k,3,\frac{k-1}k)$ with three 
type $Z$ fillings and a reducible surgery (see Table \ref{te}), the family $M_4(-2,\frac1k,-2)$ 
with four type $Z$ exceptional fillings and a toroidal filling (see Table \ref{th}), and the family 
$M_5(-2,\frac1k,3,\varnothing)$ of 2-cusped manifolds with four annular fillings on the $5^{th}$ cusp 
(see Table \ref{te}). These families are distinct from any obtained in \cite{Magic} because all hyperbolic 
fillings of $M_3$ have at least five exceptional slopes and a cyclic filling. The exceptional fillings of 
these families are written down using Proposition \ref{M_3M_5:prop} and shown in Table \ref{families_table}.

\begin{table}[htbp]
\begin{center}
\begin{tabular}{ |c|c| }\hline
\multicolumn{2}{ |c| }
{
$k\in\mathbb{Z}\backslash\{\pm1,0,2\}$, 
\qquad $E(M_5(-2, k, 3, \frac{k+1}{3k+2}))=\{-1, 0, 1, \infty\}$} \\ \hline
$\beta$ & $M_5(-2,k,3,\frac{k+1}{3k+2})(\beta)$\\ \hline \hline
$-1$ & $\seifdue D2{-1}k1\bigu0110 \seifdue D2{-1}{4k+3}{3k+2}$ \\ \hline 
$0$ & $\seiftre{S^2}{k-1}k21{8k+5}{-3k-2}$ \\ \hline 
$1$ & $\seifdue D313{-1}\bigu0110 \seifdue Dk1{2k+1}{-3k-2}$ \\ \hline
$\infty$ & $S^3$ \\ \hline\hline

\multicolumn{2}{ |c| }
{$k\in\mathbb{Z}\backslash\{-1,0,1,2\}$,
\qquad $E(M_5(-2, \frac1k, 3, \frac{k-1}{k}))=\{-1, 0, 1, \infty\}$} \\ \hline
$\beta$ & $M_5(-2,\frac1k, 3, \frac{k-1}{k})(\beta)$\\ \hline \hline
$-1$ & $\seiftre{S^2}2{-1}{2k-1}k{1-2k}2$ \\ \hline 
$0$ & $\seiftre{S^2}{1-k}1{1+2k}121$ \\ \hline 
$1$ & $\lens3{-1}\# \lens31$ \\ \hline
$\infty$ & $\seiftre{S^2}21k1{2k-3}{1-k}$ \\ \hline\hline

\multicolumn{2}{ |c| }{
$k\in\mathbb{Z}\backslash\{-2,-1,0,1,2\}$, 
\qquad $E(M_4(-2, \frac1k, -2))=\{-1, 0, 1, 2, \infty\}$}\\ \hline
$\beta$ & $M_4(-2,\tfrac1k,-2)(\beta)$\\ \hline \hline
$-1$ & $\seifuno{A}k{k-1}\bigb0110$ \\ \hline 
$0$ & $\seiftre{S^2}{2k-1}{1-k}2161$ \\ \hline 
$1$ & $\seiftre{S^2}4{-1}43{1-k}k$ \\ \hline
$2$ & $\seiftre{S^2}3{-1}3{-1}{1-2k}2$\\ \hline
$\infty$ & $\seiftre{S^2}k{1-k}2123$ \\ \hline\hline

\multicolumn{2}{ |c| }{
$k\in\mathbb{Z}\backslash\{-1,0,1,2\}$, 
\qquad $E(M_5(-2, \frac1k, 3,\varnothing))=\{-1, 0, 1, \infty\}$}\\ \hline

$\beta$ & $M_5(-2, \frac1k, 3,\varnothing)(\beta)$\\ \hline \hline
$-1$ & $\seifdue D2{-1}{1-2k}2$ \\ \hline 
$0$ & $\seifdue D21{1-k}1$ \\ \hline 
$1$ & $\seifdue D313{-1}$ \\ \hline
$\infty$ & $\seifdue D21k1$ \\ \hline\hline
\end{tabular}

\end{center}
\caption{The set of exceptional slopes and fillings for four families of cusped 
manifolds found in Tables \ref{tc}-\ref{th}.\label{families_table}}
\end{table}

\subsection{Proof of Theorem \ref{flash} (\emph{ii})--(\emph{iv})} The maximal distance between 
a lens space slope and a toroidal slope is known to be either three or four \cite{3or4}. 
Theorems \ref{5CL_thm} and \ref{4CL_thm} show that the only two $M_5(\alpha)$ with $\alpha$ not factoring 
through $M_3$ with two exceptional slopes $\beta, \gamma$ at distance greater than 3 are 
$M_4(-2,-2,-2)$ with $\beta=-2$ and $\gamma=2$ and $M_4(-2,-\frac12,-2)$ with $\beta=-1$ 
and $\gamma=3$. In all cases the fillings are toroidal (see Table \ref{tg}). So 
statement (\emph{ii}) holds. 

Theorems \ref{5CL_thm} and \ref{4CL_thm} tell us that if $e(M_5(\alpha))\geq6$ then $\alpha$ 
factors through $M_3$. So statement (\emph{iii}) holds. 

It is well known that the distance between two slopes $\frac pq$ and $\frac rs$ is $|ps-rq|$ 
(see \cite{sti}). From Theorems \ref{5CL_thm} and \ref{4CL_thm}, if the distance between two 
exceptional slopes on $M_5(\alpha)$ is greater than 4 then $\alpha$ factors through $M_3$. So 
statement (\emph{iv}) holds. 

\subsection{Proof of Theorem \ref{flash} (\emph{v})} We now consider the reducible fillings 
on $M_5$. We see directly from Tables \ref{tc}-\ref{th} that if $M_5(f)$ has two exceptional 
reducible slopes $\alpha$, $\beta$ then $\alpha,\beta\in\{0,1,\infty\}$ 
or $M_5(f)=M_4(g)$ for some filling instruction $g$ and $\alpha,\beta\in\{0,1,2,\infty\}$. 
In all cases, the reducible filling is described as a filling of $F$ by (\ref{5CLinf})-(\ref{4CL0}).
From Proposition \ref{fillingF}, $F(\tfrac ij,\tfrac kl,\tfrac nm,\tfrac tw)$ is reducible if and 
only if one of $\tfrac ij,\tfrac kl,\tfrac nm,\tfrac tw$ is zero or one of $\{\tfrac ij,\tfrac nm\}$, 
$\tfrac kl,\tfrac tw$ equals $\{\tfrac1{\eta},-\tfrac1{\eta}\}$ for some $\eta\in\mathbb{Z}$. 

If $F(\tfrac ij,\tfrac kl,\tfrac nm,\tfrac tw)=M_5(f)(\alpha)$ or 
$F(\tfrac ij,\tfrac kl,\tfrac nm,\tfrac tw)=M_4(g)(\beta)$ 
and one of $\tfrac ij,\tfrac kl,\tfrac nm,\tfrac tw$ is the zero slope then, by 
(\ref{5CLinf})-(\ref{4CL0}), one of the slopes in $f$ is in ${0,1,\infty}$ 
or one of the slopes in $g$ is in $\{0,1,2,\infty\}$ which makes $M_5(f)$ and 
$M_4(g)$ non-hyperbolic by Theorems \ref{5CL_thm} and \ref{4CL_thm}. 

By (\ref{4CLinf}) and Proposition \ref{fillingF}, if $\frac ab,\frac cd, \frac ef\neq0$ then 
$M_4(\frac ab,\frac cd, \frac ef)(\infty)$ is reducible if and only if 
$\frac{b-a}b=1\Rightarrow \frac ab=0$ making $M_4(\frac ab,\frac cd, \frac ef)$ 
non-hyperbolic by Theorem \ref{4CL_thm}. In the same way, if 
$\frac ab,\frac cd, \frac ef\neq0$ then $M_4(\frac ab,\frac cd, \frac ef)(2)$ 
is reducible if and only if $\frac ef=2$ which makes $M_4(\frac ab,\frac cd, \frac ef)$ 
non-hyperbolic by Theorem \ref{4CL_thm}. The remaining pair of possible exceptional slopes 
on $M_4(\frac ab,\frac cd, \frac ef)$ are 0 and 2. If both 0 and 2 are reducible slopes 
and $\frac ab,\frac cd, \frac ef\neq0$ then both $\frac{e-f}f=\frac1n$ and 
$\frac fe=\frac1m$ for some integers $n,m$ by (\ref{4CL2}) and (\ref{4CL0}) respectively. 
This implies that $\frac ef\in\{0,2\}$ which makes $M_4(\frac ab,\frac cd, \frac ef)$ non-hyperbolic 
by Theorem \ref{4CL_thm}. 

The final case to consider is when $M_5(\frac ab,\frac cd,\frac ef,\frac gh)$  
has a reducible pair of slopes in $\{0,1,\infty\}$ and the reducible fillings 
$F(\delta)$ and $F(\epsilon)$ have no zero slopes in $\delta$ or $\epsilon$. 
Namely, two of 
\begin{gather*} 
\text{(a)} \quad \big(\tfrac1n,\tfrac1n\big)\in
\left\{\big(-\tfrac ab,\tfrac dc\big),\big(\tfrac fe,-\tfrac gh\big)\right\} \quad 
\text{(b)} \quad \big(\tfrac1m,\tfrac1m\big)\in\left\{\big(\tfrac{a-b}b,\tfrac ef\big),
\big(\tfrac cd,\tfrac{g-h}h\big)\right\} \\  
\text{(c)} \quad \big(\tfrac1k,\tfrac1k\big)\in\left\{\big(\tfrac{b}{b-a},-\tfrac hg\big),
\big(\tfrac{c-d}c,\tfrac{e-f}f\big)\right\}
\end{gather*}
must hold. Each of (a), (b), (c) can hold in two ways. In all twelve ways that two of 
(a), (b), (c) hold, we find $\{0,1,\infty\}\cap\{\frac ab,\frac cd,\frac ef\}\neq\emptyset$ 
which makes $M_5(\frac ab,\frac cd,\frac ef,\frac gh)$ non-hyperbolic by Theorem \ref{main:teo} 
or that $\{-1,\frac12,2\}\cap\{\frac ab,\frac cd,\frac ef\}\neq\emptyset$ which makes 
$(\frac ab,\frac cd,\frac ef,\frac gh)$ factors through $M_4$ by (\ref{image-1}). 

Finally, $M_3$ does not have any exceptional reducible pairs (see Table 16 in the Appendix of \cite{Magic}). 
So statement (\emph{v}) holds. 

\subsection{Proof of Theorem \ref{flash} (\emph{i})} We first consider the $M_5(\alpha)$ and 
$M_4(\alpha)$ from Tables \ref{tc}--\ref{th} with a reducible $\beta$ not in $\{0,1,\infty\}$ 
or $\{0,1,2,\infty\}$ respectively. 
We will show that for such $M_5(\alpha)$, $\beta$ no $M_5(\alpha)(\gamma)=S^3$ for 
$\gamma\in\{0,1,\infty\}$, and for such $M_4(\alpha)$, $\beta$ no $M_4(\alpha)(\gamma)=S^3$ for 
$\gamma\in\{0,1,2,\infty\}$.

From Tables \ref{tc}--\ref{tg} we see that the only $M_5(\alpha)$ or $M_4(\alpha)$ with 
one boundary component with $\beta$ reducible not in $\{0,1,\infty\}$ 
or $\{0,1,2,\infty\}$ respectively are $M_5(-2,-n,n+3,-2)$ with $n\not\in\{-4,-3,-2,-1,0,1\}$ 
and $\beta=-1$ (found in Table \ref{tf}). It is easy to see from Table \ref{tf} that the remaining 
exceptional slopes are $0,1,\infty$ none of which give an $S^3$ filling. In particular, using Table \ref{tf}, 
$M_5(-2,n,n-3,-2, 0)$ is of type $Z$, $M_5(-2,n,n-3,-2, 1)$ is toroidal, and 
$M_5(-2,n,n-3,-2,\infty)$ is of type $T^H$. This shows that if $\beta$ is a 
reducible slope on $M_5(\alpha)$ or $M_4(\alpha)$ then $\beta$ is in $\{0,1,\infty\}$ or 
$\{0,1,2,\infty\}$ respectively.

Tables \ref{tc}--\ref{tg} show no $M_5(\alpha)$ and $M_4(\alpha)$ 
has an $S^H$ slope $\beta$ not in $\{0,1,\infty\}$ 
or $\{0,1,2,\infty\}$ respectively. So, to conclude statement (\emph{i}) for hyperbolic 
$M_5(\alpha)$ we need to show that no $M_5(\alpha)$ has reducible and $S^H$ 
slopes in $\{0,1,\infty\}$ and no $M_4(\alpha)$ has a reducible and $S^H$ 
slopes in $\{0,1,2,\infty\}$. 

We already know that an instruction $\alpha$ on $M_5$ factors through $M_3$ when 
$\alpha$ contains an instruction in $[\![((-1)_1, (-2)_2)]\!]$. If 
$\alpha=(\alpha_1,\alpha_2,\alpha_3,\alpha_4)$ is an instruction on $M_4$ 
that factors through $M_3$ then the filling instruction 
$\alpha'=(-1,\alpha_1-1,\alpha_2,\alpha_3,\alpha_4-1)$ 
(where $\alpha_i-1=\varnothing$ if $\alpha_i=\varnothing$ for $i=1,4$) 
on $M_5$ factors through $M_3$. Looking at (\ref{factor}) we see that if a 
filling instruction $\alpha$ on $M_4$ contains a slope in $\{-1, 3, \tfrac32,\tfrac12\}$ 
then $\alpha$ factors through $M_3$. 

We have 
\begin{gather*}
M_4(\tfrac ab,\tfrac cd,\tfrac ef, \tfrac{\beta_1}{\beta_2})\mathop{=}\limits_{(\ref{5CL_twist})} 
M_5(-1,\tfrac{a-b}b,\tfrac cd,\tfrac ef,\tfrac{\beta_1-\beta_2}{\beta_2})
\mathop{=}\limits_{(\ref{symeq7})}M_5\left(\tfrac f{f-e},-1,\tfrac d{d-c}, \tfrac{a-2b}{a-b},\tfrac{\beta_1-2\beta_2}{\beta_1-\beta_2}\right)
\\
\mathop{=}\limits_{(\ref{blow:eqn})} M_5(-1, \tfrac f{f-e}, \tfrac{2d-c}{d-c}, \tfrac{a-2b}{a-b}, \tfrac{\beta_2}{\beta_2-\beta_1})
\mathop{=}\limits_{(\ref{5CL_twist})} M_4\left(\tfrac{2f-e}{f-e},\tfrac{2d-c}{d-c},\tfrac{a-2b}{a-b},\tfrac{2\beta_2-\beta_1}{\beta_2-\beta_1}\right)
\end{gather*}
So, when $\beta=\tfrac{\beta_1}{\beta_2}\in \{0,1,2,\infty\}$ is a reducible slope on $M_4(\alpha)$ with 
$\alpha$ not factoring through $M_3$ we only need to consider $\beta\in\{0,1\}$.

From
\begin{equation}\label{1notS3}
M_4(n+2,\tfrac cd,-n)(1)\mathop{=}\limits_{(\ref{4CL1})}F(n,-1,\tfrac{c-d}d,-n-1)\mathop{=}\limits_{\text{Table } \ref{F4slopes}}
\seiftre{S^2}n1{n+2}{n+1}{c-d}d.
\end{equation}
we see that if $M_4(n+2,\tfrac cd,-n)(1)$ is reducible then $n=0,-2$ or $\frac cd=1$ which 
make $M_4(n+2,\tfrac cd,-n)$ non-hyperbolic. So, if $M_4(\alpha)$ is a hyperbolic knot in $S^3$ 
with a reducible filling in $\{0,1,2,\infty\}$ then we may assume that the reducible slope is 0.

In the case when 0 corresponds to a reducible filling on $M_4(\alpha)$ we have  
\[
M_4(\tfrac ab, \tfrac cd, \tfrac ef)(0)\mathop{=}\limits_{(\ref{4CL0})}F(\tfrac{2d-c}{c-d},\tfrac b{a-2b}, 2, \tfrac fe).
\]
From Proposition \ref{fillingF} we see that if $M_4(\tfrac ab, \tfrac cd, \tfrac ef)(0)$ is reducible then one of 
$2d-c,b,f=0$ or 
\begin{equation}\label{M_4red}
\tfrac b{a-2b}= \tfrac1n = -\tfrac fe
\end{equation} 
holds. If $2d-c=0$ then $\tfrac cd=2$ and $\alpha$ factors through $M_3$ (using (\ref{4CL_twist}) on the 
elements of (\ref{factor})). If $b$ or $f$ equal 0 then $\alpha$ is an exceptional filling instruction 
(see Theorem \ref{4CL_thm}). If 
(\ref{M_4red}) holds then $\tfrac ab=2+n$ and $\tfrac ef=-n$. For $M_4(\alpha)=M_4(2+n,\tfrac cd, -n)$ 
to be hyperbolic we require $n\not\in\{-2,-1,0\}$ and for $M_4(\alpha)$ to not factor through $M_3$ we require 
$n\not\in\{-3,-1,1\}$. Using Theorem \ref{4CL_thm} and the Proposition \ref{fillingF} we have 
\[
M_4(n+2,\tfrac cd,-n)(\infty)\mathop{=}\limits_{(\ref{4CLinf})}F(-1-n,\tfrac d{c-d},-1,n)\mathop{=}\limits_{\text{Table } \ref{F4slopes}}
\seiftre{S^2}{2+n}{1+n}n1d{c-d}.
\]
So $M_4(\alpha)=M_4(n+2,\tfrac cd,-n)(\infty)$ is of type $Z$ unless $n\in\{-3,-2,-1,0,1\}$ (which means $\alpha$ is exceptional 
or factors through $M_3$), or $d\in\{0,\pm1\}$. If $d=0$ then $\alpha$ is exceptional. If $d=\pm1$ then $\tfrac cd\in\mathbb{Z}$ 
so we may assume $d=1$ with out loss of generality. If $d=1$ then 
\[
M_4(\alpha)(\infty) \mathop{=}\limits_{(\ref{4CLinf})} F(-1-n,\tfrac 1{c-1},-1,n) \mathop{=}\limits_{\text{Table } \ref{F4slopes}}
\lens{cn^2+2cn+2}{\ast}.
\]
This means that if $M_4(\alpha)(\infty)=S^3$ then $n$ divides 1 or 3. The cases $n=-3,\pm1$ mean $\alpha$ is exceptional 
or factors through $M_3$. When $n=3$ we require $c\in\mathbb{Z}$ to satisfy $c(3)^2+2c(3)+2=\pm1$ which is 
impossible. 

From (\ref{1notS3}) we see that $M_4(n+2,\tfrac cd,-n)(1)$ is of type $Z$ unless $n\in\{-3,-1,0,1\}$ 
(which make $\alpha$ exceptional or factor through $M_3$) or $\frac cd=1+\frac1k$. 
When $\frac cd=1+\frac1k$ we find that  
\[
M_4(n+2,\tfrac cd,-n)(1)\mathop{=}\limits_{\text{Table } \ref{F4slopes}}
\lens{(k+1)n^2+2(k+1)n+2}{\ast}.
\]
So, if $M_4(n+2,\tfrac cd,-n)(1)=S^3$ then $n$ divides 1 or 3. The cases $\pm1,-3$ are excluded, and the case $n=3$ gives 
\[
M_4(n+2,\tfrac cd,-n)(1)=\lens{15k+17}{\ast}\neq S^3
\]
for any $k\in \mathbb{Z}$.

We have therefore shown that if $M_4(\alpha)$ is hyperbolic with $\alpha$ not factoring through $M_3$ 
with a reducible slope then $M_4(\alpha)$ does not have an $S^3$ filling.

The final case to consider is when $M_5(\alpha)(\beta)$ is reducible and $M_5(\alpha)(\gamma)=S^3$ 
with $\beta,\gamma \in\{0,1,\infty\}$ where $\alpha$ is a hyperbolic filling instruction that 
does not factor through $M_4$. There are six choices for $\beta$ and $\gamma$ but 
Lemma \ref{sym} allows us to assume that $\beta=1$ and $\gamma=0$. We have 
\[
M_5(\alpha)(1)=M_5(\tfrac ab, \tfrac cd,\tfrac ef,\tfrac gh)(1)\mathop{=}\limits_{(\ref{F_from_M_5})}
F(\tfrac{a-b}b,\tfrac cd, \tfrac ef,\tfrac{g-h}h)
\]
and 
\[
M_5(\alpha)(\infty)=M_5(\tfrac{n+1}{-n}, \tfrac n{-1},\tfrac dc,\tfrac 1{-k})(1)\mathop{=}\limits_{(\ref{5CLinf})}
F(-\tfrac{a}b,\tfrac fe, \tfrac dc,-\tfrac{g}h).
\]
From Table \ref{F4slopes} we see that if $M_5(\alpha)(1)$ is reducible then one of $a-b,c,e,g-h=0$ or 
$\tfrac{a-b}b=\tfrac1n=-\tfrac ef$ or $\tfrac{g-h}h=\tfrac1n=-\tfrac cd$. 
By a composition of (\ref{linksymeq1})--(\ref{linksymeq2}) we only need to consider the cases when one of 
$a-b,c=0$ or $\tfrac{a-b}b=\tfrac1n=-\tfrac ef$. If $c=0$ then $\alpha$ is exceptional. 
If $a-b=0$ then $\frac ab=1$ and $\alpha$ is exceptional. If $\tfrac{a-b}b=\tfrac1n=-\tfrac ef$ 
then $\tfrac ab=\tfrac{n+1}n$ and $\tfrac ef=-\tfrac1n$. Table \ref{F4slopes} tells us that in this case 
\[
M_5(\alpha)(1)=M_5(\tfrac{n+1}n, \tfrac cd,-\tfrac1n,\tfrac gh)(1)\mathop{=}\limits_{(\ref{F_from_M_5})}
F(\tfrac1n,\tfrac cd, -\tfrac1n,\tfrac{g-h}h)=\lens cd\# \lens{g-h}h.
\]
%So $M_5(\alpha)(1)$ is reducible when $c, g-h\neq\pm1$. 
%We have 

From Table \ref{F4slopes}, we see that if 
\[ M_5(\alpha)(\infty)=M_5(\tfrac{n+1}n,\tfrac cd,-\tfrac 1n,\tfrac gh)(\infty)
=F(-\tfrac{n+1}n,-n, \tfrac dc,-\tfrac{g}h)=S^3\] 
then one of $n+1=\pm1$, $n=\pm1$, $\frac cd=m\in \mathbb{Z}$ or $\frac gh=\frac 1k$ must hold. If $n+1=\pm1$ or $n=\pm1$ 
then $\alpha$ factors through $M_4$ or is exceptional. If $\frac cd=m\in \mathbb{Z}$ then, from Table \ref{F4slopes}, 
\[ M_5(\alpha)(\infty)=\seiftre{S^2}n{-1}g{-h}{(1-m)n-m}{n+1}=S^3\]
requires one of $n=\pm1$ (already excluded), $g=\pm1$, or $(1-m)n-m=\pm1$. If $(1-m)n-m=\pm1$ then $n=-1$ (which is excluded) 
or $m\in\{0,2\}$ which make $\alpha$ exceptional or factor through $M_4$. So, $\frac gh$ is necessarily of the form $\tfrac1k$. 
From Table \ref{F4slopes}, 
\[
M_5(\tfrac{n+1}{n}, \tfrac cd,-\tfrac1n,\tfrac 1{k})(\infty)\mathop{=}\limits_{(\ref{5CLinf})}
F(-\tfrac{n+1}n,-n, \tfrac dc,-\tfrac1k)=\seiftre{S^2}dc{n+1}{-n}{nk+1}n=S^3
\]
requires one of $\frac cd=m\in\mathbb{Z}$, $n+1=\pm1$ (which has been excluded), or $nk+1=\pm1$. If $nk+1=\pm1$ then 
$k=0$ or $n\in\{\pm1,\pm2\}$. The cases $k=0,n=\pm1,-2$ make $\alpha$ exceptional or factor through $M_4$, and if $n=2$ 
then $k=-1$ which makes $\alpha$ factor through $M_4$. So, if 1 is a reducible slope on a one cusped hyperbolic $M_5(\alpha)$ and 
$M_5(\alpha)(\infty)=S^3$ then we may assume that $\alpha=(\tfrac{n+1}n,m-\tfrac1n,\tfrac1k)$ where $k,m \neq\pm1,0,2$ and 
$n\neq\pm1,0,-2$. From Table \ref{F4slopes}, 
\[
M_5(\tfrac{n+1}{n}, m,-\tfrac1n,\tfrac 1{k})(\infty)\mathop{=}\limits_{(\ref{5CLinf})}
F(-\tfrac{n+1}n,-n, \tfrac1m,-\tfrac1k)=\lens{(nm+m-n)(-1-kn)-n(n+1}{\ast} =S^3
\]
if and only if 
\begin{equation}\label{last_eq}
m(1+kn+n+kn^2)=\pm1+kn^2-n^2 \rightarrow m=\frac{(k-1)n^2\pm1}{(n+1)(kn+1)}  
\end{equation}
because $n\neq1$. It is straightforward to verify that for $k>2$ and $n>1$, or $k>2$ and $n<-2$, $k<-1$ and $n>1$, 
or $k<-1$ and $n<-2$ that, in all eight cases, (\ref{last_eq}) leads to $0<m<3$. So, the only valid integer solutions 
to (\ref{last_eq}) make $\alpha$ exceptional or factor through $M_4$. 

This finishes the argument that no hyperbolic $M_5(\alpha)$ has both a reducible and $S^3$ filling.

\section{Tables}\label{tables_sec}

\begin{table}[htbp]
\begin{center}

\begin{tabular}{|c|}\hline 
$f=(-2, -\frac 12, 3, 3)$, $E(M_5(f))=\{-1, -\frac 12, 0, 1, \infty\}$, types 
$\{Z, T, S, Z, Z\}$ \phantom{\big|} \\ \hline
$f=(-2, \frac 32, \frac 32, -2)$, $E(M_5(f))=\{-1, -\frac 12, 0, 1, \infty\}$, 
types $\{T, T, T^H, T, T\}$ \phantom{\big|} \\ \hline

$f=(-2, -3, -\frac 12, -2)$, $E(M_5(f))=\{-1, 2, 0, 1, \infty\}$, 
types $\{Z, T, Z, Z, Z\}$ \phantom{\big|} \\ \hline
$f=(-2, -\frac 13, 3, \frac 23)$, $E(M_5(f))=\{-1, 2, 0, 1, \infty\}$, 
types $\{Z, T, Z, Z, Z\}$ \phantom{\big|} \\ \hline
$f=(-2, -\frac 12, 3, \frac 23)$, $E(M_5(f))=\{-1, 2, 0, 1, \infty\}$, 
types $\{Z, Z, Z, Z, Z\}$ \phantom{\big|} \\ \hline

$f=(-2, -2, -2, -2)$, $E(M_5(f))=\{-2, -1, 0, 1, \infty\}$, types 
$\{T, Z, Z, T, T^H\}$ \phantom{\big|} \\ \hline

\end{tabular}
\end{center}
\caption{Exceptional sets for $M_5(f)$ for $f$ in Table \ref{tb1}. \label{tc}}
\end{table}

\begin{table}[h!]
\begin{center}
\begin{tabular}{|c|}\hline 

$f=(-2, \frac 32, \frac 32, \frac uv)$, $\frac uv \in\mathbb{Q}\cup\{\varnothing\}\backslash\{-2, -1, 0, \frac 12, 1, 2\}$, 
$E(M_5(f))=\{-1, 0, 1, \infty\}$ \phantom{\big|}\\
\begin{minipage}{1.0\textwidth}
\begin{equation*}
\text{ types}  
\begin{cases} 
\{T,D,T,T\} \text{ if }\frac uv =\varnothing; \\
\{T, S^H, T, T\} \text{ if } |5u-14v|=1; \\ 
\{T, T^H, Z, T\} \text{ if } |u-v| = 1; \\
\{Z, T^H, T, Z \} \text{ if } |u|=1; \\
\{T, T^H, T, T\} \text{ otherwise.}
\end{cases}
\end{equation*}
\end{minipage} 
\\ \hline

$f=(-2, \frac pq, \frac 52, -\frac 12)$, $\frac pq \in\mathbb{Q}\cup\{\varnothing\}\backslash\{-1, 0, \frac 12, 1, 2\}$, 
$E(M_5(f))=\{-1, 0, 1, \infty\}$ \phantom{\big|}\\
\begin{minipage}{1.0\textwidth}
\begin{equation*}
\text{ types}  
\begin{cases}
\{T \& A, A, T \& A, A\} \text{ if } \frac pq=\varnothing ; \\ 
\{Z, T^H, T, Z\} \text{ if } \frac pq=1+ \frac 1n ; \\
\{T, Z, Z, Z\} \text{ if } |p|=1; \\ 
\{T, Z, T, T^H\} \text{ if } |q| = 1; \\
\{T, Z, T, Z\} \text{ otherwise.}
\end{cases}
\end{equation*}
\end{minipage} 
\\ \hline

$f=(-2, -2, \frac rs, -3)$, $\frac rs \in\mathbb{Q}\cup\{\varnothing\}\backslash\{-1, 0, \frac 12, 1, 2\}$, 
$E(M_5(f))=\{-1, 0, 1, \infty\}$ \phantom{\big|} \\
\begin{minipage}{1.0\textwidth}
\begin{equation*}
\text{ types} 
\begin{cases} 
\{T, A, T \& A, A\} \text{ if } \tfrac rs=\varnothing; \\
\{T, T^H, T, Z\} \text{ if } |r-s|=1; \\
\{T, Z, Z, Z\} \text{ if } |r|=1; \\
\{T, Z, T, T^H\} \text{ if } |s|=1; \\
\{T, Z, T, Z\}\text{ otherwise.}
\end{cases}
\end{equation*} 
\end{minipage} \\ \hline
\end{tabular}
\end{center}
\caption{Exceptional sets for $M_5(f)$ for $f$ in Table \ref{tb1param}, part 1/4. \label{tc1}}
\end{table}

\clearpage

\begin{table}[htbp]
\begin{center}
\begin{tabular}{|c|}\hline 

%$f=(-2, -2, \frac rs, -3)$, $\frac rs \in\mathbb{Q}\cup\{\varnothing\}\backslash\{-1, 0, \frac 12, 1, 2\}$, 
%$E(M_5(f))=\{-1, 0, 1, \infty\}$ \phantom{\big|} \\
%\begin{minipage}{1.0\textwidth}
%\begin{equation*}
%\text{ types} 
%\begin{cases} 
%\{T, A, T, A\} \text{ if } \tfrac rs=\varnothing; \\
%\{T, T^H, T, Z\} \text{ if } |r-s|=1; \\
%\{T, Z, Z, Z\} \text{ if } |r|=1; \\
%\{T, Z, T, T^H\} \text{ if } |s|=1; \\
%\{T, Z, T, Z\}\text{ otherwise.}
%\end{cases}
%\end{equation*} 
%\end{minipage} \\ \hline

\hline $f=(-2, -\frac 12, 4, \frac uv)$, $\frac uv \in\mathbb{Q}\cup\{\varnothing\}\backslash\{-1, 0, \frac 12, 1, 2\}$, 
$E(M_5(f))=\{-1, 0, 1, \infty\}$ \phantom{\big|}\\
\begin{minipage}{1.0\textwidth}
\begin{equation*}
\text{ types}  
\begin{cases}
\{T \& A, A, A, A\} \text{ if } \tfrac uv =\varnothing; \\
\{T^H,Z, Z, Z\} \text{ if } \frac uv = -2; \\ 
\{Z, S, Z, Z\} \text{ if } \frac uv = 3; \\
\{Z, T^H, Z, S\} \text{ if } \frac uv = 4; \\ 
\{T, Z, S, Z\} \text{ if } \frac uv = \frac 32; \\
\{Z, Z, Z, Z\} \text{ if } \frac uv \in \mathbb{Z}\backslash\{-2, -1, 0, 1, 2, 3, 4\}; \\
\{T, T^H, Z, Z\} \text{ if } |u-3v|=1, \,\frac uv\neq 4; \\
\{T, Z, T^H, Z\} \text{ if } |2u-3v|=1; \\
\{T, Z, Z, T^H\} \text{ if } |v-4u|=1, \frac uv\neq -\frac 72,-5\; \\
\{T, Z, Z, Z\} \text{ otherwise.}
\end{cases}
\end{equation*}
\end{minipage} 
\\ \hline

$f=(-2, \frac pq, 4, -\frac 32)$, $\frac pq \in\mathbb{Q}\cup\{\varnothing\}\backslash\{-1, 0, \frac 12, 1, 2\}$, 
$E(M_5(f))=\{-1, 0, 1, \infty\}$ \phantom{\big|}\\
\begin{minipage}{1.0\textwidth}
\begin{equation*}
\text{ types}  
\begin{cases}
\{T \& A, A, T \& A, A\} \text{ if } \tfrac pq =\varnothing; \\ 
\{Z, Z, T, T^H\} \text{ if } |q| = 1; \\
\{T, T^H, T, Z\} \text{ if } |p-q|=1; \\
\{T, Z, Z, Z\} \text{ if } |p|=1; \\ 
\{T, Z, T, Z\} \text{ otherwise.}
\end{cases}
\end{equation*}
\end{minipage} 
\\ \hline

\end{tabular}
\end{center}
\caption{ Exceptional sets for $M_5(f)$ for $f$ in Table \ref{tb1param}, part 2/4. \label{td}}
\end{table}

\begin{table}[htbp]
\begin{center}
\begin{tabular}{|c|}\hline
$f=(-2, \frac pq, 3, \frac uv)$, $E(M_5(f))=\{-1, 0, 1, \infty\}$,
$\frac pq, \frac uv\in\mathbb{Q}\cup\{\varnothing\}\backslash\{0, 1, -1, \frac 12, 2\}$ \phantom{\big|}\\ \hline \hline
\begin{minipage}{1.0\textwidth}
\begin{equation*}
-1 \text{ is of type}  
\begin{cases} 
A \text{ if } \frac pq = \varnothing \text{ and } |u+v|=1, \text{ or }\frac uv=\varnothing \text{ and }|p|=1;\\ 
T \& A \text{ if } \tfrac pq=\tfrac uv=\varnothing, \text{ or } \frac pq = \varnothing \text{ and } |u+v|\neq 1,\\ 
\phantom{T \& A } \text{ or }\frac uv=\varnothing \text{ and }|p|\neq 1;\\
T^H \text{ if } |p| = |u+v|=1; \\
Z \text{ if } |p|=1 \text{ and }|u+v|\neq 1, \text{ or }|u+v|=1 \text{ and }|p|\neq 1; \\
T \text{ otherwise} 
%\text{ if } |p|>1 \text{ and } |u+v|>1, \text{ or } \frac pq=\frac uv=\varnothing, 
%\text{ or } \frac pq=\varnothing \text{ and } |u+v|>1, \text{ or } \frac uv=\varnothing \text{ and } |p|>1.
\end{cases}
\end{equation*}
\end{minipage} 
\\ \hline

\begin{minipage}{1.0\textwidth}
\begin{equation*}
0 \text{ is of type}  
\begin{cases} 
A \text{ if } \frac pq=\varnothing \text{ and } \tfrac uv\neq3,3+\tfrac1n, 
\text{ or } \frac uv=\varnothing \text{ and }\tfrac pq\neq 1+\tfrac1n;\\
S\,\, \&\,\, D \text{ if } \frac uv=3\text{ and }\frac pq=\varnothing; \\
D \text{ if } \frac pq=1+\frac1n \text{ and } \frac uv=\varnothing, \text{ or } 
\frac uv=3+\tfrac 1n\text{ and }\frac pq=\varnothing; \\

S \text{ if } \frac uv=3\text{ and }|p-q|\neq 1; \\
S^H \text{ if } \frac pq=+\frac1n \text{ and } |(3+2n)u-(7+6n)v|=1;\\
\phantom{S^H} \text{ or } \frac uv=3+\frac1k \text{ and } |(3+2k)p-(1+2k)q|=1,\\

T^H \text{ if } \frac uv=3 \text{ and }|p-q| = 1,\\ 
\phantom{T^H} \text{ or }\frac pq=+\frac1n \text{ and } |(3+2n)u-(7+6n)v|\neq1,\\
\phantom{S^H} \text{ or } \frac uv=3+\frac1k \text{ and } |(3+2k)p-(1+2k)q|\neq1;\\

Z \text{ otherwise}.
\end{cases}
\end{equation*}
\end{minipage} 
\\ \hline

\begin{minipage}{1.0\textwidth}
\begin{equation*}
1 \text{ is of type}  
\begin{cases} 
T\,\,\&\,\, A\text{ if }\frac pq=\frac uv=\varnothing, 
\text{ or }\frac pq =\varnothing \text{ and } |u-v|\neq 1,\\
\phantom{T\,\,\&\,\, A}\text{ or } \frac uv=\varnothing\text{ and }|p|\neq 1;\\
A \text{ if }\frac pq =\varnothing \text{ and } |u-v|=1,\text{ or } \frac uv=\varnothing\text{ and }|p|=1;\\ 
S \text{ if } \frac pq=\frac1k=1-\frac uv; \\
T^H \text{ if } \frac pq=\frac1k\text{ and }|(1-k)v+ku|=1, \text{ or } \frac uv=1+\frac1k \text{ and } |kp+q|=1; \\
Z \text{ if } \frac pq=\frac1k\text{ and }|(1-k)v+ku|\neq1, \text{ or } \frac uv=1+\frac1k \text{ and } |kp+q|\neq1; \\
T \text{ otherwise.}
\end{cases}
\end{equation*}
\end{minipage} 
\\ \hline

\begin{minipage}{1.0\textwidth}
\begin{equation*}
\infty \text{ is of type}  
\begin{cases} 
A\text{ if }\frac pq=\frac uv=\varnothing, \text{ or } \frac pq=\varnothing \text{ and } |v-3u|\neq1, 
\text{ or }\frac uv=\varnothing \text{ and } |q|\neq1;\\
S\,\,\&\,\, D \text{ if } \frac uv=\frac13\text{ and } \frac pq=\varnothing;\\
S\text{ if }\frac uv=\frac13\text{ and }|q|\neq1;\\

S^H \text{ if } \frac pq=k\text{ and } |(1+2k)v-(1+6k)u|=1, \\
\phantom{S^H}\text{ or } v-3u=\epsilon =\pm1 \text{ and }|(1+2\epsilon u)q+2p|=1; \\

T^H \text{ if } \frac uv=\frac13\text{ and }|q|\neq1, \text{ or } 
\frac pq=k\text{ and } |(1+2k)v-(1+6k)u|=1, \\
\phantom{T^H }\text{ or } v-3u=\epsilon =\pm1 \text{ and }|(1+2\epsilon u)q+2p|\neq1; \\

Z \text{ otherwise}.
\end{cases}
\end{equation*}
\end{minipage} 
\\ \hline

\end{tabular}
\end{center}
\caption{Exceptional sets for $M_5(f)$ for $f$ in Table \ref{tb1param}, part 3/4. \label{te}}
\end{table}

\clearpage

\begin{table}[htbp]
\begin{center}
\begin{tabular}{|c|}\hline
$f=(-2, \frac pq, \frac rs, -2)$, 
$E(M_5(f))=\{-1, 0, 1, \infty\}$, \phantom{\big|}
$\frac pq, \frac rs\in\mathbb{Q}\cup\{\varnothing\}\backslash\{-1, 0, \frac 12, 1, 2\}$ \phantom{\big|}\\ \hline \hline

\begin{minipage}{1.0\textwidth}
\begin{equation*}
-1 \text{ is of type}  
\begin{cases} 
T \& A \text{ if } \frac pq=\tfrac rs =\varnothing;\\
A \text{ if } \frac pq=\varnothing \text{ and } |s|=1\text{, or } \frac rs=\varnothing \text{ and } |q|=1; \\
S \text{ if } \frac pq=4-\frac {r}s=n \text{ for some }n\in\mathbb{Z}\backslash\{-1,0,1,2\}; \\
T^H \text{ if } |q|=1 \text{ or } |s|=1\text{, and }4-\frac pq +\frac 1n=\frac rs \text{ for some }n\in\mathbb{Z}; \\
Z\text{ if } |q|=1 \text{ or } |s|=1\text{, and }4-\frac pq +\frac 1n\neq\frac rs \text{ for any }n\in\mathbb{Z}; \\
T \text{ otherwise}.
\end{cases}
\end{equation*}
\end{minipage} 
\\ \hline

\begin{minipage}{1.0\textwidth}
\begin{equation*}
0 \text{ is of type}  
\begin{cases}
A \text{ if } \frac pq=\tfrac rs =\varnothing,\text{ or } \frac pq=\varnothing \,\, \& \,\,  |r-s|\neq1, 
\text{ or } \frac rs=\varnothing \,\,\&\,\, |p-q|\neq1;\\
D \text{ if } \frac pq=\varnothing\text{ and } |r-s|=1, \text{ or } \frac rs=\varnothing \text{ and }|p-q|=1;\\ 
S^H \text{ if } \frac pq=1+\frac1n\text{ and }|(6-5n)s+(5n-1)r|=1, \\
\phantom{S^H} \text{ or } \frac rs=1+\frac1n\text{ and }|(5n-1)p+(6-5n)q|=1; \\
T^H \text{ if } \frac pq=1+\frac1n\text{ and }|(6-5n)s+(5n-1)r|\neq1, \\
\phantom{S^H} \text{ or } \frac rs=1+\frac1n\text{ and }|(5n-1)p+(6-5n)q|\neq1; \\
Z \text{ otherwise}.
\end{cases}
\end{equation*}
\end{minipage} 
\\ \hline

\begin{minipage}{1.0\textwidth}
\begin{equation*}
1 \text{ is of type}  
\begin{cases}
T\& A \text{ if } \frac pq=\tfrac rs =\varnothing;\\
A \text{ if } \frac pq=\varnothing\text{ and } |r|=1, \text{ or } \frac rs=\varnothing \text{ and }|p|=1;\\ 
T^H \text{ if } |p|=1\text{ and }|r|=1; \\
Z \text{ if } |p|=1\text{ and }|r|\neq1, \text{ or } |r|=1 \text{ and } |p|\neq1; \\
T \text{ otherwise}.
\end{cases}
\end{equation*}
\end{minipage} 
\\ \hline

\begin{minipage}{1.0\textwidth}
\begin{equation*}
\infty \text{ is of type}  
\begin{cases}
T\& A \text{ if } \frac pq=\tfrac rs =\varnothing;\\
A \text{ if } \frac pq=\varnothing \text{ and } |s|=1, \text{ or }\frac rs=\varnothing \text{ and } |q|=1;\\ 
T^H \text{ if } |s|=1\text{ and }|q|=1; \\
Z \text{ if } |s|=1\text{ and }|q|\neq1, \text{ or } |q|=1\text{ and }|s|\neq1; \\
T \text{ otherwise}.
\end{cases}
\end{equation*}
\end{minipage} 
\\ \hline

\end{tabular}
\end{center}
\caption{Exceptional sets for $M_5(f)$ for $f$ in Table \ref{tb1param}, part 4/4. 
\label{tf}}
\end{table}

\begin{table}[htbp]
\begin{center}
\begin{tabular}{|c||c|}
\hline $f$ & $E(M_5(f))$ and types \\ \hline \hline

$(-2, 4, 5, -\frac 32),$ $(-2, 3, 4, -\frac 43),$ & \\
$(-2, -2, 4, -\frac 53),$ $(-2, -2, -3, -5),$ & \\
$(-2, -4, -3, -3),$ $(-2, -2, -2, -7),$ & $E(M_5(\alpha))=\{-1, 0, 1, \infty\}$,\\
$(-2, -6, -2, -3),$ $(-2, -2, -5, -4),$ & \\
$(-2, -3, -5, -3),$  $(-2, -3, -2, -4),$ & types $\{T, Z, T, T^H\}$\phantom{\big|}  \\ \hline

$(-2, -\frac 12, 5, 3)$ & $E(M_5(\alpha))=\{-1, 0, 1, \infty\}$, $\{T, S, Z, Z\}$\phantom{\big|} \\ \hline
$(-2, 3, \frac 32, -\frac 12)$ & $E(M_5(\alpha))=\{-1, 0, 1, \infty\}$, $\{T, T^H, T, T^H\}$\phantom{\big|} \\ \hline
$(-2, -\frac 23, 4, -3)$ & $E(M_5(\alpha))=\{-1, 0, 1, \infty\}$, $\{T, Z, T, Z\}$\phantom{\big|} \\ \hline
%$(-2, 3, \frac 13, -3)$ & $E(M_5(\alpha))=\{-1, 0, 1, \infty\}$, $\{T, Z, Z, Z\}$\phantom{\big|} \\ \hline
%$(-2, \frac 53, \frac 32, -\frac 12),$ $(-2, \frac 32, \frac 53, -\frac 13),$ & 
$(-2, \frac 23, \frac 52, -\frac 13)$, $(-2, \frac 43, \frac 32, \frac 13)$ & $E(M_5(\alpha))=\{-1, 0, 1, \infty\}$, 
\phantom{\big|} $\{T, T^H, T, Z\}$\phantom{\big|} \\ \hline

$(-2, -2, -2, -4),$ $(-2, -3, -2, -3),$ & \\
$(-2, -2, -2, -5),$ $(-2, -2, -2, -6),$ & $E(M_5(\alpha))=\{-1, 0, 1, \infty\}$,\\
$(-2, -5, -2, -3),$ $(-2, -2, -4, -4),$ & \\
$(-2, -3, -4, -3),$ $(-2, -3, -3, -3),$ & types $\{Z, Z, T, T^H\}$\phantom{\big|} \\
$(-2, -2, -3, -4)$ & \\ \hline

$(-2, -4, -2, -3)$ & $E(M_5(\alpha))=\{-1, 0, 1, \infty\}$, $\{Z, Z, T, T^H\}$\phantom{\big|} \\ \hline
$(-2, \frac 32, \frac 52, -\frac 23)$ & $E(M_5(\alpha))=\{-1, 0, 1, \infty\}$, $\{T, T^H, T, T\}$ 
\phantom{\big|} \\ \hline
\end{tabular}
\end{center}
\caption{\,\,\,\, Exceptional sets for $M_5(f)$ for $f$ in Tables \ref{tb2} and \ref{tb3}. \label{tb}}
\end{table}

\begin{table}[htbp]
\begin{center}

\begin{tabular}{|c||c|}
\hline $f$ & $E(M_5(f))$ and types \\ \hline \hline
$(-2, -\frac 13, 3, \frac 23)$, $(-2, -\frac 23, -2, \frac 23)$,& \phantom{\big|} \\ 
$(-2, \frac 13, -3, \frac 13)$, $(-2, -\frac 13, -2, \frac 25)$, & $E(M_5(\alpha))=
\{2, 0, 1, \infty\}$ \phantom{\big|} \\ 
$(-2, -3, -\frac 12, -2)$, $(-2, -\frac 12, -\frac 32, \frac 13)$, & types $\{T, Z, Z, Z\}$ \phantom{\big|} \\ 
$(-2, -\frac 12, -3, \frac 13)$, $(-2, -\frac 12, -3, \frac 35)$ &  \phantom{\big|} \\ \hline

$(-2, -2, -\frac 13, 3)$ & $E(M_5(f))=\{2, 0, 1, \infty\}$, $\{T, S, Z, Z\}$ \phantom{\big|} \\ \hline

$(-2, -2, \frac 14, 3)$ & $E(M_5(f))=\{\frac 12, 0, 1, \infty\}$, $\{Z, S, Z, Z\}$\phantom{\big|} \\ \hline

$(-2, \frac 25, \frac 34, \frac 32)$ & $E(M_5(f))=\{\frac 12, 0, 1, \infty\}$, $\{T^H, T^H, Z, T\}$\phantom{\big|} \\ \hline

$(-2, \frac 15, \frac 43, \frac 32)$, $(-2, \frac 15, \frac 32, \frac 32)$, & $E(M_5(\alpha))
=\{\frac 12, 0, 1, \infty\}$, \phantom{\big|} \\
$(-2, \frac 16, \frac 32, \frac 32)$, $(-2, \frac 17, \frac 32, \frac 32)$, & types $\{Z, T^H, Z, T\}$\phantom{\big|} \\ 
$(-2, \frac 13, \frac 23, \frac 53)$ & \\ \hline

$(-2, \frac 14, \frac 23, \frac 53)$, $(-2, \frac 14, \frac 32, \frac 43)$, & \\
$(-2, \frac 13, \frac 32, \frac 43)$, $(-2, \frac 18, \frac 32, \frac 32)$, & 
$E(M_5(\alpha))=\{\frac 12, 0, 1, \infty\}$, \phantom{\big|} \\ 
$(-2, \frac 15, \frac 65, \frac 32)$, $(-2, \frac 38, \frac 34, \frac 32)$, & types $\{T, T^H, Z, T\}$\phantom{\big|}  \\ 
$(-2, \frac 23, \frac 34, \frac 23)$, $(-2, \frac 23, \frac 35, \frac 32)$ & \\ \hline

$(-2, 3, \frac 13, 4)$ & $E(M_5(\alpha))=\{\frac 12, 0, 1, \infty\}$, $\{T, T^H, Z, Z\}$\phantom{\big|} \\ \hline

$(-2, -2, \frac 15, 3)$ & $E(M_5(\alpha))=\{\frac 12, 0, 1, \infty\}$, $\{T, S, Z, Z\}$\phantom{\big|} \\ \hline

$(-2, \frac 35, \frac 23, \frac 43)$ & $E(M_5(\alpha))=\{\frac 12, 0, 1, \infty\}$, 
$\{T, T^H, Z, T\}$\phantom{\big|} \\ \hline

\end{tabular}
\end{center}
\caption{\,\,\,\, Exceptional sets for $M_5(f)$ for $f$ in Tables \ref{tb4} and \ref{tb5}. \label{ta}}
\end{table}

\begin{table}[htbp]
\begin{center}

\begin{tabular}{|c||c|}
\hline $f$ & $E(M_4(f))$ and types \\ \hline\hline

$(-2, -\frac 12, -2)$ & $\{-1, 0, 1, 2, 3, \infty\}$, $\{T, Z, Z, Z, T, Z\}$ \phantom{\big|} \\ \hline

$(-2, -2, -2)$ & $\{-2, -1, 0, 1, 2, \infty\}$, $\{T, T, Z, Z, T, T^H\}$ \phantom{\big|} \\ \hline

$(-3, -2, -3)$, $(-2, 4, -\frac 23)$,& \phantom{\big|} \\ 
$(-2, -3, -4)$,  $(4, 5, -\frac 12)$ & $E(M_5(f))=\{-1, 0, 1, 2, \infty\}$ \phantom{\big|} \\ 
$(-2, -5, -3)$, $(-2, -2, -6)$, & \phantom{\big|} types $\{T, Z, Z, T, T^H\}$ \phantom{\big|}
\\ \hline

$(-3,-4,-2)$, $(-2, -3, -3)$, & $E(M_5(f))=\{-1, 0, 1, 2, \infty\}$ \phantom{\big|} \\ 
$(-2, -2, -4)$, $(-2, -2, -3)$, & types $\{Z, Z, Z, T, T^H\}$ \phantom{\big|} \\
$(-2, -2, -5)$ & \\ \hline 

$(\frac23, \frac52, \frac23)$ & $\{-1, 0, 1, 2, \infty\}$, $\{T, Z, Z, Z, Z\}$  \phantom{\big|} \\ \hline

\end{tabular}
\end{center}
\caption{\,\,\,\, Exceptional sets for $M_5(f)$ for $f$ in Table 
\ref{tb6}, part 1/2. \label{tg}}
\end{table}

\begin{table}[htbp]
\begin{center}

\begin{tabular}{|c|}
\hline $f=(-2, \frac rs, -2)$, $\frac rs \in\mathbb{Q}\cup\{\varnothing\}\backslash\{-2,-1, 0, \frac 12, 1, \frac 32, 2, 3\}$, 
$E(M_4(f))=\{-1, 0, 1, 2, \infty\}$ \phantom{\big|}\\
\begin{minipage}{1.0\textwidth}
\begin{equation*}
\text{ types}  
\begin{cases} 
\{T, A, A, T \& A, A\} \text{ if }\frac rs=\varnothing ; \\
\{T,T^H, Z, T, Z \} \text{ if }\frac rs=2+\frac 1k ; \\
\{T,Z, T^H, T, Z \} \text{ if }\frac rs=1+\frac 1k ; \\
\{T, Z, Z, Z, Z\} \text{ if } \frac rs=\frac 1k; \\
\{T, Z, Z, T, T^H\} \text{ if } \frac rs\in\mathbb{Z}; \\ 
\{T, Z, Z, T, Z\} \text{ otherwise.}
\end{cases}
\end{equation*}
\end{minipage} 
\\ 
\hline \phantom{\big|} $f=(\frac pq, 4, -\frac 12)$, $\frac pq \in\mathbb{Q}\cup\{\varnothing\}\backslash
\{-1, 0, \frac 12, 1, \frac 32, 2, 3\}$, $E(M_4(f))=\{-1, 0, 1, 2, \infty\}$ \phantom{\big|}\\
\begin{minipage}{1.0\textwidth}
\begin{equation*}
\text{ types}  
\begin{cases}
\{T \& A, T \& A, A, T \& A, D\}  \text{ if }\frac pq = \varnothing ; \\
\{Z, Z, Z, T, T^H \} \text{ if } \frac pq \in \mathbb{Z}; \\
\{T, T, T^H, T, T^H \} \text{ if } \frac pq=2+\frac1k; \\
\{T, T, Z, Z, T^H \} \text{ if } \frac pq=1+\frac1k; \\ 
\{T, T, Z, T, S^H\} \text{ if } \frac pq=\frac{n}{6n-1}; \\
\{T, T, Z, T, T^H\} \text{ otherwise.}
\end{cases}
\end{equation*}
\end{minipage} 
\\ \phantom{\big|} \\ \hline
\end{tabular}

\end{center}
\caption{\,\,\,\, Exceptional sets for $M_5(f)$ for $f$ in Table \ref{tb6}, part 2/2. \label{th}}
\end{table}

\end{document}